\documentclass[11pt]{amsart}
\usepackage{amsmath,amssymb,amsthm}
\usepackage[margin=1in]{geometry}
\usepackage{array}
\usepackage{longtable}
\usepackage[T1,OT1]{fontenc}
\newcommand{\Osek}{Os{\fontencoding{T1}\fontfamily{lmr}\selectfont\k{e}}kowski}
\usepackage[hidelinks]{hyperref}
\usepackage{orcidlink}

\newtheorem{theorem}{Theorem}[section]
\newtheorem{lemma}[theorem]{Lemma}
\newtheorem{proposition}[theorem]{Proposition}

\theoremstyle{definition}

\theoremstyle{remark}
\newtheorem{remark}[theorem]{Remark}

\numberwithin{equation}{section}

\newcommand{\E}{\mathbb{E}}
\newcommand{\Prob}{\mathbb{P}}
\newcommand{\N}{\mathbb{N}}
\newcommand{\Z}{\mathbb{Z}}
\newcommand{\F}{\mathcal{F}}

\AtBeginDocument{\hypersetup{pdfauthor={Saikat Kanjilal},pdftitle={Sharp constants and martingale proofs for two-weight Hardy inequalities with general measures}}}
\begin{document}

\title[Sharp constants for two-weight Hardy inequalities]{Sharp constants and martingale proofs for two-weight Hardy inequalities with general measures}

\author[Saikat Kanjilal]{Saikat Kanjilal\,\orcidlink{0000-0002-4359-8343}}
\address{Department of Mathematics, JIS University, Kolkata 700109, West Bengal, India}
\email{saikat.kanjilal.07@gmail.com}
\thanks{ORCID 0000-0002-4359-8343}


\begin{abstract}
We study the least constant $C_*$ in the two-weight Hardy inequality from $L^p(\sigma)$ to $L^q(\omega)$ for the operator $H_\sigma g(x)=\int_{(0,x]}g\,d\sigma$, where $\sigma$ is a Borel measure on the half-line that is finite on each interval $(0,x]$ and $\omega$ is a $\sigma$-finite Borel measure there, atoms allowed. Let $\mathcal A$ be the Maz'ya--Rozin functional and let $1<p<\infty$. For $0<q<p$ we find the best factor $\kappa_{p,q}$ in $C_*\le\kappa_{p,q}\mathcal A$. For $1\le q<p$ we prove $C_*\ge\mathcal B\ge\mathcal A$, where $\mathcal B$ is a Lebesgue--Stieltjes form of $\mathcal A$. The factor $1$ is best possible. For $0<q<1$ we give a lower factor valid for all such measures and an upper bound for the best lower factor. This bound comes from atoms of $\omega$. It is smaller than the lower factor proved by Sinnamon (1991) for weights. A two-atom example shows that the bound is not always sharp. The proofs use a reverse martingale with Doob's inequality, a special function of \Osek{} (2015), a quantile coupling and a splitting of atoms. The same tools recover the theorem of Li and Mao (2020) for $1<p\le q$, in the closed-interval form of Klaassen and Wellner (2021). They also recover the two-sided estimate of Sinnamon and Stepanov (1996) for $0<q<p$, which goes back to Maz'ya (1979) for $q\ge1$ and to Sinnamon (1991) for $q<1$. We also treat the case $p=1$ and a discrete criterion.
\end{abstract}

\maketitle
\enlargethispage{2pt}

\section{Introduction}

Let $1<p<\infty$ and $p'=p/(p-1)$. The classical inequality of Hardy \cite{Ha20,Ha25} states that
\begin{equation}\label{eq:hardy-classical}
\int_0^\infty \Big(\frac1x\int_0^x f(t)\,dt\Big)^p dx \le (p')^p \int_0^\infty f(x)^p\,dx
\end{equation}
for every nonnegative measurable function $f$. The constant $(p')^p$ is sharp. The history of \eqref{eq:hardy-classical} is described in \cite{KMP07,OK90}. The weighted form asks for conditions on weights $v$ and $w$ under which
\begin{equation}\label{eq:hardy-weighted}
\Big(\int_0^\infty \Big(\int_0^x f(t)\,dt\Big)^q w(x)\,dx\Big)^{1/q}\le C\Big(\int_0^\infty f(x)^p v(x)\,dx\Big)^{1/p}.
\end{equation}
Let us write $V(x)=\int_0^x v^{1-p'}$ and $W(x)=\int_x^\infty w$. For $1<p\le q<\infty$ the inequality \eqref{eq:hardy-weighted} holds if and only if
\begin{equation}\label{eq:B-weighted}
B=\sup_{y>0}W(y)^{1/q}V(y)^{1/p'}<\infty .
\end{equation}
This is due to Muckenhoupt \cite{Mu72} for $p=q$, and to Bradley \cite{Br78}, Kokilashvili \cite{Ko79} and Maz'ya \cite[Section 1.3.2]{Ma11}. For $0<q<p<\infty$ the criterion is
\begin{equation}\label{eq:A-weighted}
\mathcal A=\Big(\int_0^\infty W(x)^{r/p}V(x)^{r/p'}w(x)\,dx\Big)^{1/r}<\infty,\qquad \frac1r=\frac1q-\frac1p .
\end{equation}
It goes back to Rozin and Maz'ya \cite{Ma79}, see \cite[Section 1.3.3]{Ma11}, for $1\le q<p$, and to Sinnamon \cite{Si91} for $0<q<1$. For $0<q<1$ Sinnamon also proved that the best constant $C$ in \eqref{eq:hardy-weighted} satisfies $(p')^{1/p'}q^{1/p}\mathcal A\le C\le(p')^{1/r}q^{1/p}\mathcal A$, see \cite[Theorem 1(3)]{Si91}. Sinnamon and Stepanov \cite{SS96} gave an elementary proof of the whole range with explicit constants. They also solved the case $p=1$, $0<q<1$. For general measures the criteria go back to Sinnamon \cite{Si87}, see \cite[Corollaries 3.4 and 3.5]{Si05}. Scales of equivalent conditions were given by Okpoti, Persson and Sinnamon \cite{OPS07} for $1<p\le q$ and \cite{OPS08} for $0<q<p$. Li and Mao \cite{LM20} found the sharp factor for $1<p\le q$ and general Borel measures. Sinnamon \cite{Si22} introduced a normal form, which reduces weighted, discrete and general-measure Hardy operators to an operator $f\mapsto\int_0^{b(x)}f$ with a nonincreasing $b$. It recovers the factor of Li and Mao and improves the upper constants of \cite{SS96}. Gogatishvili and Pick \cite{GP25} gave a new elementary proof of the criteria for weights, with constants up to equivalence. Mynbaev and Lomakina \cite{ML25} studied the two-weight Hardy inequality on topological measure spaces. On time scales, Mahmoud, Saker, Anderson and Kenawy \cite{MSAK26} gave explicit two-sided bounds for $0<q<1<p$. Scales of equivalent conditions were given in \cite{GPTU25,PSW07}, with constants that are not claimed to be sharp. For Hardy type operators the monograph \cite[Remark 2.7.5, p.~118]{KPS17} notes that the question of sharp constants is largely open.

Klaassen and Wellner \cite{KW21} wrote the Hardy average as a conditional expectation. They formulated Hardy, Copson, Carleman and Muckenhoupt type inequalities in this language for $p\le q$. In particular they extended the result of \cite{LM20} so that it also contains the discrete inequality of Liao \cite{Li15}. Related results are in \cite{Kl23,LL23}. Walker \cite{Wa15} gave a probabilistic proof of the one-weight inequality for Lebesgue measure. The proofs in \cite{Kl23,KW21} are analytic. They rest on discrete Hardy lemmas, rearrangements and H\"older's inequality. A Doob martingale of right-tail averages appears in \cite[Section 8]{KW21}. However, no inequality is proved with it. The case $q<p$ is not treated there. Shargorodsky and Sharia \cite{ShSh22} found the largest norm of $I-\E[\,\cdot\mid\mathcal G]$ on $L^p$, $1<p<\infty$, over all probability spaces and all sub-$\sigma$-algebras $\mathcal G$. They used it to find the optimal constant in the bounded compact approximation property of $L^p(0,1)$. Their operator is the conditional centring map, not a Hardy average. Their estimates are not used in this paper. Bouthat and Paris\'e \cite{BP26} reduce weighted nested averages to a discrete Hardy inequality, which is a case of Proposition~\ref{prop:forms}(2). Theorem~\ref{thm:main} then improves their factor $p^{1/q}(p')^{1/p'}$ to the factor $k_{p,q}$ of \eqref{eq:kpq} when $p<q$. \Osek{} \cite{Os15} gave a special function proof of Bliss's inequality \cite{Bl30}. Ba\~nuelos and \Osek{} \cite{BO25} deduced from Bliss's inequality the existence of an abstract Bellman function of four variables. With it they proved a sharp Carleson embedding theorem on probability spaces with a tree structure, with the constant $(p-1)^{-1/q}k_{p,q}$. They also proved a sharp two-weight bound for fractional maximal operators under Sawyer's testing condition. Cavina \cite{Ca26} found an explicit Bellman function of four variables for the dual weighted Hardy inequality on a dyadic tree. His setting has one exponent $1<p<\infty$ on both sides and a testing condition on the weights. He proved that the constant $(p')^p$ in this inequality for $p$-th powers is sharp. In the arXiv preprint 2002.07532v3 of \cite{Ca26} he also identified his function with the value function of a stochastic optimal control problem, see Theorem 3.1 there. For $p=2$ a Bellman function proof of this inequality on the dyadic tree was given earlier by Arcozzi, Holmes, Mozolyako and Volberg \cite{AHMV21}. Arcozzi, Chalmoukis, Levi and Mozolyako \cite[Theorem 9]{ACLM23} adapted Cavina's proof to general rooted trees.

The main aim of this paper is to find sharp constants in the two-sided estimates of the least constant $C_*$ in \eqref{eq:main-ineq} below by the measure form $\mathcal A$ of \eqref{eq:A-weighted}, for $0<q<p$. We prove three results of this kind. First, for every $0<q<p$ we find the best upper factor $\kappa_{p,q}$ in $C_*\le\kappa_{p,q}\mathcal A$ (Theorem~\ref{thm:sharp}). To our knowledge, the factor $(r/q)^{1/r}$ is new for $0<q<1$. For $1<q<p$ the factor $(p')^{1/p'}q^{1/p}$ is known. Its optimality is new. Second, for $1\le q<p$ we prove $C_*\ge\mathcal B\ge\mathcal A$ (Theorem~\ref{thm:lowerB}), where $\mathcal B$ is a Lebesgue--Stieltjes form of the Maz'ya--Rozin functional. The factor $1$ is best possible. For $q=1$ it follows from duality. For $1<q<p$ it appears to be new. The lower factor $(p')^{1/p'}q^{1/p}(1-q/p)$ of \cite{SS96}, also given by \cite{Si22} for atomless $\omega$, is smaller. Third, for $0<q<1$ we prove $C_*\ge q^{1/p}\big((p-q)/(p-1)\big)^{1/p'}\mathcal A$ for general measures. This factor is larger than the factor $(p')^{1/p'}q^{1/p}(1-q/p)$ of \cite{SS96}. For weights Sinnamon \cite[Theorem 1(3)]{Si91} proved the larger lower factor $(p')^{1/p'}q^{1/p}$. Geometric chains of atoms of $\omega$ give an upper bound for the best lower factor. This bound is smaller than $1$. It is also smaller than Sinnamon's factor. So atoms of $\omega$ lower the best lower factor below its value for weights. A two-atom example shows that the upper bound is not always sharp (Theorem~\ref{thm:lowerq}). All these estimates hold for every Borel measure $\sigma$ that is finite on each interval $(0,x]$ and every $\sigma$-finite Borel measure $\omega$, atoms included. None of the factors depends on the measures.

A second aim of this paper is to prove the boundedness criteria themselves by probabilistic methods, for $1\le p<\infty$ and $0<q<\infty$ and for general measures. Let us describe the idea. Let $\sigma$ be finite and nonzero. Let $Y$ be a random variable with law proportional to $\sigma$. For infinite $\sigma$ one restricts $\sigma$ to $(0,R]$ and lets $R\to\infty$. Then $M_x=\E[g(Y)\mid Y\mathbf 1_{\{Y>x\}}]$ is a reverse martingale. It equals the Hardy average $A(x)=\E[g(Y)\mid Y\le x]$ on the event $\{Y\le x\}$. Doob's maximal inequality then gives the one-weight inequality with the constant $p'$. For general measures this inequality is due to Sinnamon \cite[Theorem 1.1]{Si05}, see also Klaassen and Wellner \cite[Theorem 2.1]{KW21}. Our proof also gives the maximal form. For $1<p\le q$ the reverse martingale alone, together with a Carleson type estimate, already gives the two-weight criterion $B<\infty$ with the factor $p^{1/q}p'$ (Proposition~\ref{prop:pleq-mart}). For the sharp factor in this case we evaluate the special function of \Osek{} \cite{Os15} along the reverse martingale and its companion for $g^p$, in the spirit of Burkholder's method. Atoms of $\sigma$ are treated directly. For $q<p$ the passage to two weights uses a quantile coupling. Let $Q$ be the quantile function of $\sigma$. Then $Q$ transports Lebesgue measure on $(0,S(\infty))$ to $\sigma$. If $S(\infty)<\infty$ and $T$ is uniformly distributed on $(0,S(\infty))$, then $Q(T)$ has law $\sigma/S(\infty)$. Conditional expectation given $Q(T)$ then makes the inequality \eqref{eq:main-ineq} below exactly equivalent to an inequality with Lebesgue measure in the input. If $S(\infty)=\infty$, the same reduction is a transport of $\sigma$-finite measures. No probability law is involved in this case. A splitting of atoms complements the coupling for $q<p$. The coupling is a probabilistic form of the rearrangement argument of \cite[Theorem 1.1]{Si05}, which uses $\chi_{(-\infty,x]}^*=\chi_{(0,\Lambda(x)]}$, and of the passage to normal form in \cite[Theorems 2.4 and 2.6]{Si22}. We do not claim it as new. Alhamzi, Sheikh, Jadhav, Beleyur and Jeelani \cite{ASJBJ26} used a transport of the same kind to carry sharp one-weight bounds to an atomless Borel measure with a continuous strictly increasing distribution function. Atoms are excluded there. Finally, the first passage times of $x\mapsto\sigma((0,x])$ through the levels $\lambda^k$ give a discrete criterion for $q<p$. This is a probabilistic version of the discretization of Gogatishvili and Pick \cite{GP03}. One direction of the comparison of the discrete functional with $\mathcal A$ is a direct tail estimate. The other follows from the two characterizations of the same constant.

With these tools we prove the theorem of Li and Mao with its sharp factor (Theorem~\ref{thm:main}), in the closed-interval form of Klaassen and Wellner \cite{KW21}. We also prove the two-sided estimate of Maz'ya, Rozin, Sinnamon and Stepanov for $0<q<p$ (Theorem~\ref{thm:qlp}) and the weighted case $p=1$ (Theorem~\ref{thm:p1}). Atoms are allowed in both measures. The range $q<1$ is covered without duality. These results are known, also for general measures, see \cite[Corollaries 3.4 and 3.5]{Si05}, \cite{OPS08,Pr06,PSU18,SH24,Si22}. For them our contribution is the method. Its reductions lose nothing in the constants. This is what makes the sharp constants of Theorems~\ref{thm:sharp}, \ref{thm:lowerB} and \ref{thm:lowerq} accessible. For $1<p<\infty$ Doob's inequality gives the maximal estimate behind the criteria and the upper bound of Theorem~\ref{thm:qlp}. The weighted case $p=1$ uses a separate level-set argument. The sharp upper factor for $q\le1$ comes from H\"older's inequality and Lemma~\ref{lem:tail}. The lower bounds come from test functions and a variational argument. Table~\ref{tab:constants} at the end of Section~\ref{sec:statement} compares all constants with the known ones and records which statements are new.

The paper is organized as follows. Section~\ref{sec:statement} states the results. Section~\ref{sec:tools} contains the probabilistic tools. Section~\ref{sec:pleq} treats the case $p\le q$. Section~\ref{sec:qlp} treats the case $q<p$, the sharp factors and the discrete criterion. Section~\ref{sec:endpoints} treats the endpoint cases. Section~\ref{sec:p1} treats the case $p=1$. Section~\ref{sec:corollaries} gives the weighted, discrete and dual forms. Section~\ref{sec:remarks} contains concluding remarks.

\section{Statement of the problem and results}\label{sec:statement}

Throughout the paper $\sigma$ and $\omega$ are Borel measures on $(0,\infty)$. We assume that
\begin{equation}\label{eq:locfin}
S(x)=\sigma((0,x])<\infty\qquad\text{for every } x>0 ,
\end{equation}
and we assume that $\omega$ is $\sigma$-finite. The function $S$ is nondecreasing and right-continuous. Since $\sigma((0,1])<\infty$, continuity from above gives $S(x)\to 0$ as $x\downarrow 0$. We write $S(\infty)=\lim_{x\to\infty}S(x)\in(0,\infty]$ when $\sigma\neq 0$. We also write
\[
W(x)=\omega([x,\infty)),\qquad x>0 .
\]

For a nonnegative Borel function $g$ on $(0,\infty)$ we define the Hardy operator with respect to $\sigma$ by
\begin{equation}\label{eq:Hsigma}
H_\sigma g(x)=\int_{(0,x]} g\,d\sigma ,\qquad x>0 .
\end{equation}
Let $1<p<\infty$ and $0<q<\infty$. We study the inequality
\begin{equation}\label{eq:main-ineq}
\Big(\int_{(0,\infty)} \big(H_\sigma g\big)^q\,d\omega\Big)^{1/q}\le C\Big(\int_{(0,\infty)} g^p\,d\sigma\Big)^{1/p}
\end{equation}
for all nonnegative Borel functions $g$. We use the convention $0\cdot\infty=0$. Let $C_*$ be the infimum of the constants $C\ge0$ for which \eqref{eq:main-ineq} holds, with $C_*=\infty$ if no finite constant exists. The infimum is attained when it is finite. If $0<C_*<\infty$, this follows by letting $C\downarrow C_*$ for each $g$ with finite right side. The case of an infinite right side is trivial. If $C_*=0$, then the test functions $\mathbf 1_{(0,n]}$ give $\omega(\{S>0\})=0$. So $H_\sigma g=0$ $\omega$-almost everywhere for every $g$. The constant $C=0$ is then admissible. We call $C_*$ the least constant and set
\begin{equation}\label{eq:B}
B=\sup_{y>0}\ W(y)^{1/q}\,S(y)^{1/p'} ,
\end{equation}
and, for $q<p$ with $1/r=1/q-1/p$,
\begin{equation}\label{eq:A}
\mathcal A=\Big(\int_{(0,\infty)}W(x)^{r/p}\,S(x)^{r/p'}\,d\omega(x)\Big)^{1/r}.
\end{equation}
For $0<q<\infty$ we write $q'$ for the algebraic conjugate exponent, defined by $1/q+1/q'=1$. In particular $q'=q/(q-1)<0$ if $0<q<1$. If $q=1$, then $q'=\infty$ with $1/q'=0$. We use $q'$ only in identities between exponents, such as $r/q'=r/p'-1$. No H\"older duality in $L^q$ with $q<1$ is intended. If $\sigma=0$, then $H_\sigma g=0$ and $C_*=B=\mathcal A=0$. All statements below are then trivial. We assume $\sigma\neq0$ whenever this is convenient.

\subsection{The classical criteria}
The first part of the paper proves the classical characterizations by the reverse martingale. We state them in the measure form \eqref{eq:main-ineq}, with atoms allowed in both measures.

\begin{theorem}\label{thm:main}
Let $1<p\le q<\infty$. Then the inequality \eqref{eq:main-ineq} holds with a finite constant if and only if $B<\infty$. The least constant $C_*$ satisfies
\begin{equation}\label{eq:constants}
B\le C_*\le k_{p,q}\,B ,
\end{equation}
where $k_{p,p}=p^{1/p}(p')^{1/p'}$ and, for $p<q$ and $\varrho=q/p-1$,
\begin{equation}\label{eq:kpq}
k_{p,q}=\Big(\frac{\varrho}{\mathrm B(1/\varrho,(q-1)/\varrho)}\Big)^{1/p-1/q},
\end{equation}
with $\mathrm B$ the Beta function. The factor $k_{p,q}$ is best possible.
\end{theorem}

The sharp factor $k_{p,q}$ for general measures is due to Li and Mao \cite[Theorem 1]{LM20}, who state the inequality with the open interval $(-\infty,x)$ in the inner integral. Their form does not contain the closed interval form \eqref{eq:main-ineq} when $\sigma$ has atoms, see \cite[Remark 1]{LM20}. Klaassen and Wellner \cite[Theorem 3.1]{KW21} proved the upper bound $C_*\le k_{p,q}B$ with the closed interval, for arbitrary $\sigma$-finite measures. The lower bound $B\le C_*$ is elementary, see \cite[Remark 3.3]{KW21}. The optimality of $k_{p,q}$ follows from Bliss's inequality \cite{Bl30}. Maz'ya \cite[Section 1.3.2, Theorem 1]{Ma11} stated it in measure form, for an arbitrary measure $\omega$ and an absolutely continuous $\sigma$, with $B\le C_*\le q^{1/q}(q')^{1/p'}B$. The factor $k_{p,q}$ is smaller than this factor and than Opic's factor $(1+q/p')^{1/q}(1+p'/q)^{1/p'}$ from \cite{OK90} when $p<q$. For example, at $p=2$ and $q=3$ the three factors are about $1.468$, $1.766$ and $1.752$. Sinnamon \cite[Theorems 2.6 and 3.8]{Si22} gave a short proof through his normal form. Section~\ref{sec:pleq} proves the criterion through Doob's inequality and the sharp factor through a special function along the reverse martingale. It also records a proof through the quantile coupling as an alternative.

\begin{theorem}\label{thm:qlp}
Let $1<p<\infty$ and $0<q<p$. Then the inequality \eqref{eq:main-ineq} holds with a finite constant if and only if $\mathcal A<\infty$. The least constant $C_*$ satisfies
\begin{equation}\label{eq:constants-qlp}
(p')^{1/p'}q^{1/p}\Big(1-\frac qp\Big)\mathcal A\le C_*\le\Big(\frac rq\Big)^{1/r}p^{1/p}(p')^{1/p'}\mathcal A .
\end{equation}
\end{theorem}

Theorem~\ref{thm:qlp} goes back to Rozin and Maz'ya \cite{Ma79}, Sinnamon \cite{Si91} and Sinnamon and Stepanov \cite[Theorem 2.4]{SS96}. For $1\le q<p$, Maz'ya \cite[Section 1.3.3, Theorem 2]{Ma11} stated it in measure form, for an arbitrary measure $\omega$ and an absolutely continuous $\sigma$, with explicit constants (see the discussion after the proof of Theorem~\ref{thm:lowerB}). Sinnamon and Stepanov stated it for weights. For atomless $\omega$ its transfer to measures $\sigma$ with atoms or singular parts also follows from the normal form of \cite{Si22}. Section~\ref{sec:SS} proves it for all $\sigma$ and $\omega$, with Lemma~\ref{lem:doob} in the place of Hardy's inequality. Its upper bound is improved in Theorem~\ref{thm:sharp}. For general measures its lower bound is improved for every $0<q<p$ in Theorem~\ref{thm:lowerq}(1). It is improved further for $1\le q<p$ in Theorem~\ref{thm:lowerB}. For weights and $0<q<1$ a larger lower factor is due to Sinnamon \cite[Theorem 1(3)]{Si91}, see the discussion after Theorem~\ref{thm:lowerq}.

For $p=1$ the measure form \eqref{eq:main-ineq} is degenerate, see Proposition~\ref{prop:endpoints}(3). The natural weighted problem uses two measures in the input. Let $\mu$ and $\omega$ be $\sigma$-finite Borel measures on $(0,\infty)$ and let $v:(0,\infty)\to[0,\infty]$ be Borel. We consider
\begin{equation}\label{eq:p1-ineq}
\Big(\int_{(0,\infty)}\Big(\int_{(0,x]}f\,d\mu\Big)^q d\omega(x)\Big)^{1/q}\le C\int_{(0,\infty)}f\,v\,d\mu
\end{equation}
for all nonnegative Borel functions $f$. We denote by $C_{1,*}$ the infimum of the constants $C\ge0$ for which \eqref{eq:p1-ineq} holds, with $C_{1,*}=\infty$ if there is none. As for $C_*$, the infimum is attained when $0<C_{1,*}<\infty$. It need not be attained when $C_{1,*}=0$, because $v$ may take the value $\infty$. For example, let $\mu=\delta_1$, $\omega=\delta_2$ and $v(1)=\infty$. Then every $C>0$ is admissible, since a function $f$ with $f(1)>0$ has $\int fv\,d\mu=\infty$. However, $f=\mathbf 1_{\{1\}}$ shows that $C=0$ is not. All statements about $C_{1,*}$ below refer to this infimum. Put $u=1/v$, with $1/0=\infty$ and $1/\infty=0$, and let
\begin{equation}\label{eq:U}
U(x)=\mu\text{-}\operatorname*{ess\,sup}_{(0,x]}u ,
\end{equation}
with $U(x)=0$ when $\mu((0,x])=0$. The function $U$ is nondecreasing with values in $[0,\infty]$. For Lebesgue measure $\mu$ we have $U(x)=\big(\operatorname{ess\,inf}_{(0,x)}v\big)^{-1}$.

\begin{theorem}\label{thm:p1}
Let $\mu$, $\omega$, $v$ and $U$ be as above. We use the convention $0\cdot\infty=0$.
\begin{enumerate}
\item If $1\le q<\infty$, then $C_{1,*}=B_1$, where $B_1=\sup_{y>0}U(y)\,W(y)^{1/q}$.
\item If $q=\infty$, with the $\omega$-essential supremum on the left of \eqref{eq:p1-ineq}, then $C_{1,*}=\sup\{U(y):W(y)>0\}$, with $\sup\emptyset=0$.
\item If $0<q<1$, put $r=q/(1-q)$ and
\[
\mathcal A_1=\Big(\int_{(0,\infty)}W(x)^{r}\,U(x)^{r}\,d\omega(x)\Big)^{1/r}.
\]
Then \eqref{eq:p1-ineq} holds with a finite constant if and only if $\mathcal A_1<\infty$, and
\begin{equation}\label{eq:constants-p1}
q(1-q)\,\mathcal A_1\le C_{1,*}\le(1-q)^{-(1-q)/q}\mathcal A_1 .
\end{equation}
\end{enumerate}
\end{theorem}

When $\mu$ is Lebesgue measure, part (1) is the case $p=1$ of \cite[Section 1.3.2, Theorems 1 and 2]{Ma11}, where $\omega$ may be an arbitrary measure. For weights, part (3) is due to Sinnamon and Stepanov \cite{SS96}. For general measures, part (1) and a characterization equivalent to part (3) are due to Santacruz Hidalgo \cite{SH24}. His results assume $\sigma$-finite measures and core maps with $\mu((0,x])<\infty$ for every $x$, see \cite[Definition 1 and Theorem A]{SH24}. The general case reduces to this one. Let $\lambda$ be the restriction to $\{v<\infty\}$ of a finite measure equivalent to $\mu$. A function $f$ with $\int fv\,d\mu<\infty$ vanishes $\mu$-almost everywhere on $\{v=\infty\}$, where $u=0$. It follows that replacing $\mu$ by $\lambda$ and $f$ by $f\,d\mu/d\lambda$ changes neither $C_{1,*}$ nor $U$. For three $\sigma$-finite Borel measures on the half-line, a characterization equivalent to part (3) was obtained earlier by Johansson, Stepanov and Ushakova \cite[Theorem 3.1]{JSU08}. Section~\ref{sec:p1} gives a proof through the level sets $\{U^r>s\}$ of the running essential supremum $U$.

\begin{remark}
Let $p>1$. Under the assumptions of Proposition~\ref{prop:forms}(1) the form \eqref{eq:main-ineq} contains the weighted form \eqref{eq:hardy-weighted} through $d\sigma=v^{1-p'}dx$, $d\omega=w\,dx$ and $f=g\,v^{1-p'}$. Then $B$ and $\mathcal A$ become \eqref{eq:B-weighted} and \eqref{eq:A-weighted}, see Proposition~\ref{prop:forms}(1).
\end{remark}

\subsection{Sharp constants for \texorpdfstring{$q<p$}{q < p}}
The reductions of Section~\ref{sec:tools} lose nothing in the constants. For $q<p$ this allows us to determine the best upper factor in Theorem~\ref{thm:qlp}, and the best lower factor when $1\le q<p$.

\begin{theorem}\label{thm:sharp}
Let $1<p<\infty$, $0<q<p$ and $1/r=1/q-1/p$. Put
\begin{equation}\label{eq:kappa}
\kappa_{p,q}=\Big(\frac rq\Big)^{1/r}\quad(0<q\le1),\qquad\kappa_{p,q}=(p')^{1/p'}q^{1/p}\quad(1\le q<p).
\end{equation}
\begin{enumerate}
\item For all $\sigma$ and $\omega$ as above,
\begin{equation}\label{eq:sharp}
C_*\le\kappa_{p,q}\,\mathcal A .
\end{equation}
\item The factor $\kappa_{p,q}$ is best possible for every $0<q<p$, already when $\sigma$ is Lebesgue measure.
\item If $q=1$, then $\mathcal A\le C_*\le(p')^{1/p'}\mathcal A$. Both factors are best possible.
\end{enumerate}
\end{theorem}

The two expressions in \eqref{eq:kappa} agree at $q=1$. Theorem~\ref{thm:sharp} improves the upper bound in Theorem~\ref{thm:qlp} for every $0<q<p$. For $1\le q<p$ the lower bound $\mathcal A\le C_*$, with the best factor $1$, is part of Theorem~\ref{thm:lowerB} below. For $1\le q<p$ the bound \eqref{eq:sharp} also follows from the normal form estimate $N_{p,q}(b)\le C_0$ of Sinnamon \cite[Corollary 3.21]{Si22}. Here $N_{p,q}(b)$ is the norm of the normal form operator with parameter $b$. The quantity $C_0$ is Sinnamon's functional from \cite[Theorem 3.9]{Si22}. Let $b$ be the parameter of \eqref{eq:main-ineq} given by \cite[Theorem 2.4]{Si22}. It is the nonincreasing rearrangement of $S$ with respect to $\omega$, which requires only that $\omega$ is $\sigma$-finite. Thus the Lebesgue measure of $\{b>t\}$ equals $\omega(\{S>t\})$ for every $t>0$. Put $E_t=\{S^{r/p'}>t\}$. Since $b$ is nonincreasing, $\{b^{r/p'}>t\}$ is an interval with left endpoint $0$ and length $\omega(E_t)$. The layer cake formula and Fubini's theorem give
\[
C_0^r=(p')^{r/p'}q^{r/p}\,\frac qr\int_0^\infty\omega(E_t)^{r/q}\,dt,\qquad\mathcal A^r=\int_0^\infty\int_{E_t}W^{r/p}\,d\omega\,dt .
\]
Each $E_t$ is a terminal interval, so that $W(x)=\omega(E_t\cap[x,\infty))$ for $x\in E_t$. Let $\nu\le\omega|_{E_t}$ be a finite measure. Lemma~\ref{lem:tail} applied to $\nu$ on $E_t$ with $\gamma=r/p$ gives
\[
(q/r)\,\nu(E_t)^{r/q}\le\int_{E_t}\nu(E_t\cap[x,\infty))^{r/p}\,d\nu(x)\le\int_{E_t}W^{r/p}\,d\omega.
\]
Since $\omega$ is $\sigma$-finite, the supremum over such $\nu$ gives $(q/r)\,\omega(E_t)^{r/q}\le\int_{E_t}W^{r/p}\,d\omega$. So $C_0\le(p')^{1/p'}q^{1/p}\mathcal A$ for all $\sigma$, $\omega$ and $0<q<p$. If $\omega$ has no atoms, then equality holds, with both sides possibly infinite. If $\omega=\ell\delta_T$ with $\ell>0$ and $S(T)>0$, then $C_0=(p')^{1/p'}q^{1/p}(q/r)^{1/r}\mathcal A$. In particular, for $1\le q<p$ Sinnamon's estimate implies \eqref{eq:sharp} for every $\omega$. For some $\omega$ with atoms it is better than \eqref{eq:sharp}. For weights and $1<q<p$ the same upper estimate follows from \cite[Theorem 1(2)]{Si91}. For absolutely continuous $\sigma$ and $1\le q<p$ the bound \eqref{eq:sharp} also follows from the upper estimate of Maz'ya \cite[Section 1.3.3, Lemma and Theorem 2]{Ma11} together with the integration by parts \eqref{eq:ibp}. For $1\le q<p$ our proof follows Maz'ya's proof, with Lemma~\ref{lem:doob} in place of Hardy's inequality. For $q\le1$ it uses H\"older's inequality and Lemma~\ref{lem:tail} instead. For $1\le q<p$ the known lower bound is $C_*\ge(p')^{1/p'}q^{1/p}(1-q/p)\mathcal A$, see \cite[Theorem 2.4]{SS96}, and \cite[Theorem 3.9]{Si22} for atomless $\omega$. This factor is at most $1$. Section~\ref{sec:SS} proves this lower bound for all $\sigma$ and $\omega$. For $\omega$ with atoms Sinnamon's bound $(q/r)\,C_0\le C_*$ can give less, for example for a point mass. Maz'ya's two-sided estimate is stated for the functional $\mathcal B_\sigma=\big(\int W^{r/q}S^{r/q'}\,d\sigma\big)^{1/r}$, see \cite[Section 1.3.3, Lemma on p.~48 and Theorems 1 and 2 on pp.~50--51]{Ma11}. His lower bound is compared with Theorem~\ref{thm:lowerB} after its proof. His upper factor $q^{1/q}(p')^{1/q'}$ for $\mathcal B_\sigma$ is not stated there to be best possible. We show after the proof of Theorem~\ref{thm:lowerB} that it is best possible. The proof of Theorem~\ref{thm:sharp} is in Section~\ref{sec:sharp}. Remark~\ref{rem:sharp-compare} compares the result with the normal form bounds.

For $1\le q<p$ we put $\alpha=r/p'$ and $\beta=r/q$, and define
\begin{equation}\label{eq:Bcal}
\mathcal B=\Big(\int_{(0,\infty)}W(x)^{\beta}\,d\big(S^{\alpha}\big)(x)\Big)^{1/r},
\end{equation}
where $d(S^\alpha)$ is the Lebesgue--Stieltjes measure of the nondecreasing right-continuous function $S^\alpha$. This is the second form of the Maz'ya--Rozin functional. For weights, with $d\sigma=v^{1-p'}dx$, so that $S=V$ with $V$ as in the introduction, put $V_{\mathrm{GP}}=V^{1/p'}$. Then $\mathcal B^r$ is identical to the Lebesgue--Stieltjes functional $\int_0^\infty W^{p/(p-q)}\,d\big(V_{\mathrm{GP}}^{pq/(p-q)}\big)$, because $r/q=p/(p-q)$ and $r/p'=(pq/(p-q))/p'$. Gogatishvili and Pick \cite[Theorem 2.1]{GP25} use this functional as their criterion for every $0<q<p$, with constants up to equivalence. So the functional is not new. What Theorem~\ref{thm:lowerB} adds is the factor $1$ in $\mathcal B\le C_*$. If $\omega=\sum_{i=1}^n w_i\delta_{x_i}$ with $0<x_1<\dots<x_n$, put $x_0=0$, $S_i=S(x_i)$, $S_0=0$, $W_i=\sum_{j\ge i}w_j$ and $W_{n+1}=0$. Then $W=W_i$ on $(x_{i-1},x_i]$ and $W=0$ on $(x_n,\infty)$, so
\begin{equation}\label{eq:Bcal-discrete}
\mathcal B^r=\sum_{i=1}^nW_i^{\beta}\big(S_i^{\alpha}-S_{i-1}^{\alpha}\big)=\sum_{i=1}^nS_i^{\alpha}\big(W_i^{\beta}-W_{i+1}^{\beta}\big).
\end{equation}

\begin{theorem}\label{thm:lowerB}
Let $1<p<\infty$ and $1\le q<p$. Then $\mathcal A\le\mathcal B\le C_*$. If $\omega$ is a point mass, then $\mathcal A=\mathcal B=C_*$.
\end{theorem}

The proof is in Section~\ref{sec:lowerB}. For atomless $\sigma$ the functional $\mathcal B$ is a constant multiple of the form of the Maz'ya--Rozin condition used in \cite{Ma11}. Theorem~\ref{thm:lowerB} improves the classical lower factor, see the discussion after the proof.

For $0<q<p$ let $\lambda_{p,q}$ be the infimum of $C_*/\mathcal A$ over all pairs $(\sigma,\omega)$ with $0<\mathcal A<\infty$. Theorem~\ref{thm:lowerB} gives $\lambda_{p,q}=1$ for $1\le q<p$. For $0<q<p$ we put
\begin{equation}\label{eq:lambdas}
\lambda^{\flat}_{p,q}=q^{1/p}\Big(\frac{p-q}{p-1}\Big)^{1/p'},\qquad R_{p,q}(t)=\frac{(1-t^q)^{1/p}(1-t^{p'})^{1/p'}}{1-t},\qquad\lambda^{\sharp}_{p,q}=\inf_{0<t<1}R_{p,q}(t).
\end{equation}

\begin{theorem}\label{thm:lowerq}
Let $1<p<\infty$.
\begin{enumerate}
\item If $0<q<p$, then $C_*\ge\lambda^\flat_{p,q}\mathcal A$ for all $\sigma$ and $\omega$.
\item If $0<q<1$, then $\lambda^\flat_{p,q}\le\lambda_{p,q}\le\lambda^\sharp_{p,q}<1$.
\item $\lambda_{3,4/5}<\lambda^\sharp_{3,4/5}$.
\end{enumerate}
\end{theorem}

The factor $\lambda^\flat_{p,q}=(p')^{1/p'}q^{1/p}(1-q/p)^{1/p'}$ exceeds the lower factor in \eqref{eq:constants-qlp} by the factor $(1-q/p)^{-1/p}$. It equals $1$ at $q=1$. For weights a larger lower factor is known. Let $0<q<1$ and let $\sigma$ and $\omega$ be absolutely continuous. Then Sinnamon \cite[Theorem 1(3)]{Si91} proved $(p')^{1/p'}q^{1/p}\mathcal A\le C_*\le(p')^{1/r}q^{1/p}\mathcal A$. His functional $D$ is $\mathcal A$ in the notation \eqref{eq:A-weighted}. Part (1) holds for every Borel measure $\sigma$ that is finite on each interval $(0,x]$ and every $\sigma$-finite Borel measure $\omega$. Its factor is smaller than Sinnamon's by the factor $(1-q/p)^{1/p'}$. The upper bound in (2) comes from geometric chains of atoms of $\omega$ with $\sigma=\mathrm m$, that is, the configurations \eqref{eq:ab-config} with $a_i=t^{-i}$ and $b_i=t^i$. We have $R_{p,q}(t)\to(p')^{1/p'}q^{1/p}$ as $t\uparrow1$. The proof of (2) shows that $R_{p,q}(t)<(p')^{1/p'}q^{1/p}$ for $t$ close to $1$. Hence $\lambda^\sharp_{p,q}<(p')^{1/p'}q^{1/p}$. For $\sigma=\mathrm m$ and absolutely continuous $\omega$ the ratio $C_*/\mathcal A$ is at least $(p')^{1/p'}q^{1/p}$. For atomic $\omega$ it comes arbitrarily close to $\lambda^\sharp_{p,q}$. In this sense atoms of $\omega$ lower the best lower factor. For $1\le q<p$ the same happens with the factors $(r/q)^{1/r}$ and $1$, see the discussion after the proof of Theorem~\ref{thm:lowerB}. Part (3) shows that the bound in (2) is not always sharp. The proofs are in Section~\ref{sec:lowerq}.

Section~\ref{sec:quantile} gives a discrete characterization. For $\lambda>1$ let $(x_k)_{k\in K}$ be the quantile levels of $\sigma$ defined in \eqref{eq:levels}, with the index set $K$ and the convention $x_{k^*+1}=\infty$ given there. Let $I_k=[x_k,x_{k+1})$. Put
\begin{equation}\label{eq:D}
\mathcal D_\lambda=\Big(\sum_{k\in K}\lambda^{kr/p'}\,\omega(I_k)^{r/q}\Big)^{1/r}.
\end{equation}

\begin{theorem}\label{thm:Dlambda}
Let $1<p<\infty$, $0<q<p$ and $\lambda>1$. Then \eqref{eq:main-ineq} holds with a finite constant if and only if $\mathcal D_\lambda<\infty$. The least constant $C_*$ satisfies
\[
2^{-1/r}\Big(1-\frac1\lambda\Big)^{1/p'}\mathcal D_\lambda\le C_*\le\lambda\Big(1-\frac1\lambda\Big)^{-1/p}p'\,\mathcal D_\lambda ,
\]
and
\[
\Big(\frac qr\Big)^{1/r}\mathcal D_\lambda\le\mathcal A\le\frac{\lambda(1-1/\lambda)^{-1/p}p'}{c_{p,q}}\,\mathcal D_\lambda ,
\]
where $c_{p,q}=\lambda^\flat_{p,q}$ for $0<q<1$ and $c_{p,q}=1$ for $1\le q<p$.
\end{theorem}

The left inequality in the second display is a direct tail estimate. The right one follows from the two characterizations of the same constant $C_*$, through Theorems~\ref{thm:lowerB} and \ref{thm:lowerq}(1). Remark~\ref{rem:Dlambda-constants} optimizes the constants in $\lambda$.

\begingroup
\small
\renewcommand{\arraystretch}{1.5}
\setlength{\LTcapwidth}{\textwidth}\setlength{\LTpre}{\bigskipamount}\setlength{\LTpost}{\bigskipamount}
\begin{longtable}{>{\raggedright\arraybackslash}p{2.1cm}p{1.6cm}>{\raggedright\arraybackslash}p{5.0cm}>{\raggedright\arraybackslash}p{4.7cm}}
\caption{Two-sided estimates for the best constant. The first factor in the first row is Opic's factor from \cite{OK90}, and $k_{p,q}$ is defined in \eqref{eq:kpq}. The functional $\mathcal B$ is defined in \eqref{eq:Bcal}. Maz'ya's functional $\mathcal B_\sigma$ is defined after Theorem~\ref{thm:sharp}. The factors $\lambda^\flat_{p,q}$ and $\lambda^\sharp_{p,q}$ are defined in \eqref{eq:lambdas}. The quantity $B_\infty$ is defined in Proposition~\ref{prop:endpoints}(1). The estimates of \cite{KW21}, \cite{LM20}, \cite{Si22} and \cite{SH24} hold for general measures. Those of \cite{Ma11} hold for an arbitrary $\omega$ and an absolutely continuous $\sigma$ (Lebesgue $\mu$ when $p=1$). The other known estimates are stated for weights on the half-line. The estimates of this paper hold for general measures. In the last column, ``new'' marks a statement that is new to our knowledge.}\label{tab:constants}\\
\hline
Regime & Criterion & Known estimates & This paper \\
\hline
\endfirsthead
\multicolumn{4}{l}{\tablename~\thetable\ (continued)}\\
\hline
Regime & Criterion & Known estimates & This paper \\
\hline
\endhead
\multicolumn{4}{r}{\emph{continued on next page}}\\
\endfoot
\endlastfoot
$1<p\le q<\infty$ & $B$ & $B\le C_*\le(1+\frac q{p'})^{\frac1q}(1+\frac{p'}q)^{\frac1{p'}}B$ \cite{OK90}. $B\le C_*\le q^{\frac1q}(q')^{\frac1{p'}}B$ \cite{Ma11}. Best factor $k_{p,q}$ for general measures \cite{LM20}, in the closed-interval form of \cite{KW21}, normal form proof \cite{Si22}. & Known factor, new proof. $B\le C_*\le k_{p,q}B$, proof by a special function along the reverse martingale, alternative proof by the coupling. $C_*\le p^{1/q}p'B$ by Doob's inequality alone \\
\hline
$1<p<\infty$, $0<q<p$ & $\mathcal A$ & $(p')^{\frac1{p'}}q^{\frac1p}(1-\frac qp)\mathcal A\le C_*\le(\frac rq)^{\frac1r}p^{\frac1p}(p')^{\frac1{p'}}\mathcal A$ \cite{SS96}. Normal form, all abstract Hardy operators. Same lower bound for atomless $\omega$, better upper bound for $q\ge1$ \cite{Si22}. & Martingale proof of \cite{SS96}. Sharp upper factor $\kappa_{p,q}$ (Thm.~\ref{thm:sharp}). Sharp lower factor $1$ for $q\ge1$ (Thm.~\ref{thm:lowerB}). Discrete criterion with explicit constants (Thm.~\ref{thm:Dlambda}) \\
\hline
$1<p<\infty$, $0<q\le1$ & $\mathcal A$ & For weights and $q<1$, $(p')^{\frac1{p'}}q^{\frac1p}\mathcal A\le C_*\le(p')^{\frac1r}q^{\frac1p}\mathcal A$ \cite{Si91}. Normal form bound $N_{p,q}(b)\le\max(1,(p')^{\frac{1-q}q})C_0$ \cite{Si22}, the same upper bound for weights (Remark~\ref{rem:sharp-compare}) & New for $0<q<1$, and classical for $q=1$. $C_*\le(\frac rq)^{\frac1r}\mathcal A$, best possible. $C_*\ge\lambda^\flat_{p,q}\mathcal A$ for general measures. Best lower factor $\le\lambda^\sharp_{p,q}<(p')^{\frac1{p'}}q^{\frac1p}$ for $q<1$, by atoms of $\omega$. Bound not always sharp (Thm.~\ref{thm:lowerq}) \\
\hline
$1<p<\infty$, $1\le q<p$ & $\mathcal A$ & $q(\frac{p'}r)^{\frac1{q'}}\mathcal B_\sigma\le C_*\le q^{\frac1q}(p')^{\frac1{q'}}\mathcal B_\sigma$, $C_*=\mathcal B_\sigma$ at $q=1$ \cite{Ma11}. For weights and $q>1$ the same two-sided estimate is in \cite[Theorem 1(2)]{Si91}. Normal form bound $N_{p,q}(b)\le C_0$ \cite{Si22}, which gives $C_*\le(p')^{\frac1{p'}}q^{\frac1p}\mathcal A$ (see after Thm.~\ref{thm:sharp}) & Optimality of the upper factors new for $q>1$. $C_*\ge\mathcal B\ge\mathcal A$, best possible, new for $q>1$ (Thm.~\ref{thm:lowerB}) \\
\hline
$1\le p\le\infty$, $q=\infty$ & $B_\infty$ & $C_*=B_\infty$ for $1<p<\infty$ \cite{Ma11} & $C_*=B_\infty$ for $1\le p\le\infty$ (Prop.~\ref{prop:endpoints}) \\
\hline
$p=1$, $1\le q<\infty$ & $B_1$ & $C_{1,*}=B_1$ \cite{Ma11}. General measures \cite{SH24} & $C_{1,*}=B_1$, new proof \\
\hline
$p=1$, $0<q<1$ & $\mathcal A_1$ & $q(1-q)\mathcal A_1\le C_{1,*}\le(1-q)^{-\frac{1-q}q}\mathcal A_1$ \cite{SS96}. General measures up to equivalence \cite{JSU08,SH24} & The same constants for general measures, new. Direct proof (Thm.~\ref{thm:p1}) \\
\hline
\end{longtable}
\endgroup

\section{Probabilistic tools}\label{sec:tools}

In this section we collect the probabilistic tools. The first is a reverse martingale whose maximal function controls Hardy averages. The second is a quantile coupling that compares an arbitrary input measure with Lebesgue measure. It is a probabilistic version of the normal form reduction of Sinnamon \cite{Si22}. The third is an elementary lemma on tails of measures.

\subsection{A reverse martingale and the one-weight inequality}\label{sec:doob}

In this subsection $\sigma$ is a finite nonzero Borel measure on $(0,\infty)$. Let $(\Omega,\F,\Prob)$ be a probability space that carries a random variable $Y$ with law $\Prob(Y\in E)=\sigma(E)/S(\infty)$. Let $g\ge 0$ be a Borel function with $\int g^p\,d\sigma<\infty$. Since $\sigma$ is finite, H\"older's inequality gives $g\in L^1(\sigma)$. In particular $\E\, g(Y)<\infty$.

For $x>0$ with $S(x)>0$ we define the Hardy average
\begin{equation}\label{eq:Ax}
A(x)=\frac{H_\sigma g(x)}{S(x)}=\frac{1}{S(x)}\int_{(0,x]}g\,d\sigma=\E\big[g(Y)\mid Y\le x\big] .
\end{equation}
For every $y>0$ we define the maximal average
\begin{equation}\label{eq:Astar}
A^*(y)=\sup\{A(x)\,:\, x\ge y,\ S(x)>0\} ,
\end{equation}
with $\sup\emptyset=0$. The function $A^*$ is nonincreasing and hence Borel measurable. For $y$ with $S(y)>0$ we have $A^*(y)\le \|g\|_{L^1(\sigma)}/S(y)<\infty$.

For $x>0$ let $\F_x$ be the $\sigma$-algebra generated by the random variable $Y\mathbf 1_{\{Y>x\}}$. On the event $\{Y\le x\}$ this random variable is zero. On the event $\{Y>x\}$ it reveals $Y$.

\begin{lemma}\label{lem:filtration}
The following statements hold.
\begin{enumerate}
\item If $0<x'<x$, then $\F_x\subset\F_{x'}$.
\item If $S(x)>0$, then almost surely
\begin{equation}\label{eq:condexp}
M_x:=\E\big[g(Y)\mid\F_x\big]=g(Y)\mathbf 1_{\{Y>x\}}+A(x)\mathbf 1_{\{Y\le x\}} .
\end{equation}
\end{enumerate}
\end{lemma}

\begin{proof}
(1) Let $Z'=Y\mathbf 1_{\{Y>x'\}}$. Since $x>x'$, we have $Y\mathbf 1_{\{Y>x\}}=Z'\mathbf 1_{\{Z'>x\}}$. This Borel function of $Z'$ is $\F_{x'}$-measurable.

(2) The right-hand side of \eqref{eq:condexp} is a Borel function of $Y\mathbf 1_{\{Y>x\}}$, because $\{Y\le x\}=\{Y\mathbf 1_{\{Y>x\}}=0\}$ and $Y>x>0$ on the complement. It is integrable. Every set in $\F_x$ has the form $G=\{Y\le x\}\cup\{Y\in E\}$ or $G=\{Y\in E\}$ with a Borel set $E\subset(x,\infty)$. On $\{Y\in E\}$ both sides of \eqref{eq:condexp} are equal to $g(Y)$. On $\{Y\le x\}$ we have
\[
\E\big[g(Y)\mathbf 1_{\{Y\le x\}}\big]=\frac{1}{S(\infty)}\int_{(0,x]}g\,d\sigma=A(x)\,\Prob(Y\le x).
\]
Hence the two sides have the same integral over every $G\in\F_x$.
\end{proof}

By Lemma~\ref{lem:filtration}, the family $(M_x)$ is a martingale when $x$ decreases. More precisely, let $x_1>x_2>\dots>x_n$ be points with $S(x_n)>0$. Then $(M_{x_i})_{i=1}^n$ is a martingale with respect to the increasing filtration $(\F_{x_i})_{i=1}^n$.

\begin{lemma}\label{lem:doob}
Let $\sigma$ be a finite nonzero Borel measure on $(0,\infty)$ and let $1<p<\infty$. Then for every nonnegative Borel function $g$ with $\int g^p\,d\sigma<\infty$ (the inequality is trivial otherwise) we have
\begin{equation}\label{eq:one-weight}
\int_{(0,\infty)} A^*(y)^p\,d\sigma(y)\le (p')^p\int_{(0,\infty)} g^p\,d\sigma .
\end{equation}
In particular
\begin{equation}\label{eq:hardy-sigma}
\int_{(0,\infty)} \Big(\frac{H_\sigma g(x)}{S(x)}\Big)^p d\sigma(x)\le (p')^p\int_{(0,\infty)} g^p\,d\sigma .
\end{equation}
\end{lemma}

\begin{proof}
The claim is trivial if the right-hand side of \eqref{eq:one-weight} is infinite. Let $x_0=\inf\{x>0: S(x)>0\}$. Then $\sigma((0,x_0))=0$. Also $\sigma(\{x_0\})=0$ unless $S(x_0)>0$. Thus $S(Y)>0$ almost surely.

Let $D$ be the set of positive rational numbers $d$ with $S(d)>0$. We claim that almost surely
\begin{equation}\label{eq:claim}
A^*(Y)\le \sup_{d\in D} M_d .
\end{equation}
Let us fix an outcome with $S(Y)>0$ and a point $x\ge Y$. Then $S(x)>0$. The functions $H_\sigma g$ and $S$ are right-continuous. Hence $A(x)=\lim A(d)$ as $d\downarrow x$ through rational numbers. Since each such $d$ satisfies $d>x\ge Y$, we have $M_d=A(d)$ by \eqref{eq:condexp}. This gives $A(x)\le\sup_{d\in D}M_d$. Taking the supremum over $x\ge Y$ proves \eqref{eq:claim}.

Now let $F_1\subset F_2\subset\cdots$ be finite subsets of $D$ with union $D$. For each $n$ we order $F_n$ decreasingly. Then $(M_d)_{d\in F_n}$ is a nonnegative martingale. Its last element is $\E[g(Y)\mid\F_{\min F_n}]$. Doob's $L^p$ maximal inequality \cite[Chapter VII, Theorem 3.4]{Do53} and the contraction property of conditional expectation give
\[
\E\Big[\max_{d\in F_n}M_d^{\,p}\Big]\le (p')^p\,\E\big[M_{\min F_n}^{\,p}\big]\le (p')^p\,\E\big[g(Y)^p\big].
\]
By monotone convergence the same bound holds for $\sup_{d\in D}M_d$. Together with \eqref{eq:claim} we obtain
\[
\E\big[A^*(Y)^p\big]\le (p')^p\,\E\big[g(Y)^p\big].
\]
Multiplying by $S(\infty)$ gives \eqref{eq:one-weight}. Finally \eqref{eq:hardy-sigma} follows because $A(x)\le A^*(x)$ whenever $S(x)>0$.
\end{proof}

Inequality \eqref{eq:hardy-sigma} is \cite[Theorem 2.1]{KW21}, written in the language of measures. There it is proved by a discrete Hardy lemma of Broadbent type and an approximation of the distribution function. Our proof by Doob's inequality also gives the maximal form \eqref{eq:one-weight}.

\begin{remark}\label{rem:sharp-one-weight}
For Lebesgue measure on $(0,1]$ the inequality \eqref{eq:hardy-sigma} is the classical Hardy inequality on a bounded interval. The constant $p'$ is sharp there, as the functions $g(x)=x^{-1/p+\varepsilon}$ on $(0,1]$ with $\varepsilon>0$ and $\varepsilon\downarrow 0$ show. Lemma~\ref{lem:doob} loses nothing in the one-weight case.
\end{remark}

\subsection{Quantile coupling and splitting of atoms}\label{sec:coupling}

Let $\sigma$ be a nonzero Borel measure on $(0,\infty)$ that satisfies \eqref{eq:locfin}. Let $\mathrm m$ denote Lebesgue measure on $(0,\infty)$. Let $Q$ be the quantile function of $S$, that is,
\[
Q(t)=\inf\{y>0: S(y)\ge t\},\qquad 0<t<S(\infty).
\]
By right-continuity of $S$ we have $Q(t)\le y$ if and only if $t\le S(y)$. Hence the image of Lebesgue measure on $(0,S(\infty))$ under $Q$ is $\sigma$, and
\begin{equation}\label{eq:quantile}
\int_0^{S(\infty)} F(Q(t))\,dt=\int_{(0,\infty)}F\,d\sigma
\end{equation}
for every nonnegative Borel function $F$.

We write $Hh(s)=\int_0^sh(t)\,dt$. Let $\hat\omega$ be the image of $\omega$ under the map $x\mapsto S(x)$, a Borel measure on $[0,\infty)$. Put $\hat W(s)=\hat\omega([s,\infty))$. A possible mass of $\hat\omega$ at $0$ can be discarded, because $Hh(0)=0$ and the point $0$ enters neither $\hat B$ nor $\hat{\mathcal A}$ below. The measure $\hat\omega$ need not be $\sigma$-finite. In that case the pair $(\mathrm m,\hat\omega)$ may lie outside the standing assumptions of Section~\ref{sec:statement}. Lemma~\ref{lem:coupling} below holds for an arbitrary Borel measure $\hat\omega$, with integrals in $[0,\infty]$. The arguments applied to $(\mathrm m,\hat\omega)$ use only integrals of nonnegative functions with respect to $\hat\omega$. If $\hat\omega(\{s\})=\infty$ for some $s>0$, then $W(x_s)=\infty$ for the first point $x_s$ with $S(x_s)\ge s$. So $B=\infty$. The test function $\mathbf 1_{(0,x_s]}$ gives $C_*=\infty$. Let $\hat B$ and $\hat{\mathcal A}$ be the quantities \eqref{eq:B} and \eqref{eq:A} for the pair $(\mathrm m,\hat\omega)$, that is,
\[
\hat B=\sup_{s>0}\hat W(s)^{1/q}s^{1/p'},\qquad \hat{\mathcal A}=\Big(\int\hat W(s)^{r/p}s^{r/p'}\,d\hat\omega(s)\Big)^{1/r}.
\]

\begin{lemma}[Quantile coupling]\label{lem:coupling}
Let $1<p<\infty$ and $0<q<\infty$. Then the least constant in \eqref{eq:main-ineq} for the pair $(\sigma,\omega)$ equals the least constant for the pair $(\mathrm m,\hat\omega)$. We also have $\hat B=B$. Finally $\hat{\mathcal A}\ge\mathcal A$ if $q<p$.
\end{lemma}

\begin{proof}
Since $Q(t)\le x$ if and only if $t\le S(x)$, the image of $\mathrm m$ restricted to $(0,S(\infty))$ under $Q$ is $\sigma$. This is an identity between $\sigma$-finite measures. If $S(\infty)<\infty$ and $T$ is uniformly distributed on $(0,S(\infty))$, it says that $Q(T)$ has law $\sigma/S(\infty)$.

Let $g\ge0$ be Borel and put $h=g\circ Q$ on $(0,S(\infty))$ and $h=0$ elsewhere. Then $\int h^p\,d\mathrm m=\int g^p\,d\sigma$ and
\[
Hh(S(x))=\int_{\{t:\,Q(t)\le x\}}g(Q(t))\,dt=\int_{(0,x]}g\,d\sigma=H_\sigma g(x).
\]
Hence $\int(H_\sigma g)^q\,d\omega=\int(Hh)^q\,d\hat\omega$. This shows that the constant for $(\mathrm m,\hat\omega)$ works for $(\sigma,\omega)$.

Conversely, let $h\ge0$ be Borel with $\int h^p\,d\mathrm m<\infty$. Only the values of $h$ on $(0,S(\infty))$ matter, because $\hat\omega$ is carried by $[0,S(\infty)]$. Let $\nu$ be the image of $h\,d\mathrm m$ on $(0,S(\infty))$ under $Q$. For every $x>0$ H\"older's inequality gives $\nu((0,x])=\int_0^{S(x)}h\,dt\le S(x)^{1/p'}\big(\int h^p\,d\mathrm m\big)^{1/p}<\infty$. So $h$ is locally integrable on $(0,S(\infty))$. The measure $\nu$ is finite on bounded intervals. If $\sigma(E)=0$, then $\mathrm m(\{t:Q(t)\in E\})=\sigma(E)=0$ and $\nu(E)=0$. Thus the Radon--Nikodym derivative $g=d\nu/d\sigma$ exists. By \eqref{eq:quantile} we have $\int_{\{Q\in E\}}h\,d\mathrm m=\nu(E)=\int_Eg\,d\sigma=\int_{\{Q\in E\}}g\circ Q\,d\mathrm m$ for every Borel set $E$. In other words, $g\circ Q$ is the conditional expectation of $h$ given $Q$. Explicitly, at an atom $y$ of $\sigma$ the value $g(y)$ is the mean of $h$ over the interval $(S(y-),S(y)]$. With $E=(0,x]$ and $\{Q\le x\}=(0,S(x)]\cap(0,S(\infty))$ this gives $\int_{(0,x]}g\,d\sigma=\nu((0,x])=Hh(S(x))$. Jensen's inequality gives $\int g^p\,d\sigma\le\int h^p\,d\mathrm m$. It follows that the constant for $(\sigma,\omega)$ works for $(\mathrm m,\hat\omega)$.

Let $s>0$. The set $\{x:S(x)\ge s\}$ is empty or equals $[x_s,\infty)$ for some $x_s$, because $S$ is right-continuous. In the second case $\hat W(s)=W(x_s)$ and $S(x_s)\ge s$. So $\hat W(s)^{1/q}s^{1/p'}\le W(x_s)^{1/q}S(x_s)^{1/p'}\le B$. Finally the inclusion $\{x':S(x')\ge S(x)\}\supset[x,\infty)$ gives $\hat W(S(x))\ge W(x)$, and
\[
\hat{\mathcal A}^r=\int\hat W(S(x))^{r/p}S(x)^{r/p'}\,d\omega(x)\ge\mathcal A^r .
\]
The same inequality gives $\hat B\ge\hat W(S(y))^{1/q}S(y)^{1/p'}\ge W(y)^{1/q}S(y)^{1/p'}$ whenever $S(y)>0$. Points with $S(y)=0$ contribute nothing to $B$. Hence $\hat B=B$.
\end{proof}

The inequality $\hat{\mathcal A}\ge\mathcal A$ can be strict. For instance, if $\sigma=\mathrm m|_{(0,1]}$ and $\omega=\delta_1+\delta_2$, then $\mathcal A^r=2^{r/p}+1$ and $\hat{\mathcal A}^r=2\cdot2^{r/p}$. Here $\omega$ charges two points between which $S$ is constant. So for upper bounds in the case $q<p$ we need a second reduction, which keeps $\mathcal A$ unchanged.

\begin{lemma}[Splitting atoms]\label{lem:split}
There are an atomless Borel measure $\tilde\sigma$ with $\tilde\sigma((0,z])<\infty$ for every $z>0$ and a $\sigma$-finite Borel measure $\tilde\omega$ on $(0,\infty)$ with the following properties. The least constant in \eqref{eq:main-ineq} for $(\sigma,\omega)$ equals the least constant for $(\tilde\sigma,\tilde\omega)$. The quantities $B$ and $\mathcal A$ for $(\tilde\sigma,\tilde\omega)$ are equal to those for $(\sigma,\omega)$.
\end{lemma}

\begin{proof}
Let $y_1,y_2,\dots$ be the atoms of $\sigma$, put $\ell_n=2^{-n}$. Let $\phi(x)=x+\sum_{n:\,y_n\le x}\ell_n$. The function $\phi$ is strictly increasing. Let $G_n=(\phi(y_n)-\ell_n,\phi(y_n)]$. These intervals are disjoint and do not meet $\phi((0,\infty)\setminus\{y_1,y_2,\dots\})$. Let $\sigma_c$ be the atomless part of $\sigma$ and put
\[
\tilde\sigma=\phi_\#\sigma_c+\sum_n\frac{\sigma(\{y_n\})}{\ell_n}\,\mathrm m|_{G_n},\qquad\tilde\omega=\phi_\#\omega .
\]
The point $\phi(x)$ lies in $G_n$ or to its right when $y_n\le x$, since $\phi(y_n)$ is the right endpoint of $G_n$. It lies to the left of $G_n$ when $y_n>x$. Hence $\tilde\sigma((0,\phi(x)])=S(x)$ and $\tilde\omega([\phi(x),\infty))=W(x)$ for every $x>0$. Since $\phi(x)\ge x$ and $G_n\subset[y_n,\infty)$, we have $\tilde\sigma((0,z])\le S(z)<\infty$ for every $z>0$. As a strictly increasing Borel injection, $\phi$ maps a countable cover of $(0,\infty)$ by sets of finite $\omega$-measure onto Borel sets of finite $\tilde\omega$-measure. Since $\tilde\omega$ vanishes off $\phi((0,\infty))$, it is $\sigma$-finite. Since $\tilde\omega$ is carried by $\phi((0,\infty))$, the quantity $\mathcal A$ is unchanged. A point $z>0$ outside $\phi((0,\infty))$ satisfies $\phi(y_n)-\ell_n\le z<\phi(y_n)$ for some $n$. Then $\tilde\omega([z,\infty))=W(y_n)$ and $\tilde\sigma((0,z])\le S(y_n)$. Such points do not increase $B$. So $B$ is unchanged. Given $g\ge0$, let $\tilde g(\phi(x))=g(x)$, let $\tilde g=g(y_n)$ on $G_n$, and let $\tilde g=0$ at the remaining points, which form a $\tilde\sigma$-null and $\tilde\omega$-null set. Then $\int\tilde g^p\,d\tilde\sigma=\int g^p\,d\sigma$ and $H_{\tilde\sigma}\tilde g(\phi(x))=H_\sigma g(x)$. Since $\tilde\omega$ is carried by $\phi((0,\infty))$, the least constant for $(\sigma,\omega)$ is at most that for $(\tilde\sigma,\tilde\omega)$. Conversely, let $\tilde g\ge0$ be Borel. Put $g(x)=\tilde g(\phi(x))$ if $x$ is not an atom of $\sigma$. Let $g(y_n)$ be the $\tilde\sigma$-average of $\tilde g$ over $G_n$. By Jensen's inequality $\int g^p\,d\sigma\le\int\tilde g^p\,d\tilde\sigma$. For every $x>0$ the interval $(0,\phi(x)]$ contains $G_n$ if $y_n\le x$ and does not meet $G_n$ if $y_n>x$. This gives $H_\sigma g(x)=H_{\tilde\sigma}\tilde g(\phi(x))$. Thus $\int(H_\sigma g)^q\,d\omega=\int(H_{\tilde\sigma}\tilde g)^q\,d\tilde\omega$. The reverse inequality between the least constants follows. This completes the proof.
\end{proof}

\subsection{An elementary lemma}\label{sec:elementary}

\begin{lemma}\label{lem:tail}
Let $I\subset(0,\infty)$ be an interval, let $\nu$ be a finite Borel measure on $I$ with total mass $m$, and let $\gamma>0$. Then
\[
\int_{I}\nu(I\cap[x,\infty))^\gamma\,d\nu(x)\ge\frac{m^{\gamma+1}}{\gamma+1}.
\]
\end{lemma}

\begin{proof}
Put $T(x)=\nu(I\cap[x,\infty))$ and fix $t\ge0$. The set $J=\{x\in I:T(x)\le t\}$ is a terminal subinterval of $I$, because $T$ is nonincreasing. If $J$ contains its left endpoint $c$, then $\nu(J)=T(c)\le t$. Otherwise $\nu(J)=\lim_{x\downarrow c}T(x)\le t$. Hence $\nu(\{T>t\})\ge m-t$. By the layer cake formula
\[
\int T^\gamma\,d\nu=\int_0^\infty\gamma t^{\gamma-1}\,\nu(\{T>t\})\,dt\ge\int_0^m\gamma t^{\gamma-1}(m-t)\,dt=\frac{m^{\gamma+1}}{\gamma+1}.\qedhere
\]
\end{proof}

\section{The case \texorpdfstring{$p\le q$}{p <= q}}\label{sec:pleq}

In this section we prove Theorem~\ref{thm:main} in two steps. Section~\ref{sec:pleq-mart} proves the characterization $B\le C_*\le p^{1/q}p'B$ through the reverse martingale alone. Section~\ref{sec:pleq-sharp} obtains the sharp factor $k_{p,q}$ through a special function evaluated along the reverse martingale. It also records the proof through the quantile coupling and Bliss's inequality as an alternative.

\subsection{The characterization through the reverse martingale}\label{sec:pleq-mart}

\begin{proposition}\label{prop:pleq-mart}
Let $1<p\le q<\infty$. Then $B\le C_*\le p^{1/q}p'B$. In particular \eqref{eq:main-ineq} holds with a finite constant if and only if $B<\infty$.
\end{proposition}

The proof uses only Lemma~\ref{lem:doob}, a Carleson type estimate and the layer cake formula. It uses neither the quantile coupling nor Bliss's inequality.

\begin{proof}
\emph{Necessity.} Let \eqref{eq:main-ineq} hold with a constant $C<\infty$. Fix $y>0$ with $S(y)>0$ and put $g=\mathbf 1_{(0,y]}$. Then $\int g^p\,d\sigma=S(y)<\infty$. For $x\ge y$ we have $H_\sigma g(x)=S(y)$. Hence
\[
S(y)\,W(y)^{1/q}\le\Big(\int_{(0,\infty)}\big(H_\sigma g\big)^q\,d\omega\Big)^{1/q}\le C\,S(y)^{1/p}.
\]
Dividing by $S(y)^{1/p}$ gives $W(y)^{1/q}S(y)^{1/p'}\le C$. If $S(y)=0$, the corresponding term in \eqref{eq:B} vanishes by convention. In both cases $B\le C$. Taking the infimum over all admissible $C$ gives $B\le C_*$.

\emph{Sufficiency.} Let $B<\infty$ and let $g\ge0$ with $\int g^p\,d\sigma<\infty$, which suffices for a bound on $C_*$. First we prove a Carleson type estimate. For $a>0$, Fubini's theorem and $W\le B^qS^{-q/p'}$ on $\{S>0\}$, which carries the measure $d(S^q)$, give
\[
\int_{(0,a]}S^q\,d\omega=\int_{(0,a]}\omega([y,a])\,d(S^q)(y)\le\int_{(0,a]}W\,d(S^q)\le B^q\int_{(0,a]}S^{-q/p'}\,d(S^q)\le pB^qS(a)^{q/p}.
\]
For the last step we split $d(S^q)$ into its continuous part and its jumps. On the continuous part $S^{-q/p'}\,d(S^q)=p\,d(S^{q/p})$ by the chain rule. At a jump of $S$ from $s_0=S(y-)$ to $s_1=S(y)$ we have $s_1^{-q/p'}(s_1^q-s_0^q)\le\int_{s_0}^{s_1}qs^{q/p-1}\,ds=p(s_1^{q/p}-s_0^{q/p})$, because $s\le s_1$ and $q/p'>0$. Summing these contributions gives the last inequality. By monotone convergence $\int_LS^q\,d\omega\le pB^q\sigma(L)^{q/p}$ for every initial interval $L$. Next, $H_\sigma g\le SA^*$. Here $A$ and $A^*$ are defined by \eqref{eq:Ax} and \eqref{eq:Astar} also when $\sigma$ is infinite. For $y\le R$ the maximal average of $\sigma|_{(0,R]}$ at $y$ equals $\sup\{A(x):y\le x\le R,\ S(x)>0\}$. It increases to $A^*(y)$ as $R\to\infty$. The set $L_t=\{A^*>t\}$ is an initial interval because $A^*$ is nonincreasing. Put $\varphi(t)=\sigma(L_t)$ and $N=\int(A^*)^p\,d\sigma=\int_0^\infty pt^{p-1}\varphi(t)\,dt$. Lemma~\ref{lem:doob}, applied to $\sigma$ on $(0,R]$ with $R\to\infty$, gives $N\le(p')^p\int g^p\,d\sigma$. Put $\gamma=q/p\ge1$ and $G(t)=\int_0^tps^{p-1}\varphi(s)\,ds$. Since $\varphi$ is nonincreasing, $G(t)\ge t^p\varphi(t)$. Hence $\frac d{dt}G(t)^\gamma=\gamma G(t)^{\gamma-1}pt^{p-1}\varphi(t)\ge qt^{q-1}\varphi(t)^\gamma$ for almost every $t$. For $q=p$ we have $\gamma=1$. Then this reads $G'(t)=pt^{p-1}\varphi(t)$. Integrating gives $\int_0^\infty qt^{q-1}\varphi(t)^\gamma\,dt\le N^\gamma$. The layer cake formula now gives
\[
\int(H_\sigma g)^q\,d\omega\le\int_0^\infty qt^{q-1}\int_{L_t}S^q\,d\omega\,dt\le pB^q\int_0^\infty qt^{q-1}\varphi(t)^{q/p}\,dt\le pB^qN^{q/p}.
\]
With $N\le(p')^p\int g^p\,d\sigma$ we obtain $\int(H_\sigma g)^q\,d\omega\le p(p')^qB^q\big(\int g^p\,d\sigma\big)^{q/p}$, that is, $C_*\le p^{1/q}p'B$.
\end{proof}

The factor $p^{1/q}p'$ is larger than $k_{p,q}$, for example $2\sqrt2$ against $2$ at $p=q=2$, and $2^{4/3}\approx2.52$ against $k_{2,3}\approx1.468$ at $(p,q)=(2,3)$. The sharp factor is obtained in Section~\ref{sec:pleq-sharp}.

\subsection{The sharp factor}\label{sec:pleq-sharp}

The upper bound uses a special function of \Osek{} \cite{Os15}, evaluated along the reverse martingale of Section~\ref{sec:doob} and its companion for $g^p$. This is an argument in the spirit of Burkholder's method. It needs neither the quantile coupling nor Bliss's inequality. Atoms of $\sigma$ are treated directly. The coupling gives an alternative proof, which we record in Remark~\ref{rem:coupling-proof}. A Bellman function argument with the same constant appears in the work of Ba\~nuelos and \Osek{} \cite[Theorem 3.1]{BO25} on Carleson sequences indexed by a tree, in discrete time and under a testing condition. Their function is defined as a supremum. Its majorization is Bliss's inequality itself. The function $\mathcal U$ defined before Lemma~\ref{lem:U} is explicit. Its majorization follows from the formula for $\varphi$. It runs along a filtration indexed by $x\in(0,\infty)$ under the condition $B<\infty$. On the chain of intervals $(0,x]$ the pair $(\mathrm m,\nu)$ used for optimality below satisfies the testing condition $\int_{(0,x]}s^q\,d\nu(s)\le L^qx^{q/p}$ with $L^q=p-1$. So the two sharp constants agree. For equal exponents on a dyadic tree, under a testing condition, an explicit Bellman function with the sharp constant is due to Cavina \cite{Ca26}. It has four variables, one of which carries the testing functional of the weights. The function $\mathcal U$ has three variables. The weight $\omega$ enters only through the bound $W\le B^qS^{-q/p'}$ in \eqref{eq:Psi}. Lemma~\ref{lem:U} needs $p<q$. The upper bound for $p=q$ in Theorem~\ref{thm:main} follows from it by letting $q\downarrow p$.

Let $1<p<q$ and $\varrho=q/p-1$. By \cite[Lemmas 3.1, 3.2 and 3.3]{Os15}, applied with $k=p$, $\ell=q$ and $\alpha=\varrho$, there are a continuous function $\varphi\colon[1,\infty)\to(0,\infty)$ of class $C^1$ on $(1,\infty)$ and a function $u\colon(1,\infty)\to(0,1)$ with the following properties. First, $\varphi(1)=p/q$ and $\varphi'>0$. Second, $(p-1)(1-u)\big(s\varphi'-\tfrac qp\varphi\big)=1$ on $(1,\infty)$, by \cite[(3.2) and (3.6)]{Os15}. Third, $\varphi(s)s^{-q/p}\to C_{p,q,\varrho}$ as $s\to\infty$, where $C_{p,q,\varrho}$ is the constant of \cite[(1.3)]{Os15}. Also $(q/p')C_{p,q,\varrho}=k_{p,q}^q$. Fourth, for $s>1$ and $\xi\ge0$,
\begin{equation}\label{eq:os39}
\big(q\varphi-ps\varphi'\big)(\xi-1)+\varphi'\big(\xi^p-s\big)+\tfrac qp\varphi\ge1,
\end{equation}
by \cite[(3.9) and (3.10)]{Os15}, where $\varphi$ and $\varphi'$ are evaluated at $s$. By the second property, \eqref{eq:os39} is equivalent to $\xi_*\le u$, since the minimum over $\xi\ge0$ of its left-hand side is $(p-1)\big(s\varphi'-\frac qp\varphi\big)(1-\xi_*)$, attained at $\xi_*=\big((s\varphi'-\frac qp\varphi)/\varphi'\big)^{1/(p-1)}$. Put $\psi=1+(q/p')\varphi$ and
\[
\mathcal U(x,y,t)=t^{q/p}x^q\,\psi\big(y/x^p\big),\qquad x>0,\ y\ge x^p,\ t>0 .
\]

\begin{lemma}\label{lem:U}
Let $1<p<q$. The function $\mathcal U$ is continuous, of class $C^1$ on $\{x>0,\ y>x^p,\ t>0\}$, and has the following properties.
\begin{enumerate}
\item For $x>0$, $y>x^p$, $t>0$ and $d\ge0$,
\[
\mathcal U_x\,\frac{d-x}t+\mathcal U_y\,\frac{d^p-y}t+\mathcal U_t\ge qt^{\varrho}x^{q-1}d .
\]
\item $\mathcal U(x,x^p,t)=pt^{q/p}x^q$.
\item $\mathcal U(x,y,t)\le k_{p,q}^q(ty)^{q/p}$.
\end{enumerate}
\end{lemma}

\begin{proof}
Put $s=y/x^p$ and $\xi=d/x$. Then $\mathcal U_x=t^{q/p}x^{q-1}(q\psi-ps\psi')$, $\mathcal U_y=t^{q/p}x^{q-p}\psi'$ and $\mathcal U_t=\frac qpt^{\varrho}x^q\psi$. Dividing the left-hand side of (1) by $t^{\varrho}x^q$ and inserting $\psi=1+(q/p')\varphi$ we obtain
\[
q(\xi-1)+\frac qp+\frac q{p'}\Big(\big(q\varphi-ps\varphi'\big)(\xi-1)+\varphi'\big(\xi^p-s\big)+\frac qp\varphi\Big)\ge q\xi-q+\frac qp+\frac q{p'}=q\xi
\]
by \eqref{eq:os39}, which is (1). Part (2) follows from $\psi(1)=1+(q/p')(p/q)=p$. For (3) note that $\mathcal U(x,y,t)=(ty)^{q/p}\psi(s)s^{-q/p}$ and
\[
s^{q/p+1}\frac d{ds}\big(\psi(s)s^{-q/p}\big)=s\psi'-\frac qp\psi=\frac q{p'}\Big(s\varphi'-\frac qp\varphi\Big)-\frac qp=\frac qp\Big(\frac1{1-u}-1\Big)>0 .
\]
Thus $\psi(s)s^{-q/p}$ increases to its limit $(q/p')C_{p,q,\varrho}=k_{p,q}^q$, which gives (3).
\end{proof}

\begin{proof}[Proof of Theorem~\ref{thm:main}]
The equivalence and the lower bound $B\le C_*$ are in Proposition~\ref{prop:pleq-mart}. It remains to prove $C_*\le k_{p,q}B$ and the optimality of $k_{p,q}$.

\emph{Upper bound for $p<q$.} Let $p<q$ and $B<\infty$. Let $g\ge0$ with $G=\int g^p\,d\sigma<\infty$. Put $H=H_\sigma g$. By H\"older's inequality $H(x)\le S(x)^{1/p'}G^{1/p}<\infty$. The function $H$ is right-continuous and nondecreasing with $H=0$ on $\{S=0\}$. By Fubini's theorem and $W\le B^qS^{-q/p'}$ on $\{S>0\}$,
\begin{equation}\label{eq:Psi}
\int_{(0,\infty)}H^q\,d\omega=\int_{(0,\infty)}W\,d(H^q)\le B^q\,\Psi_\sigma(g),\qquad \Psi_\sigma(g)=\int_{\{S>0\}}S^{-q/p'}\,d(H^q).
\end{equation}
We show that $\Psi_\sigma(g)\le k_{p,q}^qG^{q/p}$. If $g\le\tilde g$, then $H_\sigma\tilde g\ge H$ and $H_\sigma\tilde g-H$ is nondecreasing. The difference $(H_\sigma\tilde g)^q-H^q$ is nondecreasing as well, because $(h+\delta)^q-h^q$ is nondecreasing in $h\ge0$ and in $\delta\ge0$ for $q\ge1$. Hence $d\big((H_\sigma\tilde g)^q\big)\ge d(H^q)$ as measures. As a result, $\Psi_\sigma$ is monotone in $g$. Replacing $g$ by $g+\varepsilon e$, where $e>0$ and $\int e^p\,d\sigma<\infty$, and letting $\varepsilon\to0$, we may assume $g>0$. Let $I=\{S>0\}$, an interval with left endpoint $a_0=\inf I$. On $I$ put
\[
A=\frac HS,\qquad P(x)=\frac1{S(x)}\int_{(0,x]}g^p\,d\sigma,\qquad Z=\mathcal U(A,P,S).
\]
In the notation of Section~\ref{sec:doob}, $A(x)$ is the value of the reverse martingale $M_x$ on $\{Y\le x\}$. Likewise $P(x)$ is the corresponding value for $g^p$. By H\"older's inequality $P\ge A^p>0$, with $P(x)=A(x)^p$ if and only if $g$ is $\sigma$-almost everywhere constant on $(0,x]$. Let $J_0$ be the set of these $x\in I$, an initial subinterval of $I$, possibly empty, on which $g=c$ $\sigma$-almost everywhere for one constant $c>0$. For an interval $L$ write $\Psi_L=\int_{L\cap I}S^{-q/p'}\,d(H^q)$. We claim that
\begin{equation}\label{eq:Zmono}
\Psi_{(a,b]}\le Z(b)-Z(a)\qquad\text{for all } a<b \text{ in } I,\qquad \Psi_{I\cap(0,b]}\le Z(b)\qquad\text{for all } b\in I .
\end{equation}

\emph{Jumps.} Let $y\in I$ be an atom of $\sigma$ with $y>a_0$, put $m=\sigma(\{y\})$, $d=g(y)$, $S_0=S(y-)>0$, $S_1=S(y)$, $X_0=A(y-)$ and $Y_0=P(y-)$. For $\tau\in[S_0,S_1]$ let $X(\tau)=(S_0X_0+(\tau-S_0)d)/\tau$ and $Y(\tau)=(S_0Y_0+(\tau-S_0)d^p)/\tau$, so that $X(S_1)=A(y)$, $Y(S_1)=P(y)$, $X'=(d-X)/\tau$ and $Y'=(d^p-Y)/\tau$. If $Y_0=X_0^p$ and $d=X_0$, then the jump of $Z$ is $pX_0^q(S_1^{q/p}-S_0^{q/p})$ by Lemma~\ref{lem:U}(2). It is at least $X_0^qS_1^{-q/p'}(S_1^q-S_0^q)$, the mass of $\Psi$ at $y$, by the estimate for jumps in the proof of Proposition~\ref{prop:pleq-mart}. Otherwise $Y(\tau)>X(\tau)^p$ for $\tau\in(S_0,S_1]$. Lemma~\ref{lem:U}(1) with $t=\tau$ gives $\frac d{d\tau}\mathcal U(X,Y,\tau)\ge q\tau^{\varrho}X^{q-1}d$ there. Since $\mathcal U$ is continuous and $\tau X(\tau)=H(y-)+(\tau-S_0)d$, integration gives
\[
Z(y)-Z(y-)\ge\int_{S_0}^{S_1}q\tau^{\varrho}X^{q-1}d\,d\tau=\int_{S_0}^{S_1}\tau^{-q/p'}\,d\big((\tau X)^q\big)\ge S_1^{-q/p'}\big(H(y)^q-H(y-)^q\big),
\]
because $\varrho-q+1=-q/p'$ and $\tau\le S_1$. The right-hand side is the mass of $\Psi$ at $y$. If $a_0\in I$, that is, $\sigma(\{a_0\})=m>0$, then $a_0\in J_0$ and $Z(a_0)=pm^{q/p}g(a_0)^q$ by Lemma~\ref{lem:U}(2). The mass of $\Psi$ at $a_0$ is $m^{q/p}g(a_0)^q$.

\emph{Continuous part.} Let $b_0=\sup J_0$, with $b_0=a_0$ if $J_0$ is empty. If $J_0=I$, the estimate on $J_0$ below already gives the first part of \eqref{eq:Zmono}. The localization that follows is then not needed. On $J_0$ we have $H=cS$ and $Z=pc^qS^{q/p}$. For $a<b$ in $J_0$ the estimate for $\int S^{-q/p'}\,d(S^q)$ in the proof of Proposition~\ref{prop:pleq-mart}, applied on $(a,b]$, gives $\Psi_{(a,b]}=c^q\int_{(a,b]}S^{-q/p'}\,d(S^q)\le pc^q\big(S(b)^{q/p}-S(a)^{q/p}\big)=Z(b)-Z(a)$. Now let $y>b_0$. Then $(0,y)$ contains a point outside $J_0$. Hence $g$ is not $\sigma$-almost everywhere constant on $(0,y)$. The strict form of H\"older's inequality gives $P(y-)>A(y-)^p$. Hence on $(b_0,\infty)$ the path $x\mapsto(A,P,S)$ and its left limits stay in the open set $\{(x,y,t):x>0,\ y>x^p,\ t>0\}$, where $\mathcal U$ is of class $C^1$. The point $b_0$ itself may belong to $I\setminus J_0$. Then $b_0$ is an atom of $\sigma$, because otherwise $A$ and $P$ would be continuous at $b_0$ and $P(b_0)=A(b_0)^p$. Its left limit lies on the boundary $\{y=x^p\}$. This jump is covered by the jump estimate above. The path is of bounded variation on compact subintervals of $I\setminus\{a_0\}$. Its continuous parts are $dA^c=(g-A)\,d\sigma^c/S$, $dP^c=(g^p-P)\,d\sigma^c/S$ and $dS^c=d\sigma^c$, where $\sigma^c$ is the continuous part of $\sigma$, singular continuous components included. Let $b_0<a'<b$. Since the path is right-continuous with left limits, its values and its left limits on $[a',b]$ form a compact subset of the open set above, on which $\mathcal U$ is of class $C^1$. By the change of variables formula for functions of bounded variation \cite[Chapter~0, Section~4]{RY99}, applied on $[a',b]$, $Z(b)-Z(a')$ is the sum of the jumps of $Z$ in $(a',b]$ and of $\int_{(a',b]}\big(\mathcal U_x\,dA^c+\mathcal U_y\,dP^c+\mathcal U_t\,dS^c\big)$. By Lemma~\ref{lem:U}(1) with $d=g$, the integrand is at least $qS^{\varrho}A^{q-1}g\,d\sigma^c=S^{-q/p'}qH^{q-1}\,dH^c$, which is the continuous part of $S^{-q/p'}\,d(H^q)$. Together with the jump estimate this gives $\Psi_{(a',b]}\le Z(b)-Z(a')$. Letting $a'\downarrow a$ for $a\ge b_0$, the right continuity of $Z$ and monotone convergence give $\Psi_{(a,b]}\le Z(b)-Z(a)$ for $b_0\le a<b$. In general we split $(a,b]$ into its parts in $J_0$, in $\{b_0\}\setminus J_0$ and in $(b_0,\infty)$. The part in $J_0$ is an interval $(a,b_1)$ or $(a,b_1]$ with $b_1\le b_0$. Its bound follows from the case $a<b$ in $J_0$ by letting the right endpoint increase. Adding the three bounds proves the first part of \eqref{eq:Zmono}. For the second part let $a\downarrow a_0$. If $a_0\notin I$, then $\Psi$ has no mass at $a_0$ and $Z\ge0$. If $a_0\in I$, the first part and the mass at $a_0$ computed above give the claim.

By \eqref{eq:Zmono} and Lemma~\ref{lem:U}(3),
\[
\Psi_{I\cap(0,b]}\le Z(b)\le k_{p,q}^q\big(S(b)P(b)\big)^{q/p}\le k_{p,q}^qG^{q/p}
\]
for every $b\in I$. Letting $b\to\infty$ gives $\Psi_\sigma(g)\le k_{p,q}^qG^{q/p}$. Then \eqref{eq:Psi} gives $C_*\le k_{p,q}B$.

\emph{Upper bound for $p=q$.} Let $B<\infty$ and let $g$ be as above. Write $\Psi^{(q)}_\sigma(g)$ for the functional in \eqref{eq:Psi} with exponent $q$. The estimate \eqref{eq:Psi} holds for $q=p$. The bound $\Psi^{(q)}_\sigma(g)\le k_{p,q}^qG^{q/p}$ for $q>p$ involves only $\sigma$ and $g$. Let $\varpi=dH^c+\sum_y\delta_y$, where $H^c$ is the continuous part of $H$ and the sum runs over the countably many jumps $y$ of $H$. This measure does not depend on $q$. With respect to $\varpi$ the measure in $\Psi^{(q)}_\sigma(g)$ has the density $qS^{-q/p'}H^{q-1}$ on the continuous part and the value $S(y)^{-q/p'}(H(y)^q-H(y-)^q)$ at each jump $y$. These densities converge pointwise as $q\downarrow p$. Fatou's lemma gives $\Psi^{(p)}_\sigma(g)\le\liminf_{q\downarrow p}k_{p,q}^qG^{q/p}$. With $\varrho=q/p-1$ we have $1/p-1/q=\varrho/q$. Stirling's formula gives $\log\mathrm B(1/\varrho,(q-1)/\varrho)=\varrho^{-1}\big((q-1)\log(q-1)-q\log q\big)+O(\log(1/\varrho))$. It follows that $\log k_{p,q}\to\log p-(1/p')\log(p-1)$, that is, $k_{p,q}\to(p-1)^{1/p}p'=p^{1/p}(p')^{1/p'}=k_{p,p}$. So $C_*\le k_{p,p}B$.

\emph{Optimality of the factor.} Let $\nu$ be the measure on $(0,\infty)$ with density $(q/p')s^{-q/p'-1}$, so that $\nu((x,\infty))=x^{-q/p'}$. Take $\sigma=\mathrm m$ and $\omega=\nu$. Then $W(x)=x^{-q/p'}$. So $B=1$. For $p<q$ and $h(t)=(1+t^\varrho)^{-(\varrho+1)/\varrho}$ both sides of \eqref{eq:main-ineq} are finite. Their ratio is $k_{p,q}$ by the case of equality in Bliss's inequality \cite{Bl30}, see \eqref{eq:bliss} below and \cite[Section 3.3]{Os15}. For $p=q$ the best constant for this pair is $(p-1)^{1/p}p'=k_{p,p}$ by the sharpness of the constant in \eqref{eq:hardy-classical}. Thus no smaller factor works for all pairs.
\end{proof}

\begin{remark}[Alternative proof through the coupling]\label{rem:coupling-proof}
The upper bound $C_*\le k_{p,q}B$ also follows from the quantile coupling and Bliss's inequality, as follows. By Lemma~\ref{lem:coupling} we may take $\sigma=\mathrm m$ and replace $\omega$ by $\hat\omega$, because $\hat B=B$. We write $\omega$, $W$ and $B$ for $\hat\omega$, $\hat W$ and $\hat B$. There is nothing to prove if $B=\infty$. Let $h\ge0$ with $\int h^p\,d\mathrm m<\infty$. Let $\nu$ be the measure on $(0,\infty)$ with density $(q/p')s^{-q/p'-1}$, so that $\nu((x,\infty))=x^{-q/p'}$. The function $F=(Hh)^q$ is finite, continuous and nondecreasing, with $F(0+)=0$, because $Hh(x)\le x^{1/p'}\|h\|_{L^p}$. Hence for $t>0$ the set $\{F>t\}$ is empty or an interval $(x,\infty)$ with $x>0$, and
\[
\omega(\{F>t\})=\lim_{y\downarrow x}W(y)\le B^q\,x^{-q/p'}=B^q\,\nu(\{F>t\}).
\]
Since $F$ vanishes at $0$, the layer cake formula gives $\int F\,d\omega\le B^q\int F\,d\nu$. For $p<q$ Bliss's inequality \cite{Bl30}, in the form of \cite[Lemma 1]{LM20}, states that
\begin{equation}\label{eq:bliss}
\Big(\int_0^\infty\big(Hh\big)^q\,d\nu\Big)^{1/q}\le k_{p,q}\Big(\int_0^\infty h^p\,d\mathrm m\Big)^{1/p},
\end{equation}
with equality for $h(t)=(1+t^\varrho)^{-(\varrho+1)/\varrho}$. A proof of \eqref{eq:bliss} by a special function argument in the spirit of Burkholder's method is given by \Osek{} \cite[Theorem 1.1 and Section 3]{Os15}. With $k=p$, $\ell=q$ and $\alpha=\varrho$ his constant satisfies $\big((q-\varrho-1)C_{k,\ell,\alpha}\big)^{1/q}=k_{p,q}$. His extremal functions are the dilates of $h$. For $p=q$ we have $d\nu=(p-1)x^{-p}dx$. Then \eqref{eq:bliss} with the factor $(p-1)^{1/p}p'=p^{1/p}(p')^{1/p'}$ is Hardy's inequality \eqref{eq:hardy-classical}. Hence $C_*\le k_{p,q}B$.
\end{remark}

\begin{remark}\label{rem:LM-proof}
Li and Mao \cite{LM20} also use the quantile transform, together with a comparison of tails as in Remark~\ref{rem:coupling-proof}. If $W$ is finite, Sinnamon \cite[Theorem 2.6]{Si22} computes the normal form parameter $b$ of \eqref{eq:main-ineq} from the distribution function of $\sigma$ and the tail of $\omega$. In general $b$ is the nonincreasing rearrangement of $S$ with respect to $\omega$, see \cite[Theorem 2.4]{Si22}. His Theorem 3.8 is the comparison with Bliss's parameter used above. Lemma~\ref{lem:coupling} expresses the same exact equivalence as a conditional expectation. If $S(\infty)<\infty$ and $T$ is uniformly distributed on $(0,S(\infty))$, it reads $\E[h(T)\mid Q(T)]=g(Q(T))$ almost surely. In general $g\circ Q$ is the conditional expectation of $h$ with respect to the $\sigma$-field generated by $Q$ on $(0,S(\infty))$ with Lebesgue measure, which is $\sigma$-finite on this $\sigma$-field because $\{Q\le x\}=(0,S(x)]\cap(0,S(\infty))$.
\end{remark}

\begin{remark}\label{rem:Os}
The proof in Remark~\ref{rem:coupling-proof} gives $C_*=k_{p,q}B$ whenever $1<p\le q$, $\sigma$ has no atoms, $S(\infty)=\infty$ and $W^{1/q}S^{1/p'}=B$ on $\{S>0\}$. Indeed, if $B<\infty$, then $\hat\omega=B^q\nu$ on $(0,\infty)$. If $B=\infty$, both $C_*$ and $k_{p,q}B$ are infinite. The optimality part of the proof of Theorem~\ref{thm:main} applies. Power weights with matched powers $W=cS^{-q/p'}$ satisfy this condition. In this way we recover the sharp inequality of \cite[Theorem 1.2]{Os15}, where the change of variables is now the quantile coupling. We also recover the sharp constant $(p/(\gamma-1))^p$ of \cite[Theorem 330]{HLP52} for $w(t)=t^{-\gamma}$ and $v(t)=t^{p-\gamma}$, $\gamma>1$.
\end{remark}

\section{The case \texorpdfstring{$q<p$}{q < p}}\label{sec:qlp}

Throughout this section $1<p<\infty$, $0<q<p$ and $1/r=1/q-1/p$. We first prove Theorem~\ref{thm:qlp}. Then we prove Theorems~\ref{thm:sharp}, \ref{thm:lowerB}, \ref{thm:lowerq} and \ref{thm:Dlambda}.

\subsection{Two reductions}\label{sec:reductions}
Let $\sigma$ be atomless. Then $S$ is continuous. The image of $\sigma|_{(0,x]}$ under $S$ is $\mathrm m|_{(0,S(x)]}$. Hence
\begin{equation}\label{eq:atomless}
\int_{(0,x]}S^\gamma\,d\sigma=\frac{S(x)^{\gamma+1}}{\gamma+1}\qquad(\gamma>-1).
\end{equation}
For $g\ge0$ and $F=H_\sigma g$, Lemma~\ref{lem:doob}, applied to $\sigma$ on $(0,R]$ with $R\to\infty$, gives
\begin{equation}\label{eq:hardy-beta}
\int_{\{S>0\}}\Big(\frac FS\Big)^pd\sigma\le(p')^p\int g^p\,d\sigma .
\end{equation}
For atomless $\sigma$ with $S(\infty)=\infty$, which is the case wherever \eqref{eq:hardy-beta} is used below, it is also the classical inequality \eqref{eq:hardy-classical}. Indeed, with $h=g\circ Q$ as in the proof of Lemma~\ref{lem:coupling} we have $\int h^p\,d\mathrm m=\int g^p\,d\sigma$ and $Hh(S(x))=F(x)$. Also $\sigma$ is carried to $\mathrm m$ by $S$. Lemma~\ref{lem:doob} gives \eqref{eq:hardy-beta} for every $\sigma$, atoms included, and also gives the maximal form.
If $\mathcal A<\infty$, then $W(y)<\infty$ whenever $S(y)>0$. Indeed, let $S(y)>0$ and let $\nu\le\omega|_{[y,\infty)}$ be a finite measure. Lemma~\ref{lem:tail} applied to $\nu$ on $[y,\infty)$ gives
\[
\nu([y,\infty))^{r/q}S(y)^{r/p'}\le\frac rq\int\nu([x,\infty))^{r/p}S(x)^{r/p'}\,d\nu(x)\le\frac rq\,\mathcal A^r .
\]
Since $\omega$ is $\sigma$-finite, the supremum over all such $\nu$ gives $W(y)^{r/q}S(y)^{r/p'}\le(r/q)\mathcal A^r$. Now let $S(u)>0$. Then $W(u)<\infty$. Lemma~\ref{lem:tail} on $[u,\infty)$ with $\gamma=r/p$ gives $W(u)^{r/q}\le(r/q)\int_{[u,\infty)}W^{r/p}\,d\omega$, because $r/p+1=r/q$. Since $\sigma$ does not charge $\{S=0\}$, Fubini's theorem and \eqref{eq:atomless} with $\gamma=r/p'-1>-1$ give
\begin{equation}\label{eq:ibp}
\int W^{r/q}S^{r/p'-1}\,d\sigma\le\frac rq\int W(x)^{r/p}\int_{(0,x]}S^{r/p'-1}\,d\sigma\,d\omega(x)=\frac{p'}{q}\,\mathcal A^r .
\end{equation}

Upper bounds of the form $C_*\le c\,\mathcal A$ reduce to this case. By Lemma~\ref{lem:split} we may assume that $\sigma$ is atomless. For $n\in\N$ let $\sigma_n$ be $\sigma$ on $(0,n]$ plus $\mathrm m$ on $(n,\infty)$. Let $\omega_n$ be the restriction of $\omega$ to $(0,n]$. Then $\mathcal A(\sigma_n,\omega_n)\le\mathcal A$. For $g$ vanishing on $(n,\infty)$ we have $\int_{(0,n]}(H_\sigma g)^q\,d\omega=\int(H_{\sigma_n}g)^q\,d\omega_n$. By monotone convergence, a bound $C_*\le c\,\mathcal A$ for atomless $\sigma$ with $S(\infty)=\infty$ gives the same bound in general.

\subsection{Proof of Theorem~\ref{thm:qlp}}\label{sec:SS}
The two bounds in \eqref{eq:constants-qlp} are proved separately. Together they show that \eqref{eq:main-ineq} holds with a finite constant if and only if $\mathcal A<\infty$.

\begin{proof}[Proof of the lower bound in \eqref{eq:constants-qlp}]
By Lemma~\ref{lem:coupling} it suffices to show $\hat{\mathcal A}\le c\,C$ for the pair $(\mathrm m,\hat\omega)$, where $c^{-1}=(p')^{1/p'}q^{1/p}(1-q/p)$ and $C$ is any constant for that pair. Let $0<a<b<\infty$, let $\omega_{ab}$ be the restriction of $\hat\omega$ to $[a,b]$. Let $W_{ab}$ and $\mathcal A_{ab}$ be the corresponding tail and integral. The constant $C$ also works for $\omega_{ab}$. The test function $\mathbf 1_{(0,a]}$ gives $a^q\,\hat\omega([a,b])\le C^qa^{q/p}$. Thus $\omega_{ab}$ is finite. We may assume $0<\mathcal A_{ab}$. Also $\mathcal A_{ab}<\infty$ because $\omega_{ab}$ is finite. Put $h(t)=W_{ab}(t)^{r/(pq)}t^{r/(pq')}$. Since $r/(pq')+1=r/(p'q)>0$ and $W_{ab}$ is nonincreasing,
\[
Hh(x)\ge W_{ab}(x)^{r/(pq)}\int_0^xt^{r/(pq')}\,dt=\frac{p'q}{r}\,W_{ab}(x)^{r/(pq)}x^{r/(p'q)},
\]
and so $\int(Hh)^q\,d\omega_{ab}\ge(p'q/r)^q\mathcal A_{ab}^r$. Now $h^p=W_{ab}^{r/q}t^{r/p'-1}$, because $r/q'=r/p'-1$. Lemma~\ref{lem:tail} on $[t,\infty)$ with $\gamma=r/p$ gives $W_{ab}(t)^{r/q}\le(r/q)\int_{[t,\infty)}W_{ab}^{r/p}\,d\omega_{ab}$, since $r/p+1=r/q$. By Fubini's theorem
\begin{equation}\label{eq:ibp-lower}
\int h^p\,d\mathrm m\le\frac rq\int W_{ab}(x)^{r/p}\int_0^xt^{r/p'-1}\,dt\,d\omega_{ab}(x)=\frac{p'}{q}\,\mathcal A_{ab}^r .
\end{equation}
The inequality for $h$ now gives $(p'q/r)\mathcal A_{ab}^{r/q}\le C(p'/q)^{1/p}\mathcal A_{ab}^{r/p}$, that is, $\mathcal A_{ab}\le cC$, since $q/r=1-q/p$. Letting $a\downarrow0$ and $b\uparrow\infty$, monotone convergence gives $\hat{\mathcal A}\le cC$. Lemma~\ref{lem:coupling} gives $\mathcal A\le\hat{\mathcal A}$.
\end{proof}

\begin{proof}[Proof of the upper bound in \eqref{eq:constants-qlp}]
Let $\mathcal A<\infty$. By Section~\ref{sec:reductions} we may assume that $\sigma$ is atomless with $S(\infty)=\infty$, so that \eqref{eq:atomless} holds and $W(y)<\infty$ whenever $S(y)>0$. Put $d\eta=S^{-p}\mathbf 1_{\{S>0\}}\,d\sigma$ and $\mathcal T(y)=\eta([y,\infty))$. By \eqref{eq:atomless}, $\mathcal T(y)=S(y)^{1-p}/(p-1)$ where $S(y)>0$. Let $F=H_\sigma g$. Then \eqref{eq:hardy-beta} reads $\int F^p\,d\eta\le(p')^p\int g^p\,d\sigma$. Let $\Phi(t)=\int_{(0,t]}W^{r/p}\mathcal T^{-r/q}\,d\eta$, a nondecreasing function. Where $\Phi(t)=0$ and $S(t)>0$, we have $W(t)=0$. The set of such $t$ is $\omega$-null. Where $S(t)=0$ we have $F(t)=0$. All integrals against $\omega$ below are taken over $\{\Phi>0\}$. H\"older's inequality with the exponents $r/q$ and $p/q$ gives
\[
\int F^q\,d\omega\le\Big(\int\Phi\,d\omega\Big)^{q/r}\Big(\int\Phi^{-p/r}F^p\,d\omega\Big)^{q/p}.
\]
By Fubini's theorem, $\int\Phi\,d\omega=\int W^{r/q}\mathcal T^{-r/q}\,d\eta=(p-1)^{r/q}\int W^{r/q}S^{r/p'-1}\,d\sigma$, because $(p-1)r/q-p=r/p'-1$. By \eqref{eq:ibp}, $\int\Phi\,d\omega\le(p-1)^{r/q}(p'/q)\mathcal A^r$.

For the second factor, let $x>0$ with $S(x)>0$. Since $\Phi$ is nondecreasing and $W$ is nonincreasing, \eqref{eq:atomless} gives
\begin{align*}
\int_{[x,\infty)}\Phi^{-p/r}\,d\omega&\le\Phi(x)^{-p/r}W(x)\le\Big(\int_{(0,x]}\mathcal T^{-r/q}\,d\eta\Big)^{-p/r}\\
&=\Big((p-1)^{r/q}\frac{p'}{r}S(x)^{r/p'}\Big)^{-p/r}=\Big(\frac rp\Big)^{p/r}\mathcal T(x).
\end{align*}
Where $S(x)=0$ the right-hand side is infinite. The function $F^p$ is nondecreasing and continuous with $F(0+)=0$. Thus $F^p(t)=\nu((0,t])$ for a Borel measure $\nu$. By Fubini's theorem
\[
\int\Phi^{-p/r}F^p\,d\omega=\int\Big(\int_{[y,\infty)}\Phi^{-p/r}\,d\omega\Big)d\nu(y)\le\Big(\frac rp\Big)^{p/r}\int\mathcal T\,d\nu=\Big(\frac rp\Big)^{p/r}\int F^p\,d\eta .
\]
Combining the two factors with \eqref{eq:hardy-beta}, we obtain
\[
C_*\le(p-1)^{1/q}\Big(\frac{p'}{q}\Big)^{1/r}\Big(\frac rp\Big)^{1/r}p'\,\mathcal A=\Big(\frac rq\Big)^{1/r}(p-1)^{1/p}p'\,\mathcal A=\Big(\frac rq\Big)^{1/r}p^{1/p}(p')^{1/p'}\mathcal A ,
\]
where we used $p'/p=1/(p-1)$ and $1/q-1/r=1/p$.
\end{proof}

\begin{remark}
The argument follows \cite[Theorem 2.4]{SS96}, with Lemma~\ref{lem:doob} in the role of Hardy's inequality. Lemmas~\ref{lem:coupling} and \ref{lem:split} transfer it to arbitrary measures, with atoms in both $\sigma$ and $\omega$, without any loss in the constants. The proof includes the quasi-Banach range $0<q<1$. No duality in the target space is used. The best upper factors are determined in Section~\ref{sec:sharp}.
\end{remark}

\subsection{The sharp upper factor}\label{sec:sharp}

We first prove \eqref{eq:sharp} for $0<q\le1$, together with the sharpness of $(r/q)^{1/r}$. In this part $0<q\le1<p$. We put $s=r/p'$. Then $0<s\le1$, because $r/p'\le1$ is equivalent to $q\le1$. We also use the identities $(p'-1)s=r/p$ and $1+r/p=r/q$.

\begin{proof}[Proof of Theorem~\ref{thm:sharp}(1) and (2) for $q\le1$]
The case $\mathcal A=\infty$ is trivial. By Lemma~\ref{lem:split} we may assume that $\sigma$ is atomless, so that $S$ is continuous. Next we localize $\omega$. Let $a_0=\inf\{x:S(x)>0\}$ and $Z=\{x:S(x)>0\}$. Since $S$ is continuous, $Z=(a_0,\infty)$. Since $H_\sigma g=0$ outside $Z$, replacing $\omega$ by its restriction to $Z$ changes neither side of \eqref{eq:main-ineq} and does not increase $\mathcal A$. For $a\in Z$ and $b>a$ let $\omega_{ab}=\omega|_{[a,b]}$. Since the interval is closed, atoms of $\omega$ at $a$ and at $b$ are kept. With the convention $W(x)=\omega([x,\infty))$, which counts an atom at $x$, the tail of $\omega_{ab}$ is
\[
W_{ab}(x)=\omega([x,\infty)\cap[a,b])=\begin{cases}\omega([a,b]),&x\le a,\\ \omega([x,b]),&a<x\le b,\\ 0,&x>b.\end{cases}
\]
Thus $W_{ab}\le W$ and $\omega_{ab}\le\omega$, so that $\mathcal A(\sigma,\omega_{ab})\le\mathcal A$. The measure $\omega_{ab}$ is also finite. Indeed, let $\nu\le\omega_{ab}$ be a finite measure. Lemma~\ref{lem:tail} applied to $\nu$ on $[a,b]$ gives
\[
\nu([a,b])^{r/q}S(a)^{r/p'}\le\frac rq\int\nu([x,\infty))^{r/p}S(x)^{r/p'}\,d\nu(x)\le\frac rq\,\mathcal A^r .
\]
Since $\omega$ is $\sigma$-finite, the supremum over all such $\nu$ gives $\omega([a,b])^{r/q}S(a)^{r/p'}\le(r/q)\mathcal A^r$. Also $S(a)>0$. As $b\uparrow\infty$ and $a\downarrow a_0$, the measures $\omega_{ab}$ increase to $\omega|_Z$. Also $\int(H_\sigma g)^q\,d\omega_{ab}$ increases to $\int(H_\sigma g)^q\,d\omega$ by monotone convergence. Hence it suffices to prove \eqref{eq:sharp} for finite $\omega$. Then $W\le\omega((0,\infty))<\infty$.

Let $g\ge0$ and $F=H_\sigma g$. The set $\{x:W(x)=0\}$ is $\omega$-null. Let $x>0$ with $W(x)>0$. For $u\le x$ we have $W(u)\ge W(x)>0$. H\"older's inequality applied to $g=gW^{-1/p}\cdot W^{1/p}$ gives
\[
F(x)\le G(x)^{1/p}R(x)^{1/p'},\qquad G(x)=\int_{(0,x]}g^pW^{-1}\,d\sigma,\quad R(x)=\int_{(0,x]}W^{p'-1}\,d\sigma .
\]
H\"older's inequality with the exponents $r/q$ and $p/q$ gives
\[
\int F^q\,d\omega\le\int R^{q/p'}G^{q/p}\,d\omega\le\Big(\int R^{s}\,d\omega\Big)^{q/r}\Big(\int G\,d\omega\Big)^{q/p}.
\]
By Fubini's theorem $\int G\,d\omega\le\int g(u)^pW(u)^{-1}\omega([u,\infty))\,d\sigma(u)=\int g^p\,d\sigma$, where the integral is taken over $\{W>0\}$. Hence
\begin{equation}\label{eq:E-bound}
C_*\le E^{1/r},\qquad E=\int R^s\,d\omega .
\end{equation}

It remains to show $E\le(r/q)\mathcal A^r$. Since $\sigma$ is atomless and $W$ is bounded, $R$ is finite, continuous and nondecreasing. Also $R(0+)=0$. Let $\rho$ be the Borel measure with $\rho((0,t])=R(t)^s$. By Fubini's theorem
\[
E=\int\rho((0,t])\,d\omega(t)=\int\omega([u,\infty))\,d\rho(u)=\int W\,d\rho .
\]
Since $R$ is continuous, $\rho$ is carried by $\{R>0\}$ and $d\rho=sR^{s-1}\,dR=sR^{s-1}W^{p'-1}\,d\sigma$ there. When $q=1$ we have $s=1$. We use the convention $0^0=1$. Then $R^{s-1}=1$ also on $\{R=0\}$ and $d\rho=dR$. Since $\rho$ does not charge $\{R=0\}$, the convention affects no integral below. The inequality $R^{s-1}\le(W^{p'-1}S)^{s-1}$ used next reads $1\le1$. Since $W$ is nonincreasing, $R(u)\ge W(u)^{p'-1}S(u)$. Because $s-1\le0$, it follows that $R(u)^{s-1}\le W(u)^{(p'-1)(s-1)}S(u)^{s-1}$ wherever $W(u)>0$. Hence
\[
E\le s\int W^{p'+(p'-1)(s-1)}S^{s-1}\,d\sigma=s\int W^{r/q}S^{r/p'-1}\,d\sigma ,
\]
because $p'+(p'-1)(s-1)=1+(p'-1)s=1+r/p=r/q$. Lemma~\ref{lem:tail} on $[u,\infty)$ with $\gamma=r/p$ gives $W(u)^{r/q}\le(r/q)\int_{[u,\infty)}W^{r/p}\,d\omega$. By Fubini's theorem and \eqref{eq:atomless} with $\gamma=s-1>-1$,
\[
E\le s\,\frac rq\int W(x)^{r/p}\int_{(0,x]}S^{s-1}\,d\sigma\,d\omega(x)=\frac rq\int W^{r/p}S^{r/p'}\,d\omega=\frac rq\,\mathcal A^r .
\]
Together with \eqref{eq:E-bound} this proves \eqref{eq:sharp}.

For sharpness let $\sigma=\mathrm m$, let $T>0$, $\ell>0$ and $0<\varepsilon<T$. Let $\omega_\varepsilon=(\ell/\varepsilon)\,\mathrm m|_{[T-\varepsilon,T]}$. The function $g=\mathbf 1_{(0,T-\varepsilon]}$ satisfies $H_\sigma g=T-\varepsilon$ on $[T-\varepsilon,T]$. Thus
\[
C_*\ge\frac{\ell^{1/q}(T-\varepsilon)}{(T-\varepsilon)^{1/p}}=\ell^{1/q}(T-\varepsilon)^{1/p'} .
\]
Next, $W(x)=\ell(T-x)/\varepsilon$ for $x\in[T-\varepsilon,T]$. The substitution $x=T-\varepsilon y$ gives
\[
\mathcal A^r\le T^{r/p'}\ell^{r/p+1}\int_0^1y^{r/p}\,dy=\frac qr\,T^{r/p'}\ell^{r/q}.
\]
Thus $C_*/\mathcal A\ge(r/q)^{1/r}(1-\varepsilon/T)^{1/p'}$. Letting $\varepsilon\to0$ shows that the factor $(r/q)^{1/r}$ cannot be decreased.
\end{proof}

We turn to $1\le q<p$. Put $\theta=(q-1)/(p-1)\in[0,1)$. Then
\begin{equation}\label{eq:theta}
\theta+\frac{1-\theta}{p}+\frac qr=1,\qquad p\theta+1-\theta=q,\qquad (p-1)\theta=q-1 .
\end{equation}

\begin{proof}[Proof of Theorem~\ref{thm:sharp}(1) for $1\le q<p$]
The argument follows Maz'ya's proof of \cite[Section 1.3.3, Lemma]{Ma11}, with Lemma~\ref{lem:doob} in place of Hardy's inequality in Step~2. Suppose $\mathcal A<\infty$. By Section~\ref{sec:reductions} we may assume that $\sigma$ is atomless and $S(\infty)=\infty$, and that $W(y)<\infty$ whenever $S(y)>0$. Then $S$ is continuous and \eqref{eq:atomless} holds. By monotone convergence it suffices to consider bounded $g\ge0$ that vanish outside $(0,n]$ for some $n$. Then $F=H_\sigma g$ is continuous, nondecreasing and bounded. Also $F=0$ where $S=0$.

\emph{Step 1. A representation.} Since $F$ is continuous and $q\ge1$, we have $F(x)^q=\int_{(0,x]}qF^{q-1}g\,d\sigma$. By Fubini's theorem
\begin{equation}\label{eq:Fq-rep}
\int F^q\,d\omega=q\int F^{q-1}g\,W\,d\sigma ,
\end{equation}
where the integrand vanishes outside $\{S>0\}$.

\emph{Step 2. The one-weight inequality.} On $\{S>0\}$ the functions $F^p$ and $S^{1-p}$ are continuous and of locally bounded variation. Also $d(F^pS^{1-p})=pF^{p-1}gS^{1-p}\,d\sigma-(p-1)F^pS^{-p}\,d\sigma$. H\"older's inequality gives $F(x)^pS(x)^{1-p}\le\int_{(0,x]}g^p\,d\sigma$, which tends to $0$ as $x$ decreases to the left endpoint of $\{S>0\}$, because $\sigma$ does not charge $\{S=0\}$. Since $F$ is bounded and $S(\infty)=\infty$, also $F^pS^{1-p}\to0$ at infinity. Integrating, and using \eqref{eq:hardy-beta}, we obtain
\begin{equation}\label{eq:martingale-identity}
\int_{\{S>0\}}F^{p-1}g\,S^{1-p}\,d\sigma=\frac1{p'}\int_{\{S>0\}}\Big(\frac FS\Big)^pd\sigma\le(p')^{p-1}\int g^p\,d\sigma .
\end{equation}

\emph{Step 3. H\"older's inequality.} On $\{S>0\}$ we factor
\[
F^{q-1}gW=\big(F^{p-1}gS^{1-p}\big)^\theta\cdot g^{1-\theta}\cdot S^{q-1}W ,
\]
which is correct by \eqref{eq:theta}. H\"older's inequality with the exponents $1/\theta$, $p/(1-\theta)$ and $r/q$ (two exponents if $\theta=0$), together with \eqref{eq:Fq-rep} and \eqref{eq:martingale-identity}, gives
\[
\int F^q\,d\omega\le q\,(p')^{(p-1)\theta}\Big(\int g^p\,d\sigma\Big)^{\theta+(1-\theta)/p}\Big(\int S^{(q-1)r/q}W^{r/q}\,d\sigma\Big)^{q/r}.
\]
Here $(q-1)r/q=r/q'=r/p'-1\ge0$. Up to this point we have proved, for the reduced pair, $C_*\le q^{1/q}(p')^{1/q'}\mathcal B_\sigma$ with $\mathcal B_\sigma^r=\int W^{r/q}S^{r/q'}\,d\sigma$. This is Maz'ya's upper bound \cite[Section 1.3.3]{Ma11}, which is stated there for absolutely continuous $\sigma$. Now \eqref{eq:ibp} bounds the last integral by $(p'/q)\mathcal A^r$.
By \eqref{eq:theta} the exponent of $\int g^p\,d\sigma$ is $q/p$. Taking $q$-th roots,
\[
C_*\le q^{1/q}(p')^{(q-1)/q}\Big(\frac{p'}q\Big)^{1/r}\mathcal A=(p')^{1/p'}q^{1/p}\mathcal A ,
\]
because $1/q-1/r=1/p$ and $1-1/q+1/r=1/p'$.
\end{proof}

\begin{remark}[Step 2 without Doob's inequality]\label{rem:bellman-step2}
The proof above uses \eqref{eq:martingale-identity} in Step 3 only through the bound $\int_{\{S>0\}}F^{p-1}gS^{1-p}\,d\sigma\le(p')^{p-1}\int g^p\,d\sigma$. This bound also follows from a pointwise inequality, without Lemma~\ref{lem:doob}. Put $ a_p=p^p/(p-1)^{p-1}$, so that $ a_p/p=(p')^{p-1}$. Put $G(x)=\int_{(0,x]}g^p\,d\sigma$ and $\mathcal U(x,y,t)=t\big( a_p y-(p-1)x^p\big)$. On $\{S>0\}$ let
\[
V=\mathcal U\big(F/S,\,G/S,\,S\big)= a_p G-(p-1)F^pS^{1-p} .
\]
With $a=F/S$ the formula for $d(F^pS^{1-p})$ in Step 2 gives
\[
dV-pF^{p-1}gS^{1-p}\,d\sigma=\big( a_p g^p-p^2a^{p-1}g+(p-1)^2a^p\big)\,d\sigma\ge0 ,
\]
because the minimum of $ a_p g^p-p^2a^{p-1}g$ over $g\ge0$ is $-(p-1)^2a^p$, attained at $g=a/p'$. As $x$ decreases to the left endpoint of $\{S>0\}$ we have $V(x)\to0$ by Step 2. Also $V\le a_p G$ everywhere. Integrating over $\{S>0\}$ we obtain $p\int_{\{S>0\}}F^{p-1}gS^{1-p}\,d\sigma\le a_p\int g^p\,d\sigma$, which is the required bound. This is the integration by parts proof of Hardy's inequality, written with the function $\mathcal U$ of the state $(F/S,G/S,S)$ in the manner of Section~\ref{sec:pleq-sharp}. It needs no a priori finiteness of $\int(F/S)^p\,d\sigma$. Together with Steps 1 and 3 it gives $C_*\le(p')^{1/p'}q^{1/p}\mathcal A$ for $1\le q<p$ from Maz'ya's factorization and an explicit function, with no maximal inequality.
\end{remark}

\begin{proof}[Proof of $\sup C_*/\mathcal A\ge(p')^{1/p'}q^{1/p}$ for all $0<q<p$, and of Theorem~\ref{thm:sharp}(2) for $1\le q<p$]
We show that $\sup C_*/\mathcal A\ge(p')^{1/p'}q^{1/p}$ for every $0<q<p$. Let $\sigma=\mathrm m$, let $N>1$, and let
\[
\omega_N=\frac q{p'}\,x^{-q/p'-1}\mathbf 1_{[1,N]}(x)\,dx+N^{-q/p'}\delta_N .
\]
Then $W(x)=x^{-q/p'}$ for $1\le x\le N$. Also $W(x)=1$ for $0<x<1$. On $[1,N]$ we have $W^{r/p}x^{r/p'}=x^{(r/p')(q/r)}=x^{q/p'}$, because $1-q/p=q/r$. Hence
\[
\mathcal A^r=\frac q{p'}\int_1^N\frac{dx}x+1=\frac q{p'}\log N+1 .
\]
Let $g(t)=t^{-1/p}\mathbf 1_{[1,N]}(t)$. Then $\int g^p\,dt=\log N$ and $H g(x)=p'(x^{1/p'}-1)$ for $1\le x\le N$. This gives
\[
\int(Hg)^q\,d\omega_N\ge q(p')^{q-1}\int_1^N\big(1-x^{-1/p'}\big)^q\frac{dx}x\ge q(p')^{q-1}(\log N-c),
\]
where $c=\int_1^\infty\big(1-(1-x^{-1/p'})^q\big)\,dx/x<\infty$, since $1-(1-y)^q\le\max(1,q)\,y$ for $0\le y\le1$. Combining the last two displays gives
\[
\frac{C_*}{\mathcal A}\ge\frac{\big(q(p')^{q-1}(\log N-c)\big)^{1/q}}{(\log N)^{1/p}\big((q/p')\log N+1\big)^{1/r}}\longrightarrow q^{1/q}(p')^{1/q'}\Big(\frac{p'}q\Big)^{1/r}=(p')^{1/p'}q^{1/p}
\]
as $N\to\infty$. For $1\le q<p$ this equals $\kappa_{p,q}$. Hence $\kappa_{p,q}$ cannot be decreased.

The same bound holds with atomless $\omega$. For $0<\varepsilon<N-1$ let $\omega_{N,\varepsilon}$ be $\omega_N$ with the atom replaced by $N^{-q/p'}\varepsilon^{-1}\mathbf 1_{[N-\varepsilon,N]}(x)\,dx$. Let $W_\varepsilon$ be its tail. Since moving mass to the left does not increase tails, $W_\varepsilon\le W$ on $[1,N]$. Hence $\mathcal A(\omega_{N,\varepsilon})^r\le\int W^{r/p}x^{r/p'}\,d\omega_{N,\varepsilon}$. The functions $W^{r/p}x^{r/p'}=x^{q/p'}$ and $(Hg)^q$ are continuous on $[1,N]$. The measures $\omega_{N,\varepsilon}$ converge weakly to $\omega_N$ there as $\varepsilon\to0$. It follows that $\limsup_{\varepsilon\to0}\mathcal A(\omega_{N,\varepsilon})\le\mathcal A(\omega_N)$ and $\int(Hg)^q\,d\omega_{N,\varepsilon}\to\int(Hg)^q\,d\omega_N$. So the lower bound above for $C_*/\mathcal A$ holds for the atomless pairs $(\mathrm m,\omega_{N,\varepsilon})$ with $\varepsilon$ small, up to an arbitrarily small error.
\end{proof}

\begin{proof}[Proof of Theorem~\ref{thm:sharp}(3)]
Let $q=1$, so that $r=p'$. The upper bound and its sharpness are parts (1) and (2). For the lower bound we may assume $C_*<\infty$. By Fubini's theorem $\int H_\sigma g\,d\omega=\int gW\,d\sigma$. Thus $C_*=\big(\int W^{p'}\,d\sigma\big)^{1/p'}$ by the duality of $L^p(\sigma)$ and $L^{p'}(\sigma)$. For absolutely continuous $\sigma$ this is the case $q=1$ of \cite[Section 1.3.3, Lemma]{Ma11}. Since $W$ is nonincreasing, $\int_{[u,\infty)}W^{p'-1}\,d\omega\le W(u)^{p'-1}\omega([u,\infty))=W(u)^{p'}$ for every $u>0$. By Fubini's theorem
\[
\mathcal A^{p'}=\int W(x)^{p'-1}S(x)\,d\omega(x)=\int\Big(\int_{[u,\infty)}W^{p'-1}\,d\omega\Big)d\sigma(u)\le\int W^{p'}\,d\sigma=C_*^{p'} .
\]
If $\omega=\ell\delta_T$ and $\sigma((0,T])>0$, then $C_*=\ell\,S(T)^{1/p'}=\mathcal A$. So the factor $1$ is best possible. If $\omega$ has no atoms, then $\int_{[u,\infty)}W^{p'-1}\,d\omega=W(u)^{p'}/p'$. The same computation then gives $C_*=(p')^{1/p'}\mathcal A$ for every $\sigma$. This means that every atomless $\omega$ with $0<\mathcal A<\infty$ attains the upper factor.
\end{proof}

\begin{remark}\label{rem:sharp-compare}
Let $\sigma=\mathrm m$ and let $\omega$ be absolutely continuous. Then Sinnamon's functional from \cite[Theorem 3.9]{Si22} is $C_0=(p')^{1/p'}q^{1/p}\mathcal A$, by an integration by parts. In normal form, Theorem~\ref{thm:sharp} reads $N_{p,q}(b)\le\rho_{p,q}C_0$ with $\rho_{p,q}=(p')^{-1/p'}q^{-1/p}(r/q)^{1/r}$ for $q\le1$. Sinnamon \cite[Corollary 3.21]{Si22} proved $N_{p,q}(b)\le\max(1,(p')^{(1-q)/q})\,C_0$. For $0<q<1$ this reads $C_*\le(p')^{1/r}q^{1/p}\mathcal A$. For weights this upper estimate was proved earlier in \cite[Theorem 1(3)]{Si91}. For $q<1$ we have $\rho_{p,q}<(p')^{(1-q)/q}$, for example $\rho_{2,1/2}\approx1.540$ against $2$. This comparison concerns absolutely continuous $\omega$. For general $\omega$ the discussion after Theorem~\ref{thm:sharp} gives $C_0\le(p')^{1/p'}q^{1/p}\mathcal A$. Hence his estimate gives $C_*\le(p')^{1/r}q^{1/p}\mathcal A$ for all $\sigma$ and $\omega$ when $0<q<1$. For some $\omega$ with atoms it is better than \eqref{eq:sharp}. Let $(p,q)=(2,1/2)$ and let $\omega$ be a point mass at a point $T$ with $S(T)>0$. Then $C_*=\mathcal A$ and $C_0=(3/4)^{3/2}\mathcal A$. His estimate gives $C_*\le2C_0=\frac{3\sqrt3}4\,\mathcal A\approx1.299\,\mathcal A$. Theorem~\ref{thm:sharp} gives $C_*\le\frac{8\sqrt3}9\,\mathcal A\approx1.540\,\mathcal A$. This does not contradict Theorem~\ref{thm:sharp}(2), which concerns the supremum of $C_*/\mathcal A$ over all pairs with $0<\mathcal A<\infty$. For $1\le q<p$, Theorem~\ref{thm:sharp}(2) shows that his bound $N_{p,q}(b)\le C_0$ cannot be improved, by the atomless measures $\omega_{N,\varepsilon}$ at the end of that proof. For general $\omega$ a point mass gives $C_*=\mathcal A$ for every $0<q<p$. Thus no lower factor valid for all $\omega$ exceeds $1$. For absolutely continuous $\omega$ the lower factor can exceed $1$. For example $C_*\ge(r/q)^{1/r}\mathcal A$ for $1\le q<p$, see the discussion after the proof of Theorem~\ref{thm:lowerB}. For $0<q<1$ Sinnamon's factor $(p')^{1/p'}q^{1/p}$ from \cite[Theorem 1(3)]{Si91} equals $1.216\ldots$ for $p=3$ and $q=4/5$. For $q<1$ the factor $1$ fails by Theorem~\ref{thm:lowerq}(2). For a numerical illustration, let $p=2$, $q=1/10$, $\sigma=\mathrm m$ and $\omega=0.3\,\delta_1+0.03\,\delta_{10}+\delta_{20}$. Then $C_*$ is the maximum of a concave function on a convex set. A numerical computation gives $C_*/\mathcal A\approx0.939$. The maximizing $g$ with $\int g^2\,d\mathrm m=1$ takes the values $0.5478\ldots$, $0.1969\ldots$ and $0.1873\ldots$ on $(0,1]$, $(1,10]$ and $(10,20]$. It vanishes on $(20,\infty)$. This is consistent with Theorem~\ref{thm:lowerq}(2), since $\lambda^\sharp_{2,1/10}=0.4455\ldots$, attained at $t=0.7396\ldots$. For $\mathbb T=\mathbb R$, $a=0$ and $0<q<1$, the functional of \cite[Theorem 1]{MSAK26} is Maz'ya's functional $\mathcal B_\sigma$. Two integrations by parts give $\mathcal B_\sigma=(p'/q)^{1/r}\mathcal A$ for absolutely continuous $\sigma$ and $\omega$. So their upper bound reads $C_*\le(p'/q)^{1/r}\mathcal A$. Here $(r/q)^{1/r}<(p'/q)^{1/r}$ because $r<p'$ when $q<1$.
\end{remark}

\subsection{The lower bound for \texorpdfstring{$1\le q<p$}{1 <= q < p}}\label{sec:lowerB}
In this subsection we prove Theorem~\ref{thm:lowerB}. The proof follows the best constant as mass is added to $\omega$ from the right. Lemma~\ref{lem:EL} is the first variation in the direction of the indicator of $(0,x]$.

\begin{lemma}\label{lem:EL}
Let $1\le q<p$ and $x>0$ with $S(x)>0$. Let $\nu$ be a finite sum of point masses at points of $[x,\infty)$. Let $M=\nu([x,\infty))>0$. Let $C(\nu)$ be the best constant in \eqref{eq:main-ineq} for the pair $(\sigma,\nu)$. Then $0<C(\nu)<\infty$. There is $g\ge0$ with $\int g^p\,d\sigma=1$ and $\int(H_\sigma g)^q\,d\nu=C(\nu)^q$ such that
\begin{equation}\label{eq:EL}
H_\sigma g(x)^q\ge C(\nu)^{-rq/p}\,M^{r/p}\,S(x)^{r/p'} .
\end{equation}
\end{lemma}

\begin{proof}
Let $x=y_0<y_1<\dots<y_m$ be the distinct points of the set consisting of $x$ and the atoms of $\nu$. Put $J_0=(0,y_0]$ and $J_k=(y_{k-1},y_k]$. Replacing $g$ on each $J_k$ by its $\sigma$-average does not change $H_\sigma g$ at the points $y_k$. By Jensen's inequality it does not increase $\int g^p\,d\sigma$. Replacing $g$ by $0$ on $(y_m,\infty)$ does not change $H_\sigma g$ at these points either. Hence $C(\nu)^q$ is the maximum of a continuous function over the compact set $\{(g_k):g_k\ge0,\ \sum_kg_k^p\sigma(J_k)=1\}$, where $k$ runs over the indices with $\sigma(J_k)>0$. It is attained at some $g$ that vanishes on $(y_m,\infty)$ and is constant on each $J_k$ with $\sigma(J_k)>0$. Clearly $C(\nu)\le S(y_m)^{1/p'}M^{1/q}<\infty$. Also $C(\nu)\ge S(x)^{1/p'}M^{1/q}>0$ by the test function $\mathbf 1_{(0,x]}$. Let $g_0$ be the value of $g$ on $J_0$, so that $H_\sigma g(x)=g_0S(x)$. For $t>-g_0$ let
\[
\Phi(t)=\frac{\int(H_\sigma g+tS(x)\mathbf 1_{[x,\infty)})^q\,d\nu}{\big(1+((g_0+t)^p-g_0^p)S(x)\big)^{q/p}},
\]
which is the ratio for $g+t\mathbf 1_{(0,x]}$. Put $\Psi=\int(H_\sigma g)^{q-1}\,d\nu$ if $q>1$ and $\Psi=\nu([x,\infty))=M$ if $q=1$, so that no power $0^0$ occurs. In both cases $\Psi>0$. For $q=1$ this is $M>0$. For $q>1$ it follows from $\int(H_\sigma g)^q\,d\nu=C(\nu)^q>0$. If $g_0=0$, then $\Phi'(0+)=qS(x)\Psi>0$, because the derivative of the denominator vanishes at $t=0$ for $p>1$. This contradicts maximality. Hence $g_0>0$ and $\Phi'(0)=0$, that is,
\[
qS(x)\Psi=C(\nu)^q\,q\,g_0^{p-1}S(x),\qquad\text{so}\qquad g_0^{p-1}=C(\nu)^{-q}\Psi .
\]
On $[x,\infty)$ we have $H_\sigma g\ge H_\sigma g(x)=g_0S(x)$. Since $q\ge1$, this gives $\Psi\ge M\,(g_0S(x))^{q-1}$. With the last display we obtain $g_0^{p-q}\ge C(\nu)^{-q}MS(x)^{q-1}$. Since $q/(p-q)=r/p$, $1/(p-q)=r/(pq)$ and $(q-1)/(p-q)+1=(p-1)/(p-q)=r/(p'q)$, this gives
\[
H_\sigma g(x)=g_0S(x)\ge C(\nu)^{-r/p}M^{r/(pq)}S(x)^{r/(p'q)},
\]
which is \eqref{eq:EL}.
\end{proof}

\begin{proof}[Proof of Theorem~\ref{thm:lowerB}]
\emph{The inequality $\mathcal A\le\mathcal B$.} For $y\ge x$ we have $W(y)\le W(x)$. Hence
\[
\int_{[x,\infty)}W^{r/p}\,d\omega\le W(x)^{r/p}W(x)=W(x)^{\beta}.
\]
Since $S(0+)=0$, we have $S(y)^\alpha=\int_{(0,y]}d(S^\alpha)$. By Fubini's theorem,
\begin{align*}
\mathcal A^r&=\int\Big(\int_{(0,y]}d(S^\alpha)\Big)W(y)^{r/p}\,d\omega(y)\\
&=\int\Big(\int_{[x,\infty)}W^{r/p}\,d\omega\Big)d(S^\alpha)(x)\le\mathcal B^r .
\end{align*}

\emph{Finite atomic $\omega$.} Let $\omega=\sum_{i=1}^nw_i\delta_{x_i}$ as in \eqref{eq:Bcal-discrete}. For $1\le i\le n+1$ let $\omega^{(i)}=\sum_{j\ge i}w_j\delta_{x_j}$, so that $\omega^{(n+1)}=0$. Write $C(\nu)$ for the best constant of the pair $(\sigma,\nu)$. We show that
\begin{equation}\label{eq:B-step}
C(\omega^{(i)})^r\ge C(\omega^{(i+1)})^r+S_i^\alpha\big(W_i^\beta-W_{i+1}^\beta\big),\qquad 1\le i\le n .
\end{equation}
Summing over $i$ and using \eqref{eq:Bcal-discrete} gives $C_*=C(\omega^{(1)})\ge\mathcal B$. If $S_i=0$, then \eqref{eq:B-step} holds because $C$ is monotone in the measure. Let $S_i>0$. For $0\le t\le w_i$ put $\nu_t=t\delta_{x_i}+\omega^{(i+1)}$ and $\varphi(t)=C(\nu_t)^q$. Then
\[
\varphi(t)=\sup_{\|g\|_{L^p(\sigma)}=1}\Big(t\,H_\sigma g(x_i)^q+\int(H_\sigma g)^q\,d\omega^{(i+1)}\Big)
\]
is a supremum of affine functions of $t$ whose slopes lie in $[0,S_i^{q/p'}]$, because $H_\sigma g(x_i)\le S_i^{1/p'}\|g\|_{L^p(\sigma)}$ by H\"older's inequality. Thus $\varphi$ is convex, nondecreasing and Lipschitz on $[0,w_i]$. Let $0<t<w_i$ be a point where $\varphi$ is differentiable. Let $g_t$ be the function of Lemma~\ref{lem:EL} for $x=x_i$ and $\nu=\nu_t$, with $M=W_{i+1}+t$. The affine function attached to $g_t$ touches $\varphi$ at $t$ and lies below it. So $\varphi'(t)\ge H_\sigma g_t(x_i)^q\ge\varphi(t)^{-r/p}(W_{i+1}+t)^{r/p}S_i^{\alpha}$. The function $\varphi^{r/q}$ is Lipschitz on $[\varepsilon,w_i]$ for every $\varepsilon>0$. Since $r/q-1-r/p=0$, at such $t$
\[
\frac{d}{dt}\varphi(t)^{r/q}=\frac rq\,\varphi(t)^{r/q-1}\varphi'(t)\ge\frac rq\,(W_{i+1}+t)^{r/p}S_i^\alpha=S_i^\alpha\,\frac{d}{dt}(W_{i+1}+t)^{\beta}.
\]
Integrate over $[\varepsilon,w_i]$ and let $\varepsilon\to0$, using the continuity of $\varphi$. This gives \eqref{eq:B-step}, because $\varphi(w_i)^{r/q}=C(\omega^{(i)})^r$ and $\varphi(0)^{r/q}=C(\omega^{(i+1)})^r$.

\emph{General $\omega$.} If $C_*=\infty$, there is nothing to prove. Let $a_0=\inf\{x:S(x)>0\}$. The measure $d(S^\alpha)$ does not charge $(0,a_0)$. So $\mathcal B$ depends only on $\omega$ restricted to $[a_0,\infty)$, or to $(a_0,\infty)$ when $S(a_0)=0$. Let $a<b$ with $S(a)>0$. Then $\omega([a,b))<\infty$, since otherwise the test function $\mathbf 1_{(0,a]}$ gives $C_*=\infty$. Let $\omega_{a,b}=\omega|_{[a,b)}$. Choose nested finite partitions $a=y_0^{(n)}<\dots<y^{(n)}_{m_n}=b$ whose mesh tends to $0$ and whose union contains every atom of $\omega_{a,b}$. Let $\omega_n$ place the mass $\omega([y^{(n)}_{k-1},y^{(n)}_k))$ at $y^{(n)}_{k-1}$. Every function $H_\sigma g$ is nondecreasing. Thus $C(\omega_n)\le C(\omega_{a,b})\le C_*$. Also $W_n(x)=\omega([y,b))$ for the least partition point $y\ge x$. This tends to $\omega_{a,b}([x,\infty))$ for every $x$. Indeed, if $x$ is an atom, then it is eventually a partition point. Otherwise $\omega(\{x\})=0$. By the finite case and Fatou's lemma, $\mathcal B(\omega_{a,b})^r\le\liminf_n\mathcal B(\omega_n)^r\le C_*^r$. Letting $b\uparrow\infty$ and $a\downarrow a_0$ (with $a=a_0$ allowed when $S(a_0)>0$), the monotone convergence theorem gives $\mathcal B\le C_*$.

\emph{Equality.} If $\omega=\ell\delta_T$, then $\mathcal A^r=\mathcal B^r=\ell^{\beta}S(T)^{\alpha}$. Also $C_*=\ell^{1/q}S(T)^{1/p'}$ by the case of equality in H\"older's inequality.
\end{proof}

For $q=1$ we have $\alpha=1$ and $\beta=p'$. Then $\mathcal B^{p'}=\int W^{p'}\,d\sigma=C_*^{p'}$ by the proof of Theorem~\ref{thm:sharp}(3). So $\mathcal B=C_*$ for every pair when $q=1$.

Let $\omega$ have no atoms and let $W(x)<\infty$ whenever $S(x)>0$. Then $\int_{[x,\infty)}W^{r/p}\,d\omega=W(x)^\beta/\beta$ for such $x$, since $W$ is nonincreasing, $W(y)\to0$ as $y\to\infty$ and $\beta-1=r/p$. In other words $-d(W^\beta)=\beta W^{r/p}\,d\omega$. The Fubini argument in the proof of $\mathcal A\le\mathcal B$ gives $\mathcal B=(r/q)^{1/r}\mathcal A$. If $W(x)=\infty$ for some $x$ with $S(x)>0$, then $C_*\ge\mathcal B=\infty$ and the bound below is trivial. Hence $C_*\ge(r/q)^{1/r}\mathcal A$ for atomless $\omega$. Together with Theorem~\ref{thm:sharp}(1), this gives $(r/q)^{1/r}\le C_*/\mathcal A\le(p')^{1/p'}q^{1/p}$ for atomless $\omega$ with $0<\mathcal A<\infty$ and $1\le q<p$. The left factor is best possible. Let $\sigma=\mathrm m$, $0<\varepsilon<T$ and $d\omega_\varepsilon=\ell\varepsilon^{-1}\mathbf 1_{[T-\varepsilon,T]}(x)\,dx$. Let $C_\varepsilon$ and $\mathcal A_\varepsilon$ be the corresponding constants. Since $H_\sigma g$ is nondecreasing, $\int(H_\sigma g)^q\,d\omega_\varepsilon\le\ell\,H_\sigma g(T)^q\le\ell\,T^{q/p'}\|g\|_{L^p}^q$. The test function $\mathbf 1_{(0,T-\varepsilon]}$ gives $C_\varepsilon\ge\ell^{1/q}(T-\varepsilon)^{1/p'}$. Further, $W_\varepsilon(x)=\ell(T-x)/\varepsilon$ on $[T-\varepsilon,T]$. Hence $\mathcal A_\varepsilon^r\to\ell^{r/q}T^{r/p'}\int_0^1u^{r/p}\,du=(q/r)\,\ell^{r/q}T^{r/p'}$ as $\varepsilon\to0$, because $r/p+1=r/q$. It follows that $C_\varepsilon/\mathcal A_\varepsilon\to(r/q)^{1/r}$. If $\sigma$ has no atoms, then $d(S^\alpha)=\alpha S^{r/q'}\,d\sigma$. So $\mathcal B=(r/p')^{1/r}\mathcal B_\sigma$ with $\mathcal B_\sigma^r=\int W^{r/q}S^{r/q'}\,d\sigma$, the form of the Maz'ya--Rozin condition used in \cite[Section 1.3.3, Theorem 2]{Ma11}. The classical two-sided estimate in this form is $q(p'/r)^{1/q'}\mathcal B_\sigma\le C_*\le q^{1/q}(p')^{1/q'}\mathcal B_\sigma$, see \cite[Section 1.3.3, Theorem 2, p.~51]{Ma11}, where it is stated for an arbitrary measure $\omega$ and an absolutely continuous $\sigma$. Theorem~\ref{thm:lowerB} improves the lower factor to $(r/p')^{1/r}$, which is best possible (point masses). The improvement is strict for $q>1$. The two factors agree at $q=1$, where $C_*=\mathcal B_\sigma$. Maz'ya's upper factor is also best possible. Let $\omega_{N,\varepsilon}$ be the atomless measures at the end of the proof of Theorem~\ref{thm:sharp}(2). Then $\omega_{N,\varepsilon}$ and $\sigma=\mathrm m$ are atomless. Thus $\mathcal B_\sigma=(p'/q)^{1/r}\mathcal A$. Also $C_*/\mathcal A$ tends to $(p')^{1/p'}q^{1/p}$ up to an arbitrarily small error. Dividing, we see that $C_*/\mathcal B_\sigma$ approaches $(p')^{1/p'}q^{1/p}(q/p')^{1/r}=q^{1/q}(p')^{1/q'}$. The proof of Theorem~\ref{thm:lowerB} does not use Lemma~\ref{lem:doob}. It is a variational argument. Lemma~\ref{lem:EL} is the Euler--Lagrange equation tested against the indicator of $(0,x]$. For $q=1$ the test function in the proof of Theorem~\ref{thm:lowerq}(1) gives $C_*\ge\mathcal A$ directly. For $1<q<p$ we know no test function proof of $C_*\ge\mathcal A$.

\subsection{The lower factor for \texorpdfstring{$0<q<1$}{0 < q < 1}}\label{sec:lowerq}
In this subsection we prove Theorem~\ref{thm:lowerq}. The first lemma reduces lower bounds to finitely many atoms. The second lemma computes $C_*$ for a family of such configurations.

\begin{lemma}\label{lem:red-lower}
Let $0<q<p$ and $c>0$. Suppose that $C_*\ge c\mathcal A$ whenever $\sigma=\mathrm m$ and $\omega=\sum_{i=1}^nw_i\delta_{x_i}$ with $0<x_1<\dots<x_n$ and $w_i>0$. Then $C_*\ge c\mathcal A$ for all $\sigma$ and $\omega$.
\end{lemma}

\begin{proof}
By Lemma~\ref{lem:coupling} we may assume $\sigma=\mathrm m$, because $C_*$ is unchanged and $\hat{\mathcal A}\ge\mathcal A$. Mass of $\hat\omega$ at $0$ contributes neither to $\mathcal A$ nor to $\int(Hh)^q\,d\hat\omega$. Such mass can be discarded. We also suppose $C_*<\infty$. For $0<a<b$ the test function $\mathbf 1_{(0,a]}$ gives $\omega([a,b))\le C_*^qa^{-q/p'}<\infty$. Let $\omega_{a,b}=\omega|_{[a,b)}$. Let $\omega_n$ be the discretizations from the proof of Theorem~\ref{thm:lowerB}, which place the mass of $[y^{(n)}_{k-1},y^{(n)}_k)$ at $y^{(n)}_{k-1}$. Then $C(\omega_n)\le C_*$. Let $\pi_n(x)=y^{(n)}_{k-1}$ for $x\in[y^{(n)}_{k-1},y^{(n)}_k)$. Then
\[
\mathcal A(\omega_n)^r=\int\omega([\pi_n(x),b))^{r/p}\,\pi_n(x)^{r/p'}\,d\omega_{a,b}(x),
\]
where $\omega([\pi_n(x),b))\ge\omega([x,b))$ and $\pi_n(x)\to x$. By Fatou's lemma, $\mathcal A(\omega_{a,b})\le\liminf_n\mathcal A(\omega_n)\le c^{-1}C_*$. Letting $a\downarrow0$ and $b\uparrow\infty$, the monotone convergence theorem gives $\mathcal A\le c^{-1}C_*$.
\end{proof}

Let $0<q<1$, let $0=a_0<a_1<\dots<a_n$ and let $b_1>\dots>b_n>b_{n+1}=0$. Put
\begin{equation}\label{eq:ab-config}
w_i=(b_i-b_{i+1})a_i^{1-q},\qquad h_i=(a_i-a_{i-1})b_i^{1-p'},\qquad x_i=h_1+\dots+h_i ,
\end{equation}
$\omega=\sum_iw_i\delta_{x_i}$ and $\sigma=\mathrm m$. Then $W(x_i)=W_i=\sum_{j\ge i}w_j$ and $S(x_i)=x_i$. Put
\[
I=\sum_{i=1}^n(b_i-b_{i+1})a_i=\sum_{i=1}^nb_i(a_i-a_{i-1}),\qquad U_i=\frac{W_i}{b_ia_i^{1-q}},\qquad V_i=\frac{x_i}{a_ib_i^{1-p'}} .
\]
Since $a_j\ge a_i$ for $j\ge i$ and $b_k\ge b_i$ for $k\le i$, we have $U_i\ge1$ and $0<V_i\le1$. Also $U_n=1$. The exponents of $a_i$ and $b_i$ cancel because $1-q/p=q/r$ and $1/p'-1=-1/p$, so that
\begin{equation}\label{eq:A-UV}
\mathcal A^r=\sum_{i=1}^nw_iW_i^{r/p}x_i^{r/p'}=\sum_{i=1}^n(b_i-b_{i+1})a_i\,U_i^{r/p}V_i^{r/p'} .
\end{equation}

\begin{lemma}\label{lem:ab}
Let $1<p<\infty$ and $0<q<1$.
For the pair defined by \eqref{eq:ab-config} we have $C_*=I^{1/r}$.
\end{lemma}

\begin{proof}
Let $g_0=b_i^{p'-1}$ on $(x_{i-1},x_i]$ and $g_0=0$ on $(x_n,\infty)$. Then $Hg_0(x_i)=a_i$, $\int g_0^p\,d\mathrm m=\sum_ih_ib_i^{p'}=I$ and $\int(Hg_0)^q\,d\omega=\sum_iw_ia_i^q=I$. Hence $C_*\ge I^{1/q-1/p}=I^{1/r}$. Conversely, let $g\ge0$ and $F_i=Hg(x_i)$, with $F_0=0$. Since $q<1$, we have $F_i^q\le a_i^q+qa_i^{q-1}(F_i-a_i)$. Also $w_ia_i^{q-1}=b_i-b_{i+1}$. Summation by parts and Young's inequality $b_ig\le g^p/p+g_0^p/p'$ on $(x_{i-1},x_i]$, where $g_0^{p-1}=b_i$, give
\[
\sum_i(b_i-b_{i+1})F_i=\sum_ib_i(F_i-F_{i-1})\le\frac1p\int g^p\,d\mathrm m+\frac1{p'}I .
\]
Thus $\int(Hg)^q\,d\omega\le(1-q)I+q\sum_i(b_i-b_{i+1})F_i\le(1-q/p)I+(q/p)\int g^p\,d\mathrm m$. Let $\int g^p\,d\mathrm m=1$ and apply this to $I^{1/p}g$. This gives $\int(Hg)^q\,d\omega\le I^{1-q/p}=I^{q/r}$. So $C_*\le I^{1/r}$.
\end{proof}

\begin{proof}[Proof of Theorem~\ref{thm:lowerq}(1)]
By Lemma~\ref{lem:red-lower} we may assume that $\sigma=\mathrm m$ and that $\omega$ is a finite sum of point masses in $(0,\infty)$. Put $\gamma=(p-1)/(p-q)$, $\Gamma(x)=(W(x)x^{p-1})^{1/(p-q)}$ and $F(x)=\sup_{0<y\le x}\Gamma(y)$. Then $F(0+)=0$. Also $F$ is constant on $[x_n,\infty)$. Near $0$ we have $F(x)=\omega((0,\infty))^{1/(p-q)}x^\gamma$. Here $\gamma>1/p'$ because $q>0$, which makes $F'$ $p$-integrable near $0$. For $0<y<z$ we have $\Gamma(z)\le\Gamma(y)(z/y)^\gamma$, because $W$ is nonincreasing. Hence $F(y)\le F(x)\le F(y)(x/y)^\gamma$ for $y<x$. In particular $F$ is locally Lipschitz and $F'\le\gamma F/x$ almost everywhere. The function $\Gamma$ is upper semicontinuous, because $W(y)\to W(x)$ as $y\uparrow x$ and $\lim_{y\downarrow x}W(y)\le W(x)$. Hence $F$ is constant near every point where $F>\Gamma$, so that $F'=0$ almost everywhere on $\{F>\Gamma\}$. On $\{F=\Gamma\}$ we have $F^{p-1}x^{1-p}=F^{q-1}W$. Combining the two cases, we obtain
\[
\int_0^\infty F'^p\,dx\le \gamma^{p-1}\int_0^\infty F'F^{q-1}W\,dx=\frac{\gamma^{p-1}}q\int W\,d(F^q)=\frac{\gamma^{p-1}}q\int F^q\,d\omega
\]
by Fubini's theorem. We also have $F\ge\Gamma$ and $\Gamma^q=W^{r/p}x^{r/p'}$, since $q/(p-q)=r/p$ and $q(p-1)/(p-q)=r/p'$. Hence $\int F^q\,d\omega\ge\mathcal A^r$. The test function $g=F'$ satisfies $Hg=F$, since $F$ is nondecreasing, locally absolutely continuous and $F(0+)=0$. It gives
\[
C_*\ge\frac{(\int F^q\,d\omega)^{1/q}}{(\int F'^p\,dx)^{1/p}}\ge q^{1/p}\gamma^{-1/p'}\Big(\int F^q\,d\omega\Big)^{1/r}\ge\lambda^\flat_{p,q}\mathcal A .\qedhere
\]
\end{proof}

\begin{proof}[Proof of Theorem~\ref{thm:lowerq}(2)]
The lower bound is (1). Fix $t\in(0,1)$ and $N\ge2$. Let $a_i=t^{-i}$ and $b_i=t^i$ for $1\le i\le N$. Then $(b_i-b_{i+1})a_i=1-t$ for $i<N$, $b_Na_N=1$, and $I=(N-1)(1-t)+1$. A direct computation gives
\[
U_i=\frac{1-t}{1-t^q}\big(1-t^{q(N-i)}\big)+t^{q(N-i)},\qquad V_i=\frac{1-t}{1-t^{p'}}\big(1-t^{(i-1)p'}\big)+t^{(i-1)p'} .
\]
These numbers are bounded. Also $U_i^{r/p}V_i^{r/p'}\to U^{r/p}V^{r/p'}=R_{p,q}(t)^{-r}$ as $\min(i,N-i)\to\infty$, where $U=(1-t)/(1-t^q)$ and $V=(1-t)/(1-t^{p'})$. By \eqref{eq:A-UV} and Lemma~\ref{lem:ab}, $(C_*/\mathcal A)^r=I/\mathcal A^r\to R_{p,q}(t)^r$ as $N\to\infty$. Thus $\lambda_{p,q}\le R_{p,q}(t)$ for every $t$. Finally, by Bernoulli's inequality, $R_{p,q}(t)\le(1-p^{-1}t^q)/(1-t)<1$ for small $t>0$, because $q<1$.

Let $t\uparrow1$ and put $s=1-t$. For $\gamma>0$ we have $(1-t^\gamma)/s=\gamma\big(1-(\gamma-1)s/2+O(s^2)\big)$. Since $(q-1)/p+(p'-1)/p'=q/p$, this gives
\[
\log R_{p,q}(t)=\log\big((p')^{1/p'}q^{1/p}\big)-\frac{q}{2p}\,s+O(s^2).
\]
Thus $R_{p,q}(t)<(p')^{1/p'}q^{1/p}$ for $t$ close to $1$.
\end{proof}

\begin{proof}[Proof of Theorem~\ref{thm:lowerq}(3)]
Let $p=3$, $q=4/5$, $a=(3/1000,1)$ and $b=(1,1/10)$. By \eqref{eq:ab-config} this is the pair $\sigma=\mathrm m$ and $\omega=\frac9{10}(3/1000)^{1/5}\delta_{x_1}+\frac1{10}\delta_{x_2}$, with $x_1=3/1000$ and $x_2=3/1000+\frac{997}{1000}\sqrt{10}$. Here $r=12/11$, $r/p=4/11$ and $r/p'=8/11$. Lemma~\ref{lem:ab} gives $C_*=I^{1/r}$ with $I=\frac9{10}\cdot\frac3{1000}+\frac1{10}=\frac{1027}{10000}$. In \eqref{eq:A-UV} we have $U_2=V_1=1$ and
\[
U_1=\frac9{10}+\frac1{10}\Big(\frac{1000}3\Big)^{1/5},\qquad V_2=\frac{997}{1000}+\frac3{1000\sqrt{10}} .
\]
Hence
\[
\frac{C_*}{\mathcal A}=\Big(\frac ID\Big)^{11/12},\qquad D=\frac{27}{10000}\,U_1^{4/11}+\frac1{10}\,V_2^{8/11} .
\]
For orientation, numerically $C_*/\mathcal A=0.9995283\ldots$ and $\lambda^\sharp_{3,4/5}=0.9996198\ldots$, the latter attained at $t=0.001701\ldots$. Put $c=99955/100000$, so that $c^3=0.9986506\ldots$. We prove $C_*/\mathcal A<c\le\lambda^\sharp_{3,4/5}$. Every step below can be checked by hand.

\emph{Step 1. The bound $C_*/\mathcal A<c$.} We use the elementary inequalities $\ln x\ge2(x-1)/(x+1)$ for $x\ge1$ and $e^y\ge1+y+y^2/2$ for $y\ge0$, and Bernoulli's inequality $x^s\le1+s(x-1)$ for $x>0$ and $0<s<1$. Applied to $1/x$, Bernoulli's inequality gives $x^s\ge1/(1+s(1/x-1))$. Since $3.195^5<1000/3$, we have $U_1>u=1.2195$. So $\ln u\ge0.439/2.2195>0.197792$. Put $y=\frac4{11}\ln u$. Then $y>0.071924$ and $U_1^{4/11}\ge1+y+y^2/2>1.07451$. Since $3.163^2>10$, we have $V_2>v=0.997+0.003/3.163>0.997948$. Hence $1/v-1<0.0020563$ and
\[
V_2^{8/11}>\frac1{1+\frac8{11}\cdot0.0020563}>0.998506 .
\]
These bounds give $D>0.0027\cdot1.07451+0.0998506>0.1027517$ and $I/D<0.1027/0.1027517<0.999497$. Bernoulli's inequality with $s=11/12$ and $x=I/D$ gives
\[
\frac{C_*}{\mathcal A}=\Big(\frac ID\Big)^{11/12}\le1-\frac{11}{12}\Big(1-\frac ID\Big)<1-\frac{11}{12}\cdot0.000503<0.999539<c .
\]

\emph{Step 2. Reduction.} For $s>0$ and $0<t<1$ let $G_s(t)=(1-t^s)/(1-t)$. This is the slope of the chord of $x\mapsto x^s$ between $t$ and $1$. When $0<s<1$, concavity of the power makes $G_s$ decreasing with $G_s(t)>s$, the derivative of $x^s$ at $1$. When $s>1$, convexity makes $G_s$ increasing. By \eqref{eq:lambdas} we have $R_{3,4/5}(t)^3=G_{4/5}(t)\,G_{3/2}(t)^2$. It remains to prove $G_{4/5}(t)G_{3/2}(t)^2>c^3$ for $0<t<1$. Let $T=(3/5)^{10}=0.0060466\ldots$

\emph{Step 3. The range $0<t\le T$.} Let $L(t)=3\ln R_{3,4/5}(t)$. For $0\le x<1$ we have $\ln(1-x)=-\sum_{k\ge1}x^k/k\ge-x-x^2/(2(1-x))$. Also $-\ln(1-t)\ge t$. Hence
\[
L(t)\ge3t-t^{4/5}-\frac{t^{8/5}}{2(1-t^{4/5})}-2t^{3/2}-\frac{t^3}{1-t^{3/2}} .
\]
Let $s=3/5$, so that $T^{1/2}=s^5$, $T^{3/5}=s^6$, $T^{4/5}=s^8$, $T^{3/2}=s^{15}$ and $T^2=s^{20}$. For $0<t\le T$ we have $t^{8/5}\le s^6t$, $t^{3/2}\le s^5t$ and $t^3\le s^{20}t$. The two denominators are at least $1-s^8$ and $1-s^{15}$. We obtain $L(t)\ge\alpha_0 t-t^{4/5}$ with
\[
\alpha_0=3-2s^5-\frac{s^6}{2(1-s^8)}-\frac{s^{20}}{1-s^{15}}>2.82071 ,
\]
since $s^5=0.07776$, $s^6=0.046656$, $s^8=0.01679616$, $s^{15}<0.00048$ and $s^{20}<0.00004$. The function $t\mapsto\alpha_0 t-t^{4/5}$ attains its minimum over $t>0$ at $t^{1/5}=4/(5\alpha_0)$. The minimum equals $-(4/5)^5/(4\alpha_0^4)=-0.08192/\alpha_0^4$. Since $\alpha_0^4>63.30$, we get $L(t)>-0.0012942$. Thus $R_{3,4/5}(t)^3=e^{L(t)}\ge1+L(t)>0.9987058>c^3$ for $0<t\le T$.

\emph{Step 4. The range $T\le t<1$.} Put $z=t^{1/10}\in[3/5,1)$. Then $G_{4/5}(t)=(1-z^8)/(1-z^{10})$ and $G_{3/2}(t)=(1-z^{15})/(1-z^{10})$. By Step 2 the first decreases and the second increases in $z$. Hence, for $z_i\le z\le z_{i+1}$,
\[
G_{4/5}(t)\,G_{3/2}(t)^2\ge G_{4/5}(z_{i+1}^{10})\,G_{3/2}(z_i^{10})^2 .
\]
Table~\ref{tab:cert} lists nodes $z_i$ from $0.60$ to $0.84$ together with truncated values of the right side. Each entry is a lower bound. Every product exceeds $c^3$. For $0.84\le z<1$ Step 2 gives $G_{4/5}(t)>4/5$ and $G_{3/2}(t)\ge G_{3/2}(0.84^{10})>1.12332$. So $G_{4/5}(t)G_{3/2}(t)^2>\frac45\cdot1.12332^2>1.0094>c^3$.

\begin{table}[ht]
\centering
\caption{Lower bounds in Step 4 of the proof of Theorem~\ref{thm:lowerq}(3). All entries are truncated, so each is a lower bound. The product is formed from the two truncated factors and truncated again.}\label{tab:cert}
\begin{tabular}{cccccc}
\hline
$i$ & $z_i$ & $z_{i+1}$ & $G_{4/5}(z_{i+1}^{10})\ge$ & $G_{3/2}(z_i^{10})\ge$ & product $\ge$\\
\hline
1 & $0.60$ & $0.61$ & $0.98787$ & $1.00561$ & $0.99898$\\
2 & $0.61$ & $0.62$ & $0.98644$ & $1.00657$ & $0.99944$\\
3 & $0.62$ & $0.63$ & $0.98488$ & $1.00768$ & $1.00006$\\
4 & $0.63$ & $0.65$ & $0.98134$ & $1.00896$ & $0.99900$\\
5 & $0.65$ & $0.67$ & $0.97720$ & $1.01206$ & $1.00091$\\
6 & $0.67$ & $0.70$ & $0.96974$ & $1.01606$ & $1.00113$\\
7 & $0.70$ & $0.75$ & $0.95358$ & $1.02418$ & $1.00025$\\
8 & $0.75$ & $0.80$ & $0.93233$ & $1.04551$ & $1.01912$\\
9 & $0.80$ & $0.84$ & $0.91155$ & $1.08087$ & $1.06494$\\
\hline
\end{tabular}
\end{table}

Steps 3 and 4 give $R_{3,4/5}(t)>c$ for $0<t<1$. Thus $\lambda^\sharp_{3,4/5}\ge c$. With Step 1 we get $\lambda_{3,4/5}\le C_*/\mathcal A<c\le\lambda^\sharp_{3,4/5}$. This completes the proof.
\end{proof}

\begin{remark}\label{rem:lowerq}
For $1<p\le2$ and $0<q<1$ we conjecture that $\lambda_{p,q}=\lambda^\sharp_{p,q}$. If this holds, then the factor $\lambda^\flat_{p,q}$ is not sharp there, because $\lambda^\flat_{p,q}<\lambda^\sharp_{p,q}$ for all $0<q<1<p$. To see this, put $a=(p-q)/(p-1)$, so that $q<a<p'$ and $q/p+a/p'=1$. For $0<t<1$, H\"older's inequality applied to $1=x^{(q-1)/p}x^{(1-q)/p}$ on $(t,1)$ gives $1-t<\big((1-t^q)/q\big)^{1/p}\big((1-t^a)/a\big)^{1/p'}$, strictly because $q\ne1$. Since $1-t^a<1-t^{p'}$, we get $R_{p,q}(t)>q^{1/p}a^{1/p'}=\lambda^\flat_{p,q}$. The limits of $R_{p,q}(t)$ as $t\downarrow0$ and $t\uparrow1$ are $1$ and $(p')^{1/p'}q^{1/p}$. Both exceed $\lambda^\flat_{p,q}$, the first by the strict weighted inequality of arithmetic and geometric means and the second because $a<p'$. Since $R_{p,q}$ is continuous on $(0,1)$, it follows that $\lambda^\sharp_{p,q}>\lambda^\flat_{p,q}$. The conjecture does not extend to all $p>1$, by part (3). For example, $\lambda^\flat_{2,1/2}=0.8660\ldots$ and $\lambda^\sharp_{2,1/2}=(2\sqrt2-2)^{1/2}=0.9101\ldots$. The proof of part (1) uses a test function that depends only on $W(x)$ and $x$. In contrast, the extremal test functions on chains depend on how the mass to the right of $x$ is distributed. At $(p,q)=(3,4/5)$ a configuration with two atoms gives a smaller ratio than $\lambda^\sharp_{3,4/5}$, the limit of the ratios of geometric chains. Whether configurations with few atoms determine $\lambda_{3,4/5}$ is open. Its value is also unknown.
\end{remark}

\subsection{A discrete criterion by first passage levels}\label{sec:quantile}

The quantile function also gives a discrete form of the criterion. It uses only the first passage times of the clock $x\mapsto S(x)$ through a geometric sequence of levels. It also gives a finite-rank structure that a study of compactness would use, see Remark~\ref{rem:further}.

Fix $\lambda>1$. Let $K=\{k\in\Z:\lambda^k<S(\infty)\}$. For $k\in K$ we define the quantile levels
\begin{equation}\label{eq:levels}
x_k=Q(\lambda^k)=\inf\{x>0: S(x)\ge\lambda^k\}.
\end{equation}
If $K$ has a largest element $k^*$, we put $x_{k^*+1}=\infty$. If $\sigma\ne0$, then $S(\infty)>0$ and $K$ is unbounded below. If $\sigma=0$, then $K$ is empty, $\mathcal D_\lambda=0$ and $C_*=0$. We exclude this case below. We put $I_k=[x_k,x_{k+1})$ for $k\in K$. We call $k\in K$ \emph{effective} if $I_k\neq\emptyset$, that is, if $x_k<x_{k+1}$. The set of effective indices is denoted by $K_e$.

\begin{lemma}\label{lem:levels}
The quantile levels have the following properties.
\begin{enumerate}
\item $0<x_k<\infty$ and $S(x_k)\ge\lambda^k$ for every $k\in K$. The sequence $(x_k)$ is nondecreasing.
\item If $k\in K$ and $x<x_{k+1}$, then $S(x)\le\lambda^{k+1}$. In particular $\lambda^k\le S(x_k)\le\lambda^{k+1}$ for $k\in K_e$.
\item The sets $I_k$ with $k\in K$ are disjoint. Their union is $\{x>0: S(x)>0\}$.
\item For every nonnegative nonincreasing function $\Phi$ on $(0,\infty)$ we have
\begin{equation}\label{eq:sum-int}
\Big(1-\frac1\lambda\Big)\sum_{k\in K}\lambda^k\,\Phi(x_k)\le\int_{(0,\infty)}\Phi\,d\sigma .
\end{equation}
\end{enumerate}
\end{lemma}

\begin{proof}
(1) Since $\lambda^k<S(\infty)$, the set in \eqref{eq:levels} is not empty. Hence $x_k<\infty$. Since $S(x)\to0$ as $x\downarrow0$, we have $x_k>0$. Right-continuity of $S$ gives $S(x_k)\ge\lambda^k$. Monotonicity in $k$ is clear.

(2) If $k+1\in K$ and $x<x_{k+1}$, then $S(x)<\lambda^{k+1}$ by the definition of $x_{k+1}$. If $k=k^*$ is the largest element of $K$, then $\lambda^{k^*+1}\ge S(\infty)\ge S(x)$. If $k\in K_e$, we may take $x=x_k$.

(3) Disjointness follows from (1). Every point of $I_k$ satisfies $S(x)\ge S(x_k)>0$. Conversely, let $S(x)>0$. Since $K$ is unbounded below, the set $\{k\in K:\lambda^k\le S(x)\}$ is not empty. It is bounded above, because $S(x)<\infty$. Let $k$ be its largest element. Then $x\ge x_k$ by \eqref{eq:levels}. If $k+1\in K$, then $\lambda^{k+1}>S(x)$. So $x<x_{k+1}$, since otherwise $S(x)\ge S(x_{k+1})\ge\lambda^{k+1}$. If $k+1\notin K$, then $x_{k+1}=\infty$. In both cases $x\in I_k$.

(4) Let $k\in K$ and $t\in(\lambda^{k-1},\lambda^k]$. Then $t\le\lambda^k\le S(x_k)$, which gives $Q(t)\le x_k$. Since $\Phi$ is nonincreasing, $\Phi(Q(t))\ge\Phi(x_k)$. The intervals $(\lambda^{k-1},\lambda^k]$ with $k\in K$ are disjoint and contained in $(0,S(\infty))$. Hence by \eqref{eq:quantile}
\[
\int_{(0,\infty)}\Phi\,d\sigma=\int_0^{S(\infty)}\Phi(Q(t))\,dt\ge\sum_{k\in K}(\lambda^k-\lambda^{k-1})\,\Phi(x_k),
\]
which is the asserted inequality.
\end{proof}

\begin{remark}
The points $x_k$ are the first passage times of the nondecreasing clock $x\mapsto S(x)$ through the levels $\lambda^k$. In the language of Section~\ref{sec:doob} they are deterministic stopping times for the reverse filtration $(\F_x)$. On each block $I_k$ the $\sigma$-mass to the left of a point stays within a factor $\lambda$. This is the probabilistic counterpart of the discretizing sequences of \cite{GP03}. Atoms of $\sigma$ cause no difficulty. An atom that carries $S$ across several levels gives coinciding points $x_k$ and empty blocks $I_k$.
\end{remark}

We shall also need disjoint sets of large $\sigma$-measure attached to the effective blocks.

\begin{lemma}\label{lem:disjoint}
Let $\lambda>1$ and let $k_1<k_2<\dots<k_n$ be effective indices. Put $k_{-1}=k_0=-\infty$ and $x_{-\infty}=0$. For $1\le m\le n$ let $E_m=(x_{k_{m-2}},x_{k_m}]$. Then the sets $E_m$ with $m$ odd are pairwise disjoint. The sets $E_m$ with $m$ even are pairwise disjoint. Also $E_m\subset(0,x_{k_m}]$ and
\begin{equation}\label{eq:Em}
\sigma(E_m)\ge\Big(1-\frac1\lambda\Big)\lambda^{k_m}\qquad(1\le m\le n).
\end{equation}
\end{lemma}

\begin{proof}
The sets $E_m$ and $E_{m+2}$ are adjacent intervals. Thus sets of the same parity are disjoint. For $m\le2$ we have $\sigma(E_m)=S(x_{k_m})\ge\lambda^{k_m}$. Let $m\ge3$. Since $k_{m-2}$ is effective, Lemma~\ref{lem:levels}(2) gives $S(x_{k_{m-2}})\le\lambda^{k_{m-2}+1}$. Since $k_{m-2}+1\le k_{m-1}\le k_m-1$, we obtain
\[
\sigma(E_m)=S(x_{k_m})-S(x_{k_{m-2}})\ge\lambda^{k_m}-\lambda^{k_m-1}.\qedhere
\]
\end{proof}

For a fixed $\lambda>1$ we write $\omega_k=\omega(I_k)$ for $k\in K$, $\Lambda_k=\lambda^{kr/p'}$ and $\tau=r/p$. Then $\tau>0$, $\tau+1=r/q$ and $\mathcal D_\lambda^r=\sum_{k\in K}\Lambda_k\,\omega_k^{\tau+1}$.

\begin{proposition}[Sufficiency]\label{prop:qlp-suff}
Let $\lambda>1$. Let $\sigma$ and $\omega$ be as in Section~\ref{sec:statement}. Then for every nonnegative Borel function $g$
\begin{equation}\label{eq:qlp-suff}
\Big(\int_{(0,\infty)}\big(H_\sigma g\big)^q\,d\omega\Big)^{1/q}\le \lambda\Big(1-\frac1\lambda\Big)^{-1/p}p'\,\mathcal D_\lambda\,\Big(\int g^p\,d\sigma\Big)^{1/p},
\end{equation}
and
\begin{equation}\label{eq:D-A}
\mathcal D_\lambda\le\Big(\frac rq\Big)^{1/r}\mathcal A .
\end{equation}
\end{proposition}

\begin{proof}
Let $\int g^p\,d\sigma<\infty$. If $S(x)=0$, then $H_\sigma g(x)=0$. So by Lemma~\ref{lem:levels}(3) we can integrate over the blocks $I_k$ only. Let $k\in K$ and $x\in I_k$. Then $x\ge x_k$ and $S(x)>0$. The definition of $A^*$ gives $A(x)\le A^*(x_k)$. Lemma~\ref{lem:levels}(2) gives
\begin{equation}\label{eq:block-bound}
H_\sigma g(x)=S(x)A(x)\le\lambda^{k+1}A^*(x_k).
\end{equation}
Here $A$ and $A^*$ are defined by \eqref{eq:Ax} and \eqref{eq:Astar} with respect to $\sigma$, also when $\sigma$ is infinite. Hence
\[
\int_{(0,\infty)}\big(H_\sigma g\big)^q\,d\omega\le\sum_{k\in K}\lambda^{(k+1)q}A^*(x_k)^q\,\omega_k=\lambda^q\sum_{k\in K}\big(\lambda^kA^*(x_k)^p\big)^{q/p}\,\lambda^{kq/p'}\omega_k .
\]
Here we used $\lambda^{kq}=\lambda^{kq/p}\lambda^{kq/p'}$. We apply H\"older's inequality with the exponents $p/q$ and $(p/q)'=r/q$. This gives
\[
\sum_{k\in K}\big(\lambda^kA^*(x_k)^p\big)^{q/p}\lambda^{kq/p'}\omega_k\le\Big(\sum_{k\in K}\lambda^kA^*(x_k)^p\Big)^{q/p}\Big(\sum_{k\in K}\lambda^{kr/p'}\omega_k^{r/q}\Big)^{q/r}.
\]
The last factor is $\mathcal D_\lambda^q$. For the first factor we apply Lemma~\ref{lem:doob} to the restriction of $\sigma$ to $(0,R]$ and let $R\to\infty$. On $(0,R]$ the maximal averages of the restrictions increase to $A^*$ as $R\to\infty$. Monotone convergence, applied to their $p$-th powers times $\mathbf 1_{(0,R]}$, gives $\int(A^*)^p\,d\sigma\le(p')^p\int g^p\,d\sigma$. The function $\Phi=(A^*)^p$ is nonnegative and nonincreasing. Lemma~\ref{lem:levels}(4) and Lemma~\ref{lem:doob} give
\begin{equation}\label{eq:doob-sum}
\sum_{k\in K}\lambda^kA^*(x_k)^p\le\Big(1-\frac1\lambda\Big)^{-1}\int_{(0,\infty)}(A^*)^p\,d\sigma\le\Big(1-\frac1\lambda\Big)^{-1}(p')^p\int_{(0,\infty)}g^p\,d\sigma .
\end{equation}
Taking $q$-th roots gives \eqref{eq:qlp-suff}.

Let us prove \eqref{eq:D-A}. We may assume $\mathcal A<\infty$. Fix $k\in K$ with $\omega_k>0$. Since $\mathcal A<\infty$ and $S(x_k)\ge\lambda^k>0$, the argument after \eqref{eq:hardy-beta} in Section~\ref{sec:reductions}, which uses only Lemma~\ref{lem:tail} and the monotonicity of $S$, gives $W(x_k)<\infty$. Since $I_k$ contains its left endpoint $x_k$, we obtain $\omega_k\le W(x_k)<\infty$. For $x\in I_k$ we have $W(x)\ge\omega([x,x_{k+1}))$ and $S(x)\ge S(x_k)\ge\lambda^k$. Lemma~\ref{lem:tail} with $I=I_k$, $\gamma=\tau$ and the identity $\tau+1=r/q$ give
\[
\Lambda_k\,\omega_k^{\tau+1}\le(\tau+1)\,\lambda^{kr/p'}\int_{I_k}\omega([x,x_{k+1}))^{\tau}\,d\omega(x)\le\frac rq\int_{I_k}W(x)^{r/p}S(x)^{r/p'}\,d\omega(x).
\]
Summing over $k$ and using the disjointness of the blocks, we obtain \eqref{eq:D-A}.
\end{proof}

\begin{proposition}[Necessity]\label{prop:qlp-nec}
Let $\lambda>1$ and let $\sigma,\omega$ be as in Section~\ref{sec:statement}. If \eqref{eq:main-ineq} holds with a constant $C<\infty$, then
\begin{equation}\label{eq:qlp-nec}
\mathcal D_\lambda\le 2^{1/r}\Big(1-\frac1\lambda\Big)^{-1/p'}C .
\end{equation}
\end{proposition}

\begin{proof}
Let $k_1<\dots<k_n$ be finitely many effective indices. Let $E_m$ be the sets of Lemma~\ref{lem:disjoint}. Let us fix a parity class $P\subset\{1,\dots,n\}$. Let $c_m\ge0$ for $m\in P$. We test \eqref{eq:main-ineq} with
\[
g=\sum_{m\in P}c_m\mathbf 1_{E_m}.
\]
The sets $E_m$ with $m\in P$ are disjoint and have finite $\sigma$-measure. Hence
\[
\int g^p\,d\sigma=\sum_{m\in P}c_m^p\,\sigma(E_m)<\infty .
\]
Let $m\in P$ and $x\in I_{k_m}$. Then $E_m\subset(0,x_{k_m}]\subset(0,x]$. Hence $H_\sigma g(x)\ge c_m\sigma(E_m)$. Summing over the disjoint blocks $I_{k_m}$ gives
\[
\sum_{m\in P}\omega_{k_m}\,c_m^q\,\sigma(E_m)^q\le\int_{(0,\infty)}\big(H_\sigma g\big)^q\,d\omega\le C^q\Big(\sum_{m\in P}c_m^p\,\sigma(E_m)\Big)^{q/p}.
\]
Taking $c_m=\mathbf 1_{\{m=m_0\}}$ shows that $\omega_{k_{m_0}}<\infty$ for every $m_0\in P$. Put $u_m=c_m^p\,\sigma(E_m)$. The last inequality becomes
\[
\sum_{m\in P}\omega_{k_m}\sigma(E_m)^{q/p'}\,u_m^{q/p}\le C^q\Big(\sum_{m\in P}u_m\Big)^{q/p}
\]
for all $u_m\ge0$. By the duality of $\ell^{p/q}$ and $\ell^{r/q}$ on the finite set $P$, this is equivalent to
\[
\sum_{m\in P}\omega_{k_m}^{r/q}\,\sigma(E_m)^{r/p'}\le C^r .
\]
By \eqref{eq:Em} we have $\sigma(E_m)^{r/p'}\ge(1-1/\lambda)^{r/p'}\Lambda_{k_m}$. Adding the two parity classes gives
\[
\sum_{m=1}^n\Lambda_{k_m}\,\omega_{k_m}^{\tau+1}\le 2\Big(1-\frac1\lambda\Big)^{-r/p'}C^r .
\]
Only effective indices contribute to $\mathcal D_\lambda^r$, because $\omega_k=0$ for $k\notin K_e$. Taking the supremum over all finite sets of effective indices gives \eqref{eq:qlp-nec}.
\end{proof}

\begin{proof}[Proof of Theorem~\ref{thm:Dlambda}]
The first display is Propositions~\ref{prop:qlp-suff} and \ref{prop:qlp-nec}. The left inequality in the second display is \eqref{eq:D-A}. The right inequality follows from $C_*\ge c_{p,q}\mathcal A$, which is Theorem~\ref{thm:lowerq}(1) for $q<1$ and Theorem~\ref{thm:lowerB} for $1\le q<p$, and from the right inequality in the first display.
\end{proof}

\begin{remark}
Theorem~\ref{thm:Dlambda} takes a different route from the anti-discretization step of the discretization method \cite{GP03}. The bound of $\mathcal A$ by $\mathcal D_\lambda$ follows from the two characterizations of the same constant $C_*$. No discrete tail estimate is needed for it. The converse bound \eqref{eq:D-A} is a direct tail estimate. The quantity $\mathcal D_\lambda$ is a criterion in its own right. It involves only the $\omega$-masses of the blocks $I_k$. For $1<q<p$ a bounded operator $H_\sigma$ is compact under the hypotheses of \cite[Theorem 3.1(b)]{ML25}, which include an input measure without atoms. In that range the tails of $\mathcal D_\lambda$ do not measure a distance to the compact operators.
\end{remark}

\begin{remark}\label{rem:Dlambda-constants}
The constants in Theorem~\ref{thm:Dlambda} can be optimized in $\lambda$. The function $\lambda\mapsto\lambda(1-1/\lambda)^{-1/p}$ attains its minimum $(p+1)^{1+1/p}/p$ at $\lambda=1+1/p$. Thus
\[
C_*\le\frac{(p+1)^{1+1/p}}{p-1}\,\mathcal D_{1+1/p}.
\]
The ratio of the upper and the lower constant in the first display of Theorem~\ref{thm:Dlambda} is $2^{1/r}p'\lambda^2/(\lambda-1)$, which is smallest at $\lambda=2$, where it equals $2^{2+1/r}p'$. Thus
\[
2^{-1/r}\Big(\frac12\Big)^{1/p'}\mathcal D_2\le C_*\le 2^{1+1/p}p'\,\mathcal D_2 .
\]
The necessity estimate can also be taken free of $\lambda$. By \eqref{eq:D-A} and $C_*\ge c_{p,q}\mathcal A$,
\[
\mathcal D_\lambda\le\min\Big\{2^{1/r}\Big(1-\frac1\lambda\Big)^{-1/p'},\ \frac{(r/q)^{1/r}}{c_{p,q}}\Big\}\,C_*,
\]
with $c_{p,q}$ as in Theorem~\ref{thm:Dlambda}. For $q=1$ this gives $\mathcal D_\lambda\le(p')^{1/p'}C_*$ for every $\lambda>1$.
\end{remark}

\section{Endpoint cases}\label{sec:endpoints}

In this section we treat the remaining exponents for the measure form \eqref{eq:main-ineq}. We interpret the left-hand side of \eqref{eq:main-ineq} as the $\omega$-essential supremum of $H_\sigma g$ when $q=\infty$, and the right-hand side as $C$ times the $\sigma$-essential supremum of $g$ when $p=\infty$.

\begin{proposition}\label{prop:endpoints}
Let $\sigma$ and $\omega$ be as in Section~\ref{sec:statement}. Let $Z=\{x>0:S(x)>0\}$.
\begin{enumerate}
\item Let $1\le p\le\infty$ and $q=\infty$. Then $C_*=B_\infty$, where $B_\infty=\sup\{S(y)^{1/p'}:\,y>0,\ S(y)>0,\ W(y)>0\}$, with $1/p'=0$ if $p=1$ and $\sup\emptyset=0$.
\item Let $p=\infty$ and $0<q<\infty$. Then $C_*=\big(\int S^q\,d\omega\big)^{1/q}$.
\item Let $p=1$ and $0<q<\infty$. Then $C_*=\omega(Z)^{1/q}$.
\end{enumerate}
\end{proposition}

\begin{proof}
(1) Let $S(x)>0$ and let $Y$ be as in Section~\ref{sec:doob} for the restriction of $\sigma$ to $(0,x]$. Jensen's inequality for conditional expectations gives
\[
H_\sigma g(x)=S(x)\,\E\big[g(Y)\mid Y\le x\big]\le S(x)\,\E\big[g(Y)^p\mid Y\le x\big]^{1/p}\le S(x)^{1/p'}\Big(\int g^p\,d\sigma\Big)^{1/p}
\]
for $1<p<\infty$. For $p=1$ and $p=\infty$ the bounds $H_\sigma g(x)\le\|g\|_{L^1(\sigma)}$ and $H_\sigma g(x)\le S(x)\|g\|_{L^\infty(\sigma)}$ are immediate. Also $H_\sigma g=0$ where $S=0$. The function $S$ is nondecreasing. The set $\{x:W(x)=0\}$ is $\omega$-null. Hence the $\omega$-essential supremum of $S^{1/p'}\mathbf 1_{\{S>0\}}$ equals $B_\infty$. This gives the upper bound. For the lower bound fix $y$ with $W(y)>0$ and $S(y)>0$. Put $g=\mathbf 1_{(0,y]}$. Then $H_\sigma g=S(y)$ on $[y,\infty)$, which has positive $\omega$-measure. Also $\|g\|_{L^p(\sigma)}=S(y)^{1/p}$, with $1/p=0$ if $p=\infty$.

(2) We have $H_\sigma g\le S\,\|g\|_{L^\infty(\sigma)}$, with equality for $g=1$.

(3) We have $H_\sigma g\le\|g\|_{L^1(\sigma)}$ everywhere and $H_\sigma g=0$ outside $Z$. This gives $C_*\le\omega(Z)^{1/q}$. If $Z=\emptyset$, then $C_*=0$. Otherwise let $x_0=\inf Z\in[0,\infty)$. If $x_0>0$ and $\sigma(\{x_0\})>0$, the function $g=\mathbf 1_{\{x_0\}}/\sigma(\{x_0\})$ has norm one and $H_\sigma g=1$ on $Z=[x_0,\infty)$. Otherwise $Z=(x_0,\infty)$. For $\varepsilon>0$ with $S(x_0+\varepsilon)>0$ the function $g=\mathbf 1_{(x_0,x_0+\varepsilon]}/S(x_0+\varepsilon)$ has norm one and $H_\sigma g=1$ on $[x_0+\varepsilon,\infty)$. Letting $\varepsilon\downarrow0$ gives $C_*\ge\omega(Z)^{1/q}$.
\end{proof}

\begin{remark}\label{rem:p1}
Part (3) shows that the measure form \eqref{eq:main-ineq} is degenerate for $p=1$. The reason is that a single measure $\sigma$ appears both in the operator and in the norm. The weighted problem with $p=1$ needs two measures in the input. It is treated in Section~\ref{sec:p1}.
\end{remark}

\section{The weighted case \texorpdfstring{$p=1$}{p = 1}}\label{sec:p1}

In this section $\mu$, $\omega$, $v$, $u$ and $U$ are as in Theorem~\ref{thm:p1}. We put $W(x)=\omega([x,\infty))$. For $f\ge0$ we write $H_\mu f(x)=\int_{(0,x]}f\,d\mu$. The maximal inequality of Doob fails for $p=1$. It is replaced by a reduction to point masses. For a finite Borel measure $\nu$ on $(0,\infty)$ we put
\[
\Psi_\nu(x)=\int_{(0,x]}U\,d\nu .
\]
If $\nu$ is the law of a random variable $Y$, then $\Psi_\nu(x)=\E\big[U(Y)\mathbf 1_{\{Y\le x\}}\big]$. We write $Z_\infty=\{x:U(x)=\infty\}$, which is a terminal interval or empty.

\begin{lemma}\label{lem:p1-reduction}
The following statements hold.
\begin{enumerate}
\item $u\le U$ holds $\mu$-almost everywhere.
\item Let $f\ge0$ with $\int fv\,d\mu\le1$. Let $d\nu=fv\,d\mu$. Then $\nu((0,\infty))\le1$ and $H_\mu f(x)\le\Psi_\nu(x)$ for every $x\notin Z_\infty$.
\item Let $y_1,\dots,y_n>0$, $c_1,\dots,c_n\ge0$ and $0\le M_k<U(y_k)$. Then there is a Borel function $f\ge0$ with $\int fv\,d\mu\le\sum_kc_k$ and $H_\mu f(x)\ge\sum_{k:\,y_k\le x}c_kM_k$ for every $x>0$.
\end{enumerate}
\end{lemma}

\begin{proof}
(1) Suppose the contrary. Then there is a rational $t\ge0$ such that the set $N_t=\{y:u(y)>t\ge U(y)\}$ has positive $\mu$-measure. The set $\{U\le t\}$ is an initial interval. The set $N_t$ is contained in it. If $y\in\{U\le t\}$, then $\mu(\{u>t\}\cap(0,y])=0$ by \eqref{eq:U}. The initial interval $\{U\le t\}$ is a countable increasing union of intervals $(0,y_n]$ with $U(y_n)\le t$. Hence $\mu(N_t)\le\mu(\{u>t\}\cap\{U\le t\})=0$, a contradiction.

(2) Let $x\notin Z_\infty$. Then $U(x)<\infty$. So $v>0$ holds $\mu$-almost everywhere on $(0,x]$. Since $\int fv\,d\mu<\infty$, $f=0$ holds $\mu$-almost everywhere on $\{v=\infty\}$. It follows that $f=(fv)\,u\le (fv)\,U$ holds $\mu$-almost everywhere on $(0,x]$, by (1). Integrating over $(0,x]$ gives $H_\mu f(x)\le\Psi_\nu(x)$.

(3) We may assume $M_k>0$ for all $k$. Since $U(y_k)>M_k$, the set $\{u>M_k\}\cap(0,y_k]$ has positive $\mu$-measure. Since $\mu$ is $\sigma$-finite, it contains a Borel set $G_k$ with $0<\mu(G_k)<\infty$. Put $f=\sum_kc_kM_k\mathbf 1_{G_k}/\mu(G_k)$. Since $v<1/M_k$ on $G_k$, we get $\int fv\,d\mu\le\sum_kc_k$. If $y_k\le x$, then $G_k\subset(0,x]$. Hence $H_\mu f(x)\ge\sum_{k:\,y_k\le x}c_kM_k$.
\end{proof}

If $\omega(Z_\infty)>0$, then $W(y)>0$ for some $y\in Z_\infty$. Also $W>0$ holds $\omega$-almost everywhere on $Z_\infty$. In this case $B_1$, $\sup\{U(y):W(y)>0\}$ and $\mathcal A_1$ are all infinite. Lemma~\ref{lem:p1-reduction}(3) with one point $y\in Z_\infty$ with $W(y)>0$ and arbitrarily large $M$ shows that $C_{1,*}=\infty$ as well. So we assume $\omega(Z_\infty)=0$ in the proofs below. A function $f$ with $\int fv\,d\mu=\infty$ satisfies \eqref{eq:p1-ineq} for every $C>0$. Such a function does not affect the infimum $C_{1,*}$. By homogeneity it suffices to consider $f$ with $\int fv\,d\mu\le1$. Then Lemma~\ref{lem:p1-reduction}(2) gives $H_\mu f\le\Psi_\nu$ $\omega$-almost everywhere.

\begin{proof}[Proof of Theorem~\ref{thm:p1}(1) and (2)]
Let $1\le q<\infty$. By Minkowski's integral inequality
\[
\Big(\int\Psi_\nu^q\,d\omega\Big)^{1/q}=\Big\|\int U(y)\mathbf 1_{[y,\infty)}\,d\nu(y)\Big\|_{L^q(\omega)}\le\int U(y)\,W(y)^{1/q}\,d\nu(y)\le B_1 .
\]
Here the points $y\in Z_\infty$ may be omitted, because $\Psi_\nu$ is only needed outside the $\omega$-null set $Z_\infty$. Thus $C_{1,*}\le B_1$. Conversely, let $y>0$ and $M<U(y)$. Lemma~\ref{lem:p1-reduction}(3) with $n=1$ and $c_1=1$ gives $f$ with $\int fv\,d\mu\le1$ and $H_\mu f\ge M$ on $[y,\infty)$. This test gives $C_{1,*}\ge MW(y)^{1/q}$. Letting $M\uparrow U(y)$ gives $C_{1,*}\ge B_1$.

Let $q=\infty$. For $x\notin Z_\infty$ we have $\Psi_\nu(x)\le U(x)\,\nu((0,x])\le U(x)$. The set $\{W=0\}$ is $\omega$-null. The function $U$ is nondecreasing. Hence the $\omega$-essential supremum of $U$ is $\sup\{U(y):W(y)>0\}$, which gives the upper bound. The lower bound follows from Lemma~\ref{lem:p1-reduction}(3) as before, because $\omega([y,\infty))>0$.
\end{proof}

For the case $0<q<1$ we use the level sets of the running essential supremum. For $s>0$ let
\[
E_s=\{x>0:U(x)^r>s\},\qquad V(s)=\omega(E_s),\qquad \tilde{\mathcal A}_1=\Big(\int_0^\infty V(s)^{r+1}\,ds\Big)^{1/r}.
\]
Each $E_s$ is a terminal interval, because $U$ is nondecreasing. Its left endpoint is the first passage time of $U$ above the level $s^{1/r}$. If $x\in E_s$, then $[x,\infty)\subset E_s$. Thus $W(x)=\omega(E_s\cap[x,\infty))$.

\begin{lemma}\label{lem:p1-layers}
Let $0<q<1$. Then $\mathcal A_1\le\tilde{\mathcal A}_1\le(r+1)^{1/r}\mathcal A_1$.
\end{lemma}

\begin{proof}
We assume $\omega(Z_\infty)=0$, as in the rest of this section. By Fubini's theorem
\[
\mathcal A_1^r=\int_0^\infty\int_{E_s}W^r\,d\omega\,ds ,
\]
because $x\in E_s$ if and only if $s<U(x)^r$. Since $W\le V(s)$ on $E_s$, we have $\int_{E_s}W^r\,d\omega\le V(s)^{r+1}$, which gives the first inequality. Lemma~\ref{lem:tail}, applied to finite measures $\nu\le\omega|_{E_s}$, gives $\nu(E_s)^{r+1}\le(r+1)\int_{E_s}W^r\,d\omega$. Since $\omega$ is $\sigma$-finite, the supremum over such $\nu$ gives $V(s)^{r+1}\le(r+1)\int_{E_s}W^r\,d\omega$. Integrating over $s$ gives the second inequality.
\end{proof}

\begin{proof}[Proof of Theorem~\ref{thm:p1}(3), sufficiency]
We show $C_{1,*}\le\tilde{\mathcal A}_1$. Then Lemma~\ref{lem:p1-layers} gives the upper bound in \eqref{eq:constants-p1}. Let $\tilde{\mathcal A}_1<\infty$. Put
\[
\Phi(x)=\int_0^{U(x)^r}V(s)^r\,ds ,
\]
a nondecreasing function. If $s<U(x)^r$, then $[x,\infty)\subset E_s$, which gives $V(s)\ge W(x)$. Integrating over $s$ gives
\begin{equation}\label{eq:Phi-lower}
\Phi(x)\ge U(x)^rW(x)^r .
\end{equation}
By Fubini's theorem $\int\Phi\,d\omega=\int_0^\infty V(s)^r\omega(E_s)\,ds=\tilde{\mathcal A}_1^r$. So $\Phi<\infty$ holds $\omega$-almost everywhere.

Let $\nu$ be as in Lemma~\ref{lem:p1-reduction}(2) and $F=\Psi_\nu$. Let $x\notin Z_\infty$ with $\Phi(x)=0$. Then $U(x)=0$ or $W(x)=0$ by \eqref{eq:Phi-lower}. In the first case $F(x)=0$. The set $\{W=0\}$ is $\omega$-null. H\"older's inequality with the exponents $1/q$ and $1/(1-q)$, applied on $\{\Phi>0\}$, gives
\[
\int F^q\,d\omega\le\Big(\int_{\{\Phi>0\}}F\,\Phi^{-1/r}\,d\omega\Big)^q\Big(\int\Phi\,d\omega\Big)^{1-q},
\]
since $(1-q)/q=1/r$. Let $y>0$ with $U(y)>0$ and $W(y)>0$. Since $\Phi$ is nondecreasing, \eqref{eq:Phi-lower} gives
\[
U(y)\int_{[y,\infty)\cap\{\Phi>0\}}\Phi^{-1/r}\,d\omega\le U(y)\,W(y)\,\Phi(y)^{-1/r}\le1 .
\]
If $U(y)=0$ or $W(y)=0$, the left-hand side vanishes. Hence, by Fubini's theorem,
\[
\int_{\{\Phi>0\}}F\,\Phi^{-1/r}\,d\omega=\int U(y)\int_{[y,\infty)\cap\{\Phi>0\}}\Phi^{-1/r}\,d\omega\,d\nu(y)\le\nu((0,\infty))\le1 .
\]
Here the points $y\in Z_\infty$ may be omitted, as in the proof of part (1). We conclude that $\int F^q\,d\omega\le\tilde{\mathcal A}_1^{r(1-q)}=\tilde{\mathcal A}_1^q$. Thus $C_{1,*}\le\tilde{\mathcal A}_1$.
\end{proof}

\begin{proof}[Proof of Theorem~\ref{thm:p1}(3), necessity]
Let \eqref{eq:p1-ineq} hold with a constant $C<\infty$. Then $\omega(Z_\infty)=0$, as shown above. We first show that
\begin{equation}\label{eq:p1-V}
V(s)\le C^qs^{-q/r}<\infty\qquad\text{for every }s>0 .
\end{equation}
We may assume $E_s\neq\emptyset$. Let $y\in E_s$. Then $M=s^{1/r}<U(y)$. Lemma~\ref{lem:p1-reduction}(3) with one point gives $f$ with $\int fv\,d\mu\le1$ and $H_\mu f\ge M$ on $[y,\infty)$. Hence $M^qW(y)\le C^q$, that is, $W(y)\le C^qs^{-q/r}$. We have $W(y)=\omega(E_s\cap[y,\infty))$ for $y\in E_s$. These sets increase to $E_s$ as $y$ decreases to the left endpoint of $E_s$ through $E_s$. Monotone convergence now gives \eqref{eq:p1-V}. Next we show the following claim. For every Borel function $\rho\ge0$ on $(0,\infty)$ with $\int\rho\,ds\le1$,
\begin{equation}\label{eq:p1-dual}
\int\Big(\int_0^{U(x)^r}s^{1/r}\rho(s)\,ds\Big)^qd\omega(x)\le C^q .
\end{equation}
Let $0<s_1<\dots<s_n$ and $c_k\ge0$ with $\sum_kc_k\le1$. Let $\delta>0$. If $E_{s_k}=\emptyset$, then no $x$ satisfies $s_k<U(x)^r$. So the index $k$ does not occur in the inequality below. We discard it. For the remaining indices choose $y_k\in E_{s_k}$ with $\omega(E_{s_k}\cap(0,y_k))<\delta$. This is possible by \eqref{eq:p1-V}. Indeed, $V(s_k)<\infty$. Also $E_{s_k}$ is a nonempty terminal interval. If $E_{s_k}$ contains its left endpoint, we take $y_k$ to be this endpoint. Otherwise the sets $E_{s_k}\cap(0,y)$ decrease to the empty set as $y$ decreases to the left endpoint through $E_{s_k}$. Their $\omega$-measures tend to $0$ by continuity from above. Lemma~\ref{lem:p1-reduction}(3) with $M_k=s_k^{1/r}<U(y_k)$ gives a function $f$ with $\int fv\,d\mu\le1$ and $H_\mu f(x)\ge\sum_{k:\,y_k\le x}c_ks_k^{1/r}$. Let $x\in E_{s_k}$. If $x\ge y_k$, then the term $c_ks_k^{1/r}$ appears in the lower bound for $H_\mu f(x)$. Otherwise $x$ lies in $E_{s_k}\cap(0,y_k)$. These exceptional sets have total $\omega$-measure less than $n\delta$. Since $t\mapsto t^q$ is subadditive for $q<1$, we obtain
\[
\int\Big(\sum_{k:\,s_k<U(x)^r}c_ks_k^{1/r}\Big)^qd\omega(x)\le C^q+n\delta\Big(\sum_kc_ks_k^{1/r}\Big)^q ,
\]
and letting $\delta\to0$ removes the last term. Now let $\rho$ be given and let $m\in\N$. Let $J_1,\dots,J_n$ be the intervals $(j2^{-m},(j+1)2^{-m}]$ contained in $(2^{-m},m]$. Let $s_k$ be the right endpoint of $J_k$ and $c_k=\int_{J_k}\rho\,ds$. Then $c_ks_k^{1/r}\ge\int_{J_k}s^{1/r}\rho(s)\,ds$. Also $s_k<T$ for every $J_k\subset(0,T)$. The intervals $J_k$ contained in $(0,T)$ cover $(2^{-m},T-2^{1-m}]\cap(0,m]$. Thus the discrete inequality bounds the integral of the $q$-th power of $\int_{(2^{-m},T-2^{1-m}]\cap(0,m]}s^{1/r}\rho\,ds$ with $T=U(x)^r$. Letting $m\to\infty$, Fatou's lemma gives \eqref{eq:p1-dual}.

Next we may assume that $\omega$ is finite and that $U\le n$ holds $\omega$-almost everywhere for some $n$. Indeed, let $\omega_n$ be the restriction of $\omega$ to $\{x:1/n<U(x)\le n\}$. It is finite by \eqref{eq:p1-V}, because $\{U>1/n\}=E_{n^{-r}}$. It satisfies \eqref{eq:p1-ineq} with the same constant. If $U(x)>0$, then $[x,\infty)\subset\{U\ge U(x)\}$. Thus $\omega_n([x,\infty))\uparrow W(x)$, since $\omega(Z_\infty)=0$. Since the integrand of $\mathcal A_1$ vanishes where $U=0$, monotone convergence shows that $\mathcal A_1$ for $\omega_n$ increases to $\mathcal A_1$. Then $0\le\tilde{\mathcal A}_1<\infty$. We may assume $\mathcal A_1>0$, so that $\tilde{\mathcal A}_1>0$ by Lemma~\ref{lem:p1-layers}. Put $\rho=V^{r+1}/\tilde{\mathcal A}_1^r$. Let $x$ satisfy $T=U(x)^r>0$. For $s<T$ we have $V(s)\ge\omega(\{U^r\ge T\})\ge W(x)$. Since $\int_0^Ts^{1/r}\,ds=qT^{1/q}$, we obtain
\[
\int_0^Ts^{1/r}\rho(s)\,ds\ge q\,\frac{W(x)^{r+1}U(x)^{r/q}}{\tilde{\mathcal A}_1^r}.
\]
Since $(r+1)q=r$, the claim \eqref{eq:p1-dual} gives $q^q\mathcal A_1^r\le C^q\tilde{\mathcal A}_1^{rq}\le C^q(r+1)^q\mathcal A_1^{rq}$, by Lemma~\ref{lem:p1-layers}. Since $r(1-q)=q$ and $r+1=1/(1-q)$, this gives $q(1-q)\mathcal A_1\le C$. Taking the infimum over $C$ gives the lower bound in \eqref{eq:constants-p1}.
\end{proof}

\begin{remark}
For Lebesgue measure $\mu$ and $d\omega=w\,dx$, Theorem~\ref{thm:p1}(1) is the classical formula $C_{1,*}=\sup_{t>0}\big(\int_t^\infty w\big)^{1/q}\operatorname{ess\,sup}_{(0,t)}v^{-1}$, see \cite{SS96} and the references there. Part (3) recovers the condition of Sinnamon and Stepanov \cite[Theorem 3.3]{SS96}, namely
\[
\int_0^\infty\Big(\operatorname*{ess\,inf}_{(0,x)}v\Big)^{-q/(1-q)}\Big(\int_x^\infty w\Big)^{q/(1-q)}w(x)\,dx<\infty .
\]
For the counting measure $\mu$ on $\N$ we have $U(n)=\max_{j\le n}v_j^{-1}$. Theorem~\ref{thm:p1} characterizes the discrete inequality $\big(\sum_nw_n(\sum_{j\le n}a_j)^q\big)^{1/q}\le C\sum_ja_jv_j$. For general measures, and indeed for abstract Hardy operators, the case $p=1$ was treated by Santacruz Hidalgo \cite{SH24}. He obtained the exact constant for $q\ge1$. Hence Theorem~\ref{thm:p1}(1) is not new in that generality. For $q<1$ his characterization is up to equivalence. It extends \cite[Theorem 3.1]{JSU08}, which treats three $\sigma$-finite Borel measures on the half-line. To our knowledge, the explicit constants in Theorem~\ref{thm:p1}(3) are new for general measures. The results of \cite{JSU08,SH24} for $q<1$ are equivalences with unspecified constants. Our argument reaches the constants of \cite[Theorem 3.3]{SS96} directly on the first passage sets $E_s$. That argument does not use the level function or Theorem 3.2 of \cite{SS96}. The intermediate quantity $\tilde{\mathcal A}_1$ of Lemma~\ref{lem:p1-layers} also satisfies $C_{1,*}\le\tilde{\mathcal A}_1$. Also $\tilde{\mathcal A}_1\le(1-q)^{-(1-q)/q}\mathcal A_1$ by Lemma~\ref{lem:p1-layers}.
\end{remark}

\section{Weighted, discrete and dual forms}\label{sec:corollaries}

For $1<p<\infty$ we use the following forms of \eqref{eq:main-ineq}. In part (2) we let $(s_n)_{n\ge1}$ and $(w_n)_{n\ge1}$ be nonnegative sequences, $S_n=\sum_{j\le n}s_j$, $W_n=\sum_{m\ge n}w_m$, $B_d=\sup_nW_n^{1/q}S_n^{1/p'}$ and $\mathcal A_d=\big(\sum_nW_n^{r/p}S_n^{r/p'}w_n\big)^{1/r}$. In part (3) we use the Copson operator $H^*_\omega h(y)=\int_{[y,\infty)}h\,d\omega$.

\begin{proposition}\label{prop:forms}
Let $1<p<\infty$ and $0<q<\infty$. Each of the following inequalities holds with a finite constant if and only if $B<\infty$ in the case $p\le q$, and if and only if $\mathcal A<\infty$ in the case $q<p$. Its best constant satisfies \eqref{eq:constants} and \eqref{eq:constants-qlp}, respectively, together with the bounds of Theorem~\ref{thm:sharp} for $0<q<p$, of Theorem~\ref{thm:lowerq}(1) for $0<q<1$ and of Theorem~\ref{thm:lowerB} for $1\le q<p$.
\begin{enumerate}
\item The weighted inequality \eqref{eq:hardy-weighted}, where $0<v<\infty$ and $0\le w<\infty$ almost everywhere, $v^{1-p'}\in L^1_{\rm loc}([0,\infty))$, and $B$, $\mathcal A$ are given by \eqref{eq:B-weighted} and \eqref{eq:A-weighted}.
\item The discrete inequality
\begin{equation}\label{eq:discrete}
\Big(\sum_{n\ge1}w_n\Big(\sum_{j\le n}s_jb_j\Big)^q\Big)^{1/q}\le C\Big(\sum_{j\ge1}s_jb_j^p\Big)^{1/p},
\end{equation}
with $B_d$ and $\mathcal A_d$ in place of $B$ and $\mathcal A$.
\item If $q>1$, the dual inequality $\big(\int(H^*_\omega h)^{p'}d\sigma\big)^{1/p'}\le C\big(\int h^{q'}d\omega\big)^{1/q'}$. Its best constant equals $C_*$.
\end{enumerate}
\end{proposition}

\begin{proof}
(1) Put $d\sigma=v^{1-p'}dx$, $d\omega=w\,dx$ and $g=fv^{p'-1}$. Then $H_\sigma g(x)=\int_0^xf$ and $g^pv^{1-p'}=f^pv$, since $(p-1)(p'-1)=1$. Every $g\ge0$ arises in this way. Since $\omega$ has no atoms, \eqref{eq:B} and \eqref{eq:A} become \eqref{eq:B-weighted} and \eqref{eq:A-weighted}. For the same reason the lower bound $C_*\ge\mathcal A$ improves to $C_*\ge(r/q)^{1/r}\mathcal A$ for $1\le q<p$, see the discussion after the proof of Theorem~\ref{thm:lowerB}. For $0<q<1$ the larger lower factor $(p')^{1/p'}q^{1/p}$ of \cite[Theorem 1(3)]{Si91} applies.

(2) Put $\sigma=\sum_js_j\delta_j$, $\omega=\sum_nw_n\delta_n$ and $g(j)=b_j$. For $x\in[n,n+1)$ we have $S(x)=S_n$ and $W(x)\le W_n$, with equality at $x=n$. Hence $B=B_d$ and $\mathcal A=\mathcal A_d$.

(3) By Tonelli's theorem $\int(H_\sigma g)h\,d\omega=\int g\,(H^*_\omega h)\,d\sigma$ for $g,h\ge0$. Since $\sigma$ and $\omega$ are $\sigma$-finite and $1<p,q<\infty$, the norms of $L^p(\sigma)$ and $L^q(\omega)$ are attained by duality against nonnegative functions. So the two best constants coincide, including the case where both are infinite.

In each case Theorems~\ref{thm:main}, \ref{thm:qlp}, \ref{thm:sharp}, \ref{thm:lowerB} and \ref{thm:lowerq}(1) complete the proof.
\end{proof}

With $a_j=s_jb_j$, so that $a_j=0$ when $s_j=0$, and with the convention $0\cdot\infty=0$, the inequality \eqref{eq:discrete} becomes the usual weighted discrete Hardy inequality with the weight $s_j^{1-p}$ on the right. For $q<p$ it was studied by Bennett \cite{Be91}. Mixed cases, in which one measure is discrete and the other is absolutely continuous, are covered in the same way. Compare \cite{LM20} for $p\le q$.

\section{Concluding remarks}\label{sec:remarks}

\begin{remark}[Probabilistic form]
Let $\sigma$ and $\omega$ be finite and nonzero. Let $Y$ and $X$ be independent random variables with laws $\sigma/S(\infty)$ and $\omega/\omega((0,\infty))$. We put $A(x)=0$ when $S(x)=0$. By independence
\[
\frac{\E\big[g(Y)\mathbf 1_{\{Y\le X\}}\mid X\big]}{\mathbb P(Y\le X\mid X)}=A(X)
\]
almost surely on the event $\{S(X)>0\}$. Since $H_\sigma g=S\,A$, the inequality \eqref{eq:main-ineq} is equivalent to
\[
\Big(\E\big[S(X)^q A(X)^q\big]\Big)^{1/q}\le C\,\omega((0,\infty))^{-1/q}S(\infty)^{1/p}\Big(\E\big[g(Y)^p\big]\Big)^{1/p}.
\]
This is the form of the weighted inequality studied by Klaassen and Wellner \cite{KW21} for $p\le q$. The reductions of Section~\ref{sec:tools} stay inside this language. The coupling realizes $Y$ as $Q(T)$ with $T$ uniformly distributed. The passage from functions of $T$ to functions of $Y$ is a conditional expectation. In this sense the coupling is the probabilistic reading of the normal form of \cite{Si22}.
\end{remark}

\begin{remark}[Empirical distributions]
None of the factors in our two-sided estimates depends on the measures. Let $\sigma_n$ be the empirical distribution of a sample $y_1,\dots,y_n>0$, ties allowed. Then the least constant in $\big(\int R_n^q\,d\omega\big)^{1/q}\le C\big(n^{-1}\sum_ia_i^p\big)^{1/p}$, where $R_n(x)=n^{-1}\sum_{y_i\le x}a_i$ is the cumulative sum of nonnegative marks, is the constant $C_*$ of the pair $(\sigma_n,\omega)$. Indeed, replacing the marks at tied sample points by their average leaves $R_n$ unchanged and, by Jensen's inequality, does not increase $n^{-1}\sum_ia_i^p$. This averaging is a conditional expectation, as in Lemma~\ref{lem:coupling}. Thus the two least constants are equal. As a result, all bounds of this paper hold for every sample, with $B$ and $\mathcal A$ computed from the order statistics.
\end{remark}

\begin{remark}[Further directions]\label{rem:further}
The blocks $I_k$ of Section~\ref{sec:quantile} also give operators of finite rank. For $n\in\N$ let $\mathcal P_ng=\sum_k\big(\sigma(I_k)^{-1}\int_{I_k}g\,d\sigma\big)\mathbf 1_{I_k}$, where the sum runs over the $k\in K$ with $|k|\le n$ and $0<\sigma(I_k)<\infty$. If $H_\sigma$ is bounded, then $H_\sigma\circ\mathcal P_n$ is a bounded operator of rank at most $2n+1$, because $\mathcal P_n$ is a contraction of $L^p(\sigma)$. For fixed $\lambda$ it need not converge to $H_\sigma$ even strongly (take $\sigma=\mathrm m$, $\lambda=2$, $g=\mathbf 1_{(1,3/2)}$, $d\omega=\mathbf 1_{(1,2)}dx$). Hence blocks must be refined. For $1<p\le q$ compactness of Hardy operators is characterized by Sinnamon \cite[Theorems 3.23 and 5.6]{Si22}, through finite rank approximation and the normal form parameter $b$. The limits of $W(y)^{1/q}S(y)^{1/p'}$ where $S(y)\downarrow0$ or $S(y)\uparrow S(\infty)$ are not the right quantity. For $d\sigma=\mathbf 1_{(0,1)}dx$ and $\omega=\delta_1$ this product tends to $1$ as $y\uparrow1$. However, $H_\sigma$ has rank one. A probabilistic proof of that criterion and a formula for the essential norm through the blocks $I_k$ are not attempted here. For $1<q<p$ boundedness implies compactness under the hypotheses of \cite[Theorem 3.1(b)]{ML25}.

For $q<p$ one may also ask for a Bellman function in the style of Section~\ref{sec:pleq-sharp}. Let $\mathbb B(x,y,t)$ be the supremum of $\int_{(0,t]}(Hh)^q\,d\hat\omega$ over all $h\ge0$ on $(0,t]$ with $Hh(t)=tx$ and $\int_0^th^p\,d\mathrm m=ty$, after the coupling of Lemma~\ref{lem:coupling}. We expect, but do not prove, that the splicing argument of \cite[Lemma 2.2]{Os15} characterizes \eqref{eq:main-ineq} by nonnegative functions on the attainable states. These are the states with $y\ge x^p$ in which $x=0$ forces $y=0$. The functions should have the domination $C^q(ty)^{q/p}$ and the pathwise monotonicity of \cite[Section 2]{Os15}. Every such function dominates $\mathbb B$. If $\mathbb B$ itself belongs to this class, it is its least member. Lemma~\ref{lem:EL} is a related first-variation calculation for the operator norm. The best lower factor for $0<q<1$ is a further open problem. Sharp constants for bilinear operators $H_{\sigma_1}g_1\cdot H_{\sigma_2}g_2$ with general measures are open as well. For weights, characterizations of the bilinear inequality on arbitrary intervals are given in \cite{Ka26a}. Compactness and essential norms of weighted bilinear Hardy operators are studied in \cite{Ka26b} for $1<p_1,p_2<\infty$ and in \cite{Ka26c} when at least one input exponent is $1$ or $\infty$. Compactness of weighted bilinear Hardy--Steklov operators is studied in \cite{Ka26d}. None of these questions is used in the paper.
\end{remark}

\section*{Declarations}
\noindent\textbf{Funding.} The author received no funding for this work.\\
\textbf{Competing interests.} The author has no competing interests to declare.\\
\textbf{Data availability.} No external data were used. No proof in this paper depends on a computer calculation. A script that rechecks the numerical entries in the proof of Theorem~\ref{thm:lowerq}(3) in exact rational arithmetic is provided as the supplementary file \texttt{check\_thm2\_7\_3.py}.

\section*{Declaration of AI}
During the preparation of this manuscript, the author used Hark and OpenAI's ChatGPT Work to assist with language editing, to improve the clarity and organization of the exposition, to prepare the \LaTeX{} source and to provide supportive checks of certain mathematical calculations and arguments. All mathematical statements, proofs, derivations, and conclusions were written and verified by the author, who takes full responsibility for the content.

\end{document}